\documentclass[12pt,letterpaper]{amsart}

\usepackage{amssymb,latexsym, amsmath, amsxtra,psfig}
\usepackage[dvips]{graphics}

\theoremstyle{plain}
\newtheorem{theorem}{Theorem}[section]
\newtheorem{corollary}[theorem]{Corollary}
\newtheorem{conjecture}[theorem]{Conjecture}
\newtheorem{lemma}[theorem]{Lemma}

\theoremstyle{definition}

\theoremstyle{remark}

\numberwithin{equation}{subsection}
\numberwithin{theorem}{subsection}
\numberwithin{table}{subsection}

\def\({\left(}
\def\){\right)}

\newcommand{\ontop}[2]{\genfrac{}{}{0pt}{}{#1}{#2}}

\renewcommand{\mod}{\bmod}

\def\sumhar{\operatornamewithlimits{\sum\nolimits^h}}
\def\sumstar{\operatornamewithlimits{\sum\nolimits^*}}

\def\dfn{\it} 
\def\sgn{\hbox{sgn}}

\def\myitem#1{\noindent{\it #1 }}

\begin{document}

\title{
Integral moments of $L$-functions
 }

\abstract
We give a new heuristic for all of the main terms in the
integral moments of various families of primitive $L$-functions.  The
results agree with previous conjectures for the leading
order terms.  Our conjectures also have an almost identical form to
exact expressions for the corresponding moments of the
characteristic polynomials of either unitary, orthogonal,
or symplectic matrices, where the moments are defined by
the appropriate group averages.  This lends support to the
idea that arithmetical $L$-functions have a spectral
interpretation, and that their value distributions can be
modeled using Random Matrix Theory.  Numerical examples
show good agreement with our conjectures.
\endabstract
\author{
J. B. Conrey, D. W. Farmer, J. P. Keating,
M. O. Rubinstein, N. C. Snaith
 }

\address{
{\parskip 0pt
American Institute of Mathematics\endgraf
360 Portage Ave.\endgraf
Palo Alto, CA 94306\endgraf
\null
Department of Mathematics\endgraf
Oklahoma State University\endgraf
Stillwater, OK 74078-0613\endgraf
\null
School of Mathematics\endgraf
University of Bristol \endgraf
Clifton, Bristol\endgraf
BS8 1TW\endgraf
United Kingdom\endgraf
\null
Pure Mathematics\endgraf
University of Waterloo \endgraf
Waterloo, Ontario\endgraf
N2L 3G1\endgraf
Canada\endgraf
}
  }

\thanks{
Research partially supported by the American Institute of
Mathematics and a Focused Research Group grant from the National
Science Foundation. The last author was also supported by a Royal
Society Dorothy Hodgkin Fellowship.
 }

\maketitle

 \tableofcontents

\newpage

\parskip 4pt

\section{  Introduction and Statement of Results
}\label{sec:intro}

Random Matrix Theory (RMT) has recently become a fundamental tool
for understanding $L$-functions. Montgomery \cite{Mon} showed that
the two-point correlations between the non-trivial zeros of the
Riemann $\zeta$-function, on the scale of the mean zero spacing,
are similar to the corresponding correlations between the
eigenvalues of random unitary matrices in the limit of large
matrix size \cite{Meh} and conjectured that these correlations
are, in fact, identical to each other.  There is extensive
numerical evidence \cite{Odl} in support of this conjecture.
Rudnick and Sarnak \cite{RS} extended Montgomery's analysis to all
$n$-point correlations, and to the zeros of other principal
$L$-functions.  Katz and Sarnak \cite{KSa} introduced the idea of
studying zero distributions within families of $L$-functions (see
also \cite{OS,Rub}) and have conjectured that these coincide with
the eigenvalue distributions of the classical compact groups. In
this context symmetries of an $L$-function family determine the
associated classical group. We shall here be concerned with the
distribution of values taken by $L$-functions, either individually
(i.e.~along the appropriate critical line), or with respect to
averages over families.  Specifically, we shall calculate the
integral moments of these distributions.

Keating and Snaith \cite{KS1} suggested that the value
distribution of the Riemann $\zeta$-function (or any other
principal $L$-function) on its critical line is related to
that of the characteristic polynomials of random unitary
matrices.  This led them to a general conjecture for the
leading-order asymptotics of the moments of this
distribution in the limit of large averaging range.   Their
conjecture agrees with a result of Hardy and Littlewood~\cite{HL} for
the second moment and a result of Ingham~\cite{I} for the fourth
moment (see, for example \cite{T}). It also agrees with
conjectures, based on number-theoretical calculations, of
Conrey and Ghosh \cite{CG2} and Conrey and Gonek \cite{CGo} for
the sixth and eighth moments.  General conjectures for the
leading-order asymptotics of the moments of $L$-functions
within families, based on random-matrix calculations for
the characteristic polynomials of matrices from the
orthogonal and unitary-symplectic groups, were developed by
Conrey and Farmer \cite{CF} and Keating and Snaith~\cite{KS2}.
These are also in agreement with what is known, and with
previous conjectures.  

Our purpose here is, for the integral moments of a family of
primitive $L$-functions, to go beyond
the leading order asymptotics previously investigated: we
give conjectures for the full main terms.
We propose a refined definition of ``conductor'' of an $L$-function,
which to leading order is the (logarithm of) the ``usual'' conductor.
We find that often, but not always, the mean values can be
expressed as polynomials in the conductor.
Importantly, our conjectures show a striking formal
similarity with analogous expressions for the
characteristic polynomials of random matrices.  This
provides a strong measure of the depth of the connection
between $L$-functions and RMT.  We also perform numerical
calculations which show very good agreement with our
conjectures.  Non-primitive families can also
be handled by our methods, but we do not treat those here.

The conjectures we develop here can also be obtained by
techniques of multiple Dirichlet
series, as described by Diaconu, Goldfeld,
and Hoffstein~\cite{DGH}.  In their formulation, one
considers Dirichlet series in several complex variables.
The mean values we conjecture would then follow from
a plausible conjecture about the polar divisors of the
function.  An interesting feature of their approach is that
for higher moments it seems to predict  lower order terms
of the form $c T^A$ with $\frac12<A<1$, while in this paper
we conjecture that our main terms are valid with an error
of size~$O(T^{\frac12+\varepsilon})$.  The cubic moment
of quadratic Dirichlet $L$-functions is a specific case
for which there is a conjectured lower order term~\cite{Z}
which possibly could be tested numerically.

There are many theorems dealing with moments of $L$-functions in
particular families.  The technique to prove these theorems usually
involves invoking an approximate functional equation and averaging
the coefficients of the $L$-function over the family. The averaging process
behaves like a harmonic detection device. This harmonic detector usually
presents itself as a formula with a relatively simple part and a somewhat
more complicated part that is smaller in the first approximation.
In the theorems in the literature it is often the case that the simple
part of the harmonic detector is sufficiently good to determine the first
or second moment of the family. The terms involved here are usually
called the ``diagonal'' terms. But invariably the more complicated
version is needed to determine the asymptotics of the third or fourth
moments; in these situations one has gone ``beyond the diagonal.''
In at least one situation (fourth moment of cusp form $L$-functions)
it has been necessary to identify three separate stages of more subtle
harmonic detection: the first featuring diagonal term contributions and
the second and third featuring contributions to the main terms by two
different types of off-diagonal terms.  We believe that as one steps up
the moments of a family then at every one or two steps a new type of
off-diagonal contribution will emerge. The whole process is poorly
understood; we only have glimpses of a mechanism but no clear idea of how
or why it works.

        It is remarkable that all of these complicated harmonic detection
devices ultimately lead to very simple answers, as detailed in this paper.
It is also remarkable that there are only three or four different types
of symmetries; families with the same symmetry type often have different
harmonic detectors, with different wrinkles at each new stage of
off-diagonal, but somehow lead to answers which are structurally the same.
It would be worthwhile to understand how this works.

        Finally, we comment that the recipe we develop in this paper only
uses the simplest diagonal harmonic detectors. Our formulas are expressed
as combinatorial sums arising only from diagonal terms. We are well aware
of the off-diagonal pieces, and we do not understand how they
cancel and combine. What we do understand and what we are presenting
here is a conjecture for the final simple answer that should emerge after 
all of the complicated
cancellations between the increasingly subtle off-diagonal terms
are taken into account.
The reader needs to be aware of this to understand the goals and contents
of this paper.


The paper is organized as follows. In the remainder of this section we
give a detailed comparison between $L$-functions and
characteristic polynomials of unitary matrices,
summarize our previous work on the leading terms in the mean values
of $L$-functions,
and
describe the more general moments considered in this
paper.  This allows
us to state our main results and conjectures, which are
given in Section~\ref{sec:mainconjectures}.  We then give a detailed comparison
with known results for the Riemann $\zeta$-function.

In Section~\ref{sec:taspectrecipe} we give a detailed derivation
of our conjectures in the case of moments on 
the critical line
of a single $L$-function.  We first write the conjecture in
terms of a function defined by an infinite sum, and then write
it as an Euler product and identify the leading-order poles.
The local factors are also written in a concise form which
is more suitable for computation.
Both the $L$-function and random matrix calculations lead to
expressions involving a sum over a set of partitions.  These sums
can be written in a concise form involving contour integrals,
as described in Section~\ref{sec:concisesums}.  We also show that the original
results of Keating and Snaith \cite{KS1,KS2} for the leading order
term can be re-derived from the present work.  In addition, we
express the arithmetic factor in the moments of the Riemann
zeta-function in an explicit form.

In Section~\ref{sec:families} we describe a particular notion of a family of
$L$-functions which can be used to give a unified treatment of all
of the mean values we have considered. These families are central
to our method of conjecturing mean values and we give a detailed
description of the method in Section~\ref{sec:recipe}. 
 As explicit examples we give the
details of the calculations for $L$-function families with
Unitary, Symplectic, and Orthogonal symmetry.

In Section~\ref{sec:numerics} we give numerical approximations for the
coefficients in our conjectured mean values.
We then report on numerical calculations of representative cases of
the conjectures.  Good agreement is found.

The calculations of the random matrix averages, which are
based in part on \cite{BF} and \cite{BH}, are
complicated but elementary.  Those results have been presented
in~\cite{CFKRS}.  In subsequent papers we will also
present a fuller discussion of the terms which appear in our
conjectures, give some more general conjectures,
and describe the algorithms behind our numerical calculations.

The authors thank P.~Forrester, R.~Heath-Brown, C.~Hughes, N.~Katz, P.~Michel,
and P.~Sarnak for many helpful discussions.

\subsection{  Properties of $L$-functions }
\label{sec:propofLfunctions}

We present the definition and key properties of $L$-functions.
These properties are familiar, but a summary will be useful in our discussion
of mean values and for the comparison with
the characteristic polynomials of random matrices.

The definition of an $L$-function which we give below is a slight
modification of what has come to be called the ``Selberg class''
\cite{S,CG3,Mur} of Dirichlet series. Let $s=\sigma+it$ with
$\sigma$ and $t$ real. An {\dfn $L$-function} is a Dirichlet
series
\begin{equation}
L(s)=\sum_{n=1}^\infty \frac{a_n}{n^s},
\end{equation}
with $a_n \ll_\varepsilon n^\varepsilon$ for every $\varepsilon>0$, which
has three additional properties.

\myitem{Analytic continuation:}  $L(s)$ continues to a meromorphic
function of finite order with at most finitely many poles, and all
poles are located on the $\sigma=1$ line.

\myitem{Functional equation:} There is a number $\varepsilon$
with $|\varepsilon|=1$, and a function $\gamma_L(s)$ of the form
\begin{equation}
\gamma_L(s)= P(s) Q^s \prod_{j=1}^w \Gamma(w_j s+ \mu_j) ,
\label{eqn:gammafactors}
\end{equation}
where $Q>0$, $w_j >0$, $\Re \mu_j \ge 0$, and $P$ is a
polynomial whose only zeros in $\sigma>0$ are at the poles of
$L(s)$, such that
\begin{equation}
\xi_L(s):= \gamma_L(s)L(s)
\end{equation}
is entire, and
\begin{equation}
\xi_L(s)    =\varepsilon \overline{\xi_L}(1-s),
\end{equation}
where $\overline{\xi_L}(s) = \overline{\xi_L (\overline{s})}$ and
$\overline{s}$ denotes the complex conjugate of~$s$ .

The number $2 \sum_{j=1}^w w_j$ is called the {\it degree} of the $L$-function,
and this is conjectured to be an integer.  It is conjectured furthermore 
that each $w_j$ can be taken to equal $\frac12$, so $w$ equals
the degree of the $L$-function. 

For the calculations we do in this paper, it is
convenient to write the functional equation in
asymmetric form:
\begin{equation}
L(s)=\varepsilon X_L(s) \overline{L}(1-s),
\end{equation}
where $ \displaystyle X_L(s) =
\frac{\overline{\gamma_L}(1-s)}{\gamma_L(s)}  . $ Also we define
the ``$Z$-function'' associated to an $L$-function:
\begin{equation}
Z_L(s):= \varepsilon^{-\frac 1 2}X_L^{-\frac12}(s) L(s), \label{eqn:ZL}
\end{equation}
which satisfies the functional equation
\begin{equation}
Z_L(s)= \overline{Z_L}(1-s).
\end{equation}
Note that here we define $Z_L$ as a function of a complex
variable, which is slightly different from the standard notation.
Note also that $Z_L(\tfrac{1}{2} +it)$ is real when $t$ is real,
$X_L(\frac12)=1$, and $|X_L(\frac12 + it)|=1$ if $t$ is real.

\myitem{Euler product:}  For $\sigma>1$ we have
\begin{equation}
L(s)=\prod_p L_p(1/p^s),
\end{equation}
where the product is over the primes~$p$, and
\begin{equation}
L_p(1/p^s)=\sum_{k=0}^\infty \frac{a_{p^k}}{p^{k s}} =
\exp\left(\sum_{k=1}^\infty \frac{b_{p^k}}{p^{k s}} \right) ,
\end{equation}
where $b_n\ll n^\theta$ with $\theta <\frac12$.

Note that $L(s)\equiv 1$ is the only constant $L$-function,
 the set of $L$-functions is closed under products,
and if $L(s)$ is an $L$-function then so is $L(s+iy)$ for
any real~$y$.  An $L$-function is called {\dfn primitive}
if it cannot be written as a nontrivial product of
$L$-functions, and it can be shown, assuming Selberg's orthonormality conjectures,
 that any $L$-function
has a unique representation as a product of primitive
$L$-functions. See~\cite{CG3}. It is believed that $L$-functions
only arise from arithmetic objects, such as characters~\cite{Dav},
automorphic forms~\cite{Iw1,Iw2}, and automorphic
representations~\cite{BC,Bu}. Very little is known about
$L$-functions beyond those cases which have been shown to be
arithmetic.

There are several interesting consequences of the above properties,
and many conjectures which have been established in
 a few (or no) cases.  We highlight some additional properties
of $L$-functions and then discuss their random matrix analogues.

\myitem{Location of zeros:} Since $\xi_L(s)$ is entire,
$L(s)$ must vanish at the poles
of the $\Gamma$-functions in the $\gamma_L$ factor.  These are
known as the {\dfn trivial zeros} of the $L$-function.  By the
functional equation and the Euler product, the only other possible
zeros of $L(s)$ lie in the {\dfn critical strip} $0\le \sigma\le 1$.
By the argument
principle, the number of nontrivial zeros with $0<t<T$ is
asymptotically $(W/\pi)T\log T$, where $W=\sum w_j$.
The {\dfn Riemann Hypothesis} for $L(s)$ asserts that the nontrivial zeros
of $L(s)$ lie on the {\dfn critical line} $\sigma=\frac12$.
The much weaker
(but still deep) assertion that $L(s)\not=0$ on $\sigma=1$ has
been proven for arithmetic $L$-functions~\cite{JS},
which can be viewed as a generalization of the prime number theorem.

\myitem{Average spacing of zeros:} By the zero counting result
described above, the average gap between consecutive zeros of
$L(s)$ with imaginary part around~$T$ is $\pi/(W\log T)$.

\myitem{Zeros of derivatives:} If the Riemann Hypothesis is
true then all zeros of the
derivative $\xi^{\prime}(s)$ lie on the critical line, while all
zeros of $\zeta^{\prime}(s)$ lie to the right of the
critical line~\cite{LM}.

\myitem{Critical values:} The value $L(\frac12)$ is called the
{\dfn critical value} of the $L$-function.  The significance of
$s=\frac12$ is that it is the symmetry point of the functional
equation.  The mean values we study in this paper are averages of
(powers of) critical values of $L$-functions, where the average is
taken over a ``family'' of $L$-functions.  Examples of families
and their corresponding mean values are given in Section~\ref{sec:examplefamilies}.

Note.  If the set $\{\mu_j\}$ is stable under complex conjugation
and the $a_n$ are real, then $\varepsilon$ is commonly called the
{\dfn sign of the functional equation}. If the sign is $-1$ then
$L(s)$ has an odd order zero at $s=\frac12$; more generally, if
the sign is not~$1$ then $L(\frac12)=0$. When $L(\frac12)$
vanishes, it is common to use the term `critical value' for the
first nonzero derivative $L^{(j)}(\frac12)$, but in this paper we
use `critical value' to mean `value at the critical point.'

\myitem{Log conductor:} We measure the ``size'' of an $L$-function
by its \emph{log conductor}, defined as 
$c(L)=cond(L)=|X_L'(\frac12)|$.
The conductor of an $L$-function has a conventional meaning in
many contexts, and the log  conductor is a 
 simple function of the (logarithm of the) usual conductor.
Other authors use similar names, such as ``analytic conductor'',
for similar quantities.
By the argument principle, the density of zeros near the critical
point is $2\pi c(L)^{-1}$.

\myitem{Approximate functional equation:}
A standard tool for studying analytic properties of $L$-functions
is an approximate functional equation for $L(s)$, which expresses
the $L$-function as a sum of two Dirichlet series involving
the Dirichlet coefficients of $L$ multiplied by a smoothing
function. See for example~\cite[5.3]{IK}. 
For the purposes
of the heuristics we develop, we use a sharp cutoff
and don't concern ourselves with the remainder,
\begin{equation}
L(s)= \sum_{m < x} \frac{a_m}{m^s} + \varepsilon X_L(s)
\sum_{n < y} \frac{\overline{a_n}}{n^{1-s}} + remainder.
\end{equation}
Here the product $xy$ depends on parameters in the functional
equation. The name comes from the fact that the right side looks
like $L(s)$ if $x$ is large, and like $\varepsilon X_L(s)
\overline{L}(1-s)$ if $x$ is small, which suggests the asymmetric
form of the functional equation.

The approximate functional equation is the starting point of our approach
to conjecturing the moments of $L$-functions.  This is described in
Sections~\ref{sec:therecipe} and~\ref{sec:generalrecipe}.


\subsection{  Properties of characteristic polynomials
}\label{sec:propofcharpolys}

With the exception of the Euler product, all of the properties
of $L$-functions have a natural analogue in the  characteristic
polynomials of unitary matrices.
We note each property in turn.

Let
\begin{equation}\Lambda (s)=\Lambda_A(s)=\det(I-A^*s)
=\prod_{n=1}^N \(1-s e^{-i\theta_n}\)
\end{equation}
denote the characteristic polynomial of an $N\times N$ matrix $A$.
Throughout the paper we assume that
$A$ is unitary
(i.e. $A^{*}A=I$ where $A^{*}$ is the Hermitian
conjugate of $A$),  so the eigenvalues of $A$ lie on the unit circle
and can be denoted by~$e^{i\theta_n}$.

Note: in our previous paper~\cite{CFKRS} we used a different 
definition of the characteristic polynomial.  

We can express $\Lambda (s)$ in expanded form:
\begin{equation}
\Lambda (s)=\sum_{n=0}^N a_n s^n, \label{eqn:lambdasum}
\end{equation}
 which corresponds to the Dirichlet series representation for
$L$-functions.

\myitem{Analytic continuation:} Since $\Lambda (s)$ is a
polynomial, it is an entire function.

\myitem{Functional equation:}  Since $A$ is unitary, we have
\begin{equation}\Lambda_A (s)=(-1)^N
\det A^* \
s^N
\det(I-A s^{-1}),
\end{equation}
and so, writing
\begin{equation}\det A =e^{i\phi}
\end{equation}
(where unitarity implies that $\phi\in {\mathbb R}$), we have
\begin{eqnarray}\Lambda_A (s)&=&(-1)^N \det A^* s^N \Lambda_{A^*}(\tfrac 1{s})\cr 
&=&(-1)^N e^{-i\phi}s^N \, \overline{\Lambda_{A} }(\tfrac 1{s}).
\label{eqn:lambdafe}
\end{eqnarray}
This plays the same role for $\Lambda (s)$ as the functional
equation for $L$-functions: it represents a symmetry with respect
to the unit circle ($s=re^{i\alpha} \to \frac 1{s}=\frac 1 r
e^{-i\alpha}$).

Also let
\begin{equation}
\mathcal{Z}_A(s)=\((-1)^N e^{-i\phi}\)^{-\frac12} s^{-N/2} \Lambda_A(s),\label{eqn:ZLambda}
\end{equation}
in direct analogy to (\ref{eqn:ZL}), 
the sign of the functional equation $\varepsilon$ being identified with
$(-1)^Ne^{-i\phi}=(-1)^N\det A^*$.  The functional equation becomes
\begin{equation}
\mathcal{Z}_A(s)=\overline{\mathcal{Z}}_A(\tfrac{1}{s}).
\end{equation}
Note that this implies $\mathcal{Z}$ is real on the unit circle,
and in analogy to the $X_L$ factor from
an $L$-function, the factor $s^{-N/2}$ equals 1 at
the critical point $s=1$,
and has absolute value~1 on the unit circle.  

\myitem{Location of zeros:} Since  $A$ is unitary, its eigenvalues
all have modulus~1, so the zeros of $\Lambda (s)$ lie on the unit
circle (i.e.~the Riemann Hypothesis is true).  The unit circle
is the ``critical line'' for~$\Lambda (s)$.

\myitem{Average spacing of zeros:} Since the $N\times N$ matrix
$A$ has $N$ eigenvalues on the unit circle, the average spacing
between zeros of $\Lambda _A(s)$ is $2\pi/N$.  

\myitem{Zeros of derivatives:} Since the zeros of $\Lambda (s)$
lie on the unit circle, the zeros of the derivative $\Lambda
^{\prime}$ lie inside the unit circle. This follows from the
general fact that the zeros of the derivative of a polynomial lie
in the convex hull of the zeros of the polynomial.

\myitem{Critical values:}  The critical point for $\Lambda (s)$ is
the symmetry point of the functional equation
$s=1=e^{i\cdot 0}$, and $\Lambda (1)$ is the critical value.

\myitem{Conductor:} In analogy with the case of $L$-functions,
we define the conductor of $\Lambda$ to be the (absolute value of the)
derivative of 
the factor in the asymmetric form of the functional equation, evaluated
at the critical point $s=1$.  That is, the conductor of $\Lambda$ is~$N$.
Also in analogy to the case of $L$-functions, the density of
zeros on the unit circle is $2\pi/N$.

When modeling a
family of $L$-functions, we choose $N$ so that $L$ and $\Lambda$ have
the same conductor.  Equivalently, $L$ and $\Lambda$ have
the same density of zeros near the critical point.

\myitem{Approximate functional equation:}
Substituting the polynomial (\ref{eqn:lambdasum}) into the functional equation
(\ref{eqn:lambdafe}), we have
\begin{equation}\sum_{n=0}^N a_n s^n =(-1)^N e^{-i\phi}\sum_{n=0}^N \overline{a_n}
 s^{N-n},
\end{equation}
and so
\begin{equation}a_n =(-1)^N e^{-i\phi}\, \overline{a_{N-n}}.
\end{equation}
Hence, when $N$ is odd, we have
\begin{equation}
\Lambda (s) =\sum_{m=0}^{\frac {N-1}{2}}a_m s^m 
+ (-1)^N e^{-i\phi} s^N \sum_{n=0}^{\frac{N-1}{2}}\overline{a_n} s^{-n} ,
\end{equation}
which corresponds to the approximate functional equation for $L$-functions.
When $N$ is even, there is an additional
term: $a_{\frac N 2 } s^{\frac N 2}$.

Although we use the approximate functional equation in our
calculations for $L$-functions,
in our previous paper \cite{CFKRS}
we use other methods for the characteristic polynomials.
In principle,
it would be possible to use
the approximate functional equation and compute
averages of products of the coefficients~$a_n$.
Such a calculation would, presumably, mirror that for
the $L$-functions.
This would appear to be more cumbersome than the
approach taken in \cite{CFKRS}, but might merit further investigation.

The above discussion applies to any unitary matrix.  We also
consider matrices which, in addition to being unitary, are either
symplectic or orthogonal. We use these three ensembles of matrices
to model families of $L$-functions.  While the notion of ``family
of $L$-functions''  has not yet been made precise, we give several
natural examples in the next section.

Associated to each family is a ``symmetry type'' which identifies
the matrix ensemble that will be used to model the family. This
correspondence is most easily seen in terms of the sign of the
functional equation, which is analogous to the determinant of the
matrix. If $A$ is unitary symplectic, then $\det A=1$ (i.e.
$\phi=0$), and if $A$ is orthogonal, then $\det A=\pm 1$.
Correspondingly, the functional equations for $L$-functions with
unitary symmetry involve a (generally complex) phase factor,
whereas for $L$-functions with symplectic symmetry this phase
factor is unity, and in the case of orthogonal symmetry it is
either $+1$ or $-1$.

While the sign of the functional equation can sometimes 
suggest the symmetry type of the family, in general it requires
a calculation to determine the symmetry type.  One possible
calculation is to determine
the moments of the family near the critical point,
as described in this paper.
Comparison with the corresponding
random matrix average can then be used to determine the symmetry type.
Another possibility is to determine the density
of the low-lying zeros of the family.


\subsection{  Example families and moments of $L$-functions
}\label{sec:examplefamilies}

We now give examples of families of primitive  
$L$-functions and describe the
associated mean values.  The families we consider here are of a
special form, 
which is described in Section~\ref{sec:families}. In preparation
for the comparison with random matrices in the next section, we
will classify the example families according to their symmetry
type: Unitary, Orthogonal, and Symplectic.
For the Orthogonal symmetry type we recognize three cases:
$SO$, $O^-$, and $O$, corresponding respectively to Orthogonal families 
in which the functional equation has $\varepsilon=1$, or
$\varepsilon=-1$, or $\varepsilon=\pm 1$ equally often.
 Note that each
family is a partially ordered set, and the order 
is determined by a quantity called the ``conductor'' of the
$L$-function. The mean values given below are conjectural for all
but a few small values of $k$.
For a general discussion of these mean values
and some more examples, see \cite{CF}.

Unitary examples:

\myitem{1)} $\{L(s+ iy)\ |\ y\ge 0\}$, ordered by~$y$,
where $L(s)$ is any primitive $L$-function.
These are the only known continuous families of $L$-functions (Sarnak's rigidity conjecture).

\myitem{2)} $\{L(s,\chi)\ |\ q \hbox{ a positive integer,
        $\chi$ a primitive character $\bmod q$}\}$ ordered by~$q$.

An example conjectured  mean value for integer $k$ is:
\begin{equation}  \int_0^T |\zeta(\tfrac12 +it)|^{2k}\,dt  = T \, {\mathcal P}_{k}(\log T) +O(T^{\frac12 +\varepsilon}),
\label{eqn:zeta2k}
\end{equation}
for some polynomial ${\mathcal P}_{k}$ of degree $k^2$ with
leading coefficient $g_k a_k/{k^2!}$, where
\begin{eqnarray}
a_k& =  &\prod_p \left(1-\tfrac{1}{p}\right)^{k^2}
\sum_{m=0}^{\infty} \left( \frac{\Gamma(m+k)} {m! \Gamma(k)}
\right)^2 p^{-m} \nonumber \\ & =  &\prod_p\left(1-\tfrac 1p\right)^{(k-1)^2}
\ \sum_{j=0}^{k-1}\binom{k-1}{j}^2p^{-j}  \label{eqn:zeta2kak} \end{eqnarray}
and
\begin{equation}
g_k=  k^2 ! \prod_{j=0}^{k-1}\frac{j!}{(k+j)!}.
\label{eqn:zeta2kgk}
\end{equation}
(The placement of $k^2!$ is to ensure that $g_k$ is an integer
\cite{CF}.)  The above conjecture 
has been proven for $k=1,\,2$.
See \cite{HL,I,A,Ko,H-B1,C,Mot}.  
When $k=2$ our conjectured error term of $O(T^{\frac12+\varepsilon})$
has only been obtained in the case of a smooth weight function~\cite{Iv1}.

Symplectic examples:

\myitem{3)} $\{L(s,\chi_d)\ |\ d \hbox{ a fundamental discriminant, }
 \chi_d(n)=(\tfrac{d}{n})\}$ ordered by~$|d|$.

\myitem{4)} $\{L(s, \hbox{sym}^2 f)\ |\ f\in S_k(\Gamma_0(1)),
                \hbox{ $k$ a positive even integer}\}$, ordered by~$k$.

An example  conjectured mean value is:
\begin{equation}
\sumstar_{|d|\le D}L(\tfrac12 ,\chi_d)^k= \frac{6}{\pi^2}D\, {\mathcal
Q}_{k}(\log D) +O(D^{\frac12 +\varepsilon}), \label{eqn:sumchidk}
\end{equation}
where $\sumstar$ is over fundamental discriminants,
$\chi_d(n)=(\tfrac{d}{n})$ is the Kronecker symbol,
and the sum is over all real,
primitive Dirichlet characters of conductor up to $D$.
Here ${\mathcal Q}_k$  is a polynomial of degree $k(k+1)/2$,
with leading coefficient $g_k a_k/(k(k+1)/2)!$, where
\begin{equation}
a_k=\prod_p \frac{(1-\frac 1p)^{\frac{k(k+1)}{2}}}
{1+\frac 1p}
\left(\frac{(1-\frac{1}{\sqrt{p}})^{-k}+
(1+\frac{1}{\sqrt{p}})^{-k}}{2}+\frac{1}{p}\right)
\end{equation}
and
\begin{equation}
g_k=  (k(k+1)/2)! \prod_{j=1}^k\frac{j!}{(2j)!} .
\end{equation}
The main term of this conjecture has been proven for
$k=1,\,2,\,3$, and the case of $k=4$ is almost within reach of
current methods. See \cite{GV,J,So}.

Orthogonal examples:

\myitem{5)} $\{L(s,f)\ |\ f\in S_k(\Gamma_0(N)),
        \hbox{ $N$ fixed, $k$ a positive even integer}\}$,
ordered by~$k$.

\myitem{6)} $\{L(s,f)\ |\ f\in S_k(\Gamma_0(N)),
                \hbox{ $k$ fixed, $N$ a positive integer}\}$, ordered by~$N$.

An example  conjectured mean value is:
\begin{equation} 
\sum_{f\in H_2(q)}L_f(\tfrac12,f )^k =\frac{q}{3} \, {\mathcal R}_k(\log q)
+O(q^{\frac12 +\varepsilon}) , \label{eqn:sumlfk}
\end{equation}
where $H_2(q)$ is the collection of Hecke newforms
of weight 2 and squarefree level $q$.
Here ${\mathcal R}_k$ is a polynomial of degree $k(k-1)/2$, with leading
coefficient $g_k a_k/(k(k-1)/2)!$, where
\begin{equation} a_k=\prod_{p\nmid q}
\left(1-\tfrac{1}{p}\right)^{\frac{k(k-1)}{2}}
\frac{2}{\pi}\int_0^\pi  \sin^2\theta
\left(\frac{e^{i\theta}\left(1-\frac{e^{i\theta}}{\sqrt{p}}\right)^{-1}
-e^{-i\theta}\left(1-\frac{e^{-i\theta}}{\sqrt{p}}\right)^{-1}}
{e^{i\theta}-e^{-i\theta}}\right)^k \,d\theta
\end{equation}
and
\begin{equation}g_k = 2^{k-1}(k(k-1)/2)!\prod_{j=1}^{k-1}\frac{j!}{(2j)!} .
\end{equation}
The main term of this conjecture has been proven for
$k=1,\,2,\,3,\,4$, in the case that $q$ is prime. See
\cite{D,DFI,KMV}.  Also see Ivi\'c~\cite{Iv2} for the analogous mean values
for Maass forms.

The above examples are merely meant to give a flavor of
the types of families which are  of current interest.

The above cases, and their random matrix analogues, have been
extensively discussed from the perspective of the {\it leading
terms} in the asymptotic expansions. See \cite{CF,KS1,KS2}. In the
present paper we extend that work to include all of the terms in
the above mean values (ie., all coefficients in the conjectured
polynomials), which we recover from a more general mean value
involving a product of $L$-functions whose arguments are free
parameters. In the next two sections we describe these more
general mean values, discuss their random matrix analogues, and
then state our results and conjectures.


\subsection{  Shifted moments
}\label{sec:shiftedmoments}

A key point in this paper is that the structure of mean values of
$L$-functions is more clearly revealed if one considers the average
of a product of $L$-functions, where each $L$-function is evaluated
at a location slightly shifted from the critical point.  
The example mean values given in the previous section can be obtained
by allowing the shifts to tend to zero.

 Let $\alpha=(\alpha_1,\ldots, \alpha_{2k})$, where throughout
the paper we assume
$ |\Re \alpha_j|< \frac{1}{2}$, 
 and suppose that $g(t)$ is a suitable weight function. The
mean values we consider are
\begin{equation}
I_k(L,\alpha,g)= \int_{-\infty}^\infty
        Z_L({\textstyle \frac12}+\alpha_1+it)\cdots
            Z_L({\textstyle \frac12}+\alpha_{2k}+it)g(t)\,dt ,
\end{equation}
and, with $\alpha=(\alpha_1,\ldots,\alpha_{k})$,
\begin{equation}
S_k({\mathcal F},\alpha,g)= \sum_{L\in{\mathcal F}} Z_L({\textstyle
\frac12}+\alpha_1)\cdots
                    Z_L({\textstyle \frac12}+\alpha_{k})g(c(L)) .
\end{equation}
In the first case it is assumed that $L(s)$ is a primitive $L$-function,
and in the second
$\mathcal F$ is family of primitive $L$-functions partially ordered
by log conductor $c(L)$.

We refer to $g$ as a ``suitable'' weight function, but we leave
that term undefined. An example of a suitable weight function is
$g(x) =  f(x/T)$, where $f$ is real, nonnegative, bounded, and
integrable on the positive real line.

The random matrix analogs of the above expressions are
\begin{equation}
J_k(U(N),\alpha) =\int_{U(N)} {\mathcal Z}_A(e^{-\alpha_{1}}) \cdots
{\mathcal Z}_A(e^{-\alpha_{2k}})\,dA,
\end{equation}
where $\alpha=(\alpha_1,\ldots,\alpha_{2k})$ and the average is
over Haar measure on $U(N)$, and
\begin{equation}
J_k(G(N),\alpha)=\int_{G(N)} {\mathcal Z}_A(e^{-\alpha_1})\cdots {\mathcal
Z} _A(e^{-\alpha_k}) \,dA,
\end{equation}
where $G(N)$ is $USp(2N)$, $O^-(2N)$, or $SO(2N)$ and
$\alpha=(\alpha_1,\ldots, \alpha_{ k})$.  $O^-(2N)$ is defined as
the collection of orthogonal $2N\times 2N$ matrices with
determinant negative one.  Haar measure on $USp(2N)$ and $O(2N)$
determines the weighting for the averages.

In the next section we compare our conjectures for the
$L$-function mean values with exact formulae for the random matrix
averages.


\subsection{  Main results and example conjectures
}\label{sec:mainconjectures}

We state our main results and conjectures here. We give example
conjectures for the full main term in shifted mean values of
number theoretic interest; these examples illustrate our methods
and cover the three symmetry types of families of $L$-functions.
We also give a corresponding theorem about the random matrix
analogue of these mean values for each of the three compact matrix
ensembles.

We present our results in pairs: a conjecture for an $L$-function
mean value, followed by a theorem, quoted from \cite{CFKRS},
for the corresponding average of
characteristic polynomials.  
For each pair the parts of each
formula match according to the scheme described in Section~\ref{sec:propofcharpolys}. 
In particular, the scaling of the large parameter is determined
by equating log conductors.
In
the random matrix formula the integrand contains a term $(1-e^{\pm
z_m - z_\ell})^{-1}$, which has a simple pole at $z_{\ell}= \pm
z_m$. In the $L$-function formula this corresponds to the term
containing all of the arithmetic information, which is of the form
$\zeta(1+z_i \mp z_j)$ times an Euler product, and so also has
a simple pole at $z_i =\pm z_j$.

The formulae are written in terms of contour integrals and involve
the Vandermonde:
\begin{equation}
\Delta(z_1,\dots,z_{m})=\prod_{1\le i < j\le m}(z_j-z_i).\label{eqn:vandermonde}
\end{equation}
We also set $e(z)=e^{2\pi i z}$.

\begin{conjecture} \label{thm:zeta2kconjecture}
 Suppose $g(t)$ is a suitable weight function. Then
\begin{equation}
\int_{-\infty}^\infty |\zeta(\tfrac12 +it)|^{2k}g(t)\,dt= \int_{-\infty}^\infty
P_k\left(\log \tfrac{t}{2 \pi}\right)
(1+O(t^{-\frac{1}{2}+\varepsilon})) g(t) \,dt ,
\end{equation}
where $P_k$ is the polynomial of degree $k^2$ given by the
$2k$-fold residue
\begin{equation}
P_k(x)= \frac{(-1)^k}{k!^2}\frac{1}{(2\pi i)^{2k}} \oint\cdots
\oint \frac{G(z_1, \dots,z_{2k})\Delta^2(z_1,\dots,z_{2k})}
{\displaystyle \prod_{j=1}^{2k} z_j^{2k}} 
e^{\tfrac x2 \sum_{j=1}^{k}z_j-z_{k+j}}\,dz_1\dots dz_{2k} ,
\end{equation}
where one integrates over small circles about $z_i=0$, with
\begin{equation}
G(z_1,\dots,z_{2k})= A_k(z_1,\dots,z_{2k})
\prod_{i=1}^k\prod_{j=1}^k\zeta(1+z_i-z_{k+j}) ,
\end{equation}
and $A_k$ is the Euler product
\begin{equation}
A_k(z) =\prod_p \prod_{i=1}^k\prod_{j=1}^k
\left(1-\frac{1}{p^{1+z_i-z_{k+j}}}\right) \int_0^1 \prod_{j=1}^k
\left(1-\frac{e(\theta)}{p^{\frac12 +z_j}}\right)^{-1}
\left(1-\frac{e(-\theta)}{p^{\frac12 -z_{k+j}}}\right)^{-1}\,d\theta.
\end{equation}
More generally,
\begin{equation}
I_k(\zeta,\alpha,g)= \int_{-\infty}^\infty 
P_k\left(\log \tfrac{t}{2 \pi} , \alpha\right) (1+O(t^{-\frac{1}{2}+\varepsilon})) g(t)\,dt ,
\end{equation}
where
\begin{equation}
P_k(x,\alpha)=
 \frac{(-1)^k}{k!^2}\frac{1}{(2\pi i)^{2k}}
\oint \cdots \oint \frac{G(z_1,
\dots,z_{2k})\Delta^2(z_1,\dots,z_{2k})} {\displaystyle
\prod_{j=1}^{2k} \prod_{i=1}^{2k}(z_j-\alpha_i)} 
e^{\tfrac x2 \sum_{j=1}^{k}z_j-z_{k+j}}\,dz_1\dots dz_{2k},
\end{equation}
with the path of integration being small circles surrounding the poles
$\alpha_i$.
\end{conjecture}

A general version of the above conjecture is given in
Conjecture~\ref{thm:general2kconjecture}.

\begin{theorem} \label{thm:UN2k}

  In the notation of Section~\ref{sec:shiftedmoments} we have
\begin{eqnarray}
J_k(U(N),0) =\prod_{j=0}^{k-1} \left(\frac{j!}{(k+j)!}
\prod_{i=1}^{k}(N+i+j) \right). \end{eqnarray}
More generally, with 
$$
G(z_1,\dots,z_{2k})=\prod_{i=1}^{k}\prod_{j=1}^k (1-e^{-z_i+z_{j+k}})^{-1},
$$ 
we have
\begin{equation}
J_k(U(N),\alpha)  =  \frac{(-1)^{k}}{k!^2} \frac {1}{(2\pi i)^{2k}}
 \oint\cdots \oint \frac {
G(z_1,\dots,z_{2k})
\Delta^2(z_1,\ldots,z_{2k})} { \displaystyle \prod_{i=1}^{2k}
\prod_{j=1}^{2k} (z_j-\alpha_i)} 
e^{\frac{N}{2}
\sum_{j=1}^k z_j-z_{k+j} }dz_1\cdots dz_{2k}  . 
\end{equation}
\end{theorem}

Comments on the formulae:
\begin{enumerate}

\item Let $\alpha_i\to0$ in the second part of the Conjecture to obtain
the first part of the Conjecture.

\item 
The structures of $J_k(U(N),\alpha)$ and $P_k(x,\alpha)$ are identical
in that the functions $G(z_1,\dots,z_{2k})$ have simple poles 
at~$z_i=z_{k+j}$.

\item The local factors of $A_k(\alpha)$ are polynomials
in $p^{-1}$ and $p^{-\alpha_i}$, $i=1,\ldots,k$, as seen from
Theorem~\ref{thm:zeta2kAk}. Since $A_k(\alpha)$ comes from a symmetric
expression, it is also a polynomial in $p^{\alpha_{i+k}}$,
$i=1,\ldots,k$. This is discussed in Section~\ref{sec:explicitarithmeticfactor}.  Note also that
$a_k$ in (\ref{eqn:zeta2kak}) equals $A_k(0,\ldots,0)$,
as shown in Section~\ref{sec:leadingorder}.

\item That $P_k(x)$ is actually a polynomial of degree $k^2$
can be seen by considering the order of the pole at $z_j=0$. We wish to extract
from the numerator of the integrand, the coefficient of $\prod z_i^{2k-1}$, a polynomial
of degree $2k(2k-1)$.
The Vandermonde determinant squared is a homogeneous polynomial of degree
$2k(2k-1)$. However, the poles coming from the $\zeta(1+z_i-z_{k+j})$
cancel $k^2$ of the Vandermonde factors. This requires us, in computing the residue,
to take, in the Taylor expansion of $\exp(\frac{x}{2}\sum_1^k z_j-z_{k+j})$, terms up
to degree $k^2$.

\item The fact that $P_k(\log\tfrac{t}{2\pi})$ is a polynomial in $\log\tfrac{t}
{2\pi}$ of degree $k^2$ corresponds nicely to the formula for
$J_k(U(N),0)$, which is a polynomial of degree $k^2$ in~$N$.
Equating the density of the Riemann zeros at height $t$ with the
density of the random matrix eigenvalues suggests the familiar equivalence
$N=\log\tfrac{t}{2\pi}$.  In this paper we view this as equating conductors.

\item The leading term of
$P_k(x)$ coincides with the leading term conjectured by Keating
and Snaith (see Section~\ref{sec:leadingorder}).
The full polynomial $P_k(x)$ agrees,
when $k=1$ and $k=2$, with known theorems (see Sections~\ref{sec:secondmoment}
and~\ref{sec:fourthmoment}).

\item We can recover the polynomial ${\mathcal P}_k$ in (\ref{eqn:zeta2k}) from $P_k$ by taking
$g(t) = \chi_{[0,T]}(t)$ in the conjecture.

\item The multiple integrals in the Theorem and the Conjecture
can be written as combinatorial sums.  See Section~\ref{sec:taspectrecipe}
where a detailed derivation of our conjecture is given.

\item Our conjecture concerning the order of the error term is
based on our numerical calculations (see Section~\ref{sec:numerics}) and
examination of examples in the literature.

\end{enumerate}

\begin{conjecture} \label{thm:Lhalfchidconjecture}

  Suppose $g(u)$ is a suitable weight
function with support in either $(0,\infty)$ or $(-\infty,0)$,
and let $X_d(s)=|d|^{\frac12-s}X(s,a)$ where
$a=0$ if $d>0$   and  $a=1$ if $d<0$, and
\begin{equation}
X(s,a) = \pi^{s-\frac12}
        \frac{\Gamma\(\frac{1+a-s}{2}\)}{\Gamma\(\frac{s+a}{2}\)} .
\end{equation}
That is, $\chi_d(s)$ is the factor in the functional equation
$L(s,\chi_d)=\varepsilon_d X_d(s)L(1-s,\chi_d)$.
Summing over fundamental discriminants $d$ we have
\begin{equation}
\sumstar_{\!d}L(\tfrac12 ,\chi_d)^kg(|d|)= \sumstar_{\!d}\,Q_k(\log
{|d|})(1+O(|d|^{-\frac{1}{2}+\varepsilon})) g(|d|)
\end{equation}
where $Q_k$ is the polynomial of degree $k(k+1)/2$ given by
the $k$-fold residue
\begin{equation}
Q_k(x)= \frac{(-1)^{k(k-1)/2}2^k}{k!}
\frac{1}{(2\pi i)^{k}}
\oint \cdots \oint
\frac{G(z_1, \dots,z_{k})\Delta(z_1^2,\dots,z_{k}^2)^2}
{\displaystyle \prod_{j=1}^{k} z_j^{2k-1}}
e^{\tfrac x2 \sum_{j=1}^{k}z_j}\,dz_1\dots dz_{k} ,
\end{equation}
where
\begin{equation}
G(z_1,\dots,z_k)=A_k(z_1,\dots,z_k) 
\prod_{j=1}^k X(\tfrac12+z_j,a)^{-\frac12}
\prod_{1\le i\le j\le k}\zeta(1+z_i+z_j),
\end{equation}
and $A_k$ is the Euler product, absolutely convergent for
$|\Re z_j|<\frac12 $, defined by
\begin{eqnarray}
A_k(z_1,\dots,z_k) & =  &\prod_p \prod_{1\le i \le j \le k}
\left(1-\frac{1}{p^{1+z_i+z_j}}\right) \nonumber \\&   &\times \left(\frac
12 \left(\prod_{j=1}^k\left( 1-\frac{1}{p^{\frac 12+z_j}}\right)^{-1} +
\prod_{j=1}^k\left(1+\frac{1}{p^{\frac12+z_j}}\right)^{-1}
\right)+\frac 1p \right) \left( 1+ \frac{1}{p}\right)^{-1}. \end{eqnarray}
More generally, if $\mathcal F$ is the family of real primitive Dirichlet
$L$-functions then
\begin{equation}
S_k({\mathcal F},\alpha,g) = \sumstar_{\!d}\,
Q_k(\log {|d|} , \alpha ) (1+O(|d|^{-\frac12 + \varepsilon}))  g(|d|),
\end{equation}
in which
\begin{eqnarray}
Q_k(x, \alpha)& =  &\frac{(-1)^{k(k-1)/2}2^k}{k!} \frac{1}{(2\pi
i)^{k}} \nonumber \\
&   &\ \ \times \oint \cdots \oint 
\frac{G(z_1,
\dots,z_{k})\Delta(z_1^2,\dots,z_{k}^2)^2 \prod_{j=1}^{k} z_j}
{\displaystyle \prod_{i =1}^{k} \prod_{j=1}^{k} (z_j -
\alpha_i)(z_j+\alpha_{i}) } 
e^{\tfrac x2 \sum_{j=1}^{k}z_j}\,dz_1\dots dz_{k}, \end{eqnarray}
where the path of integration encloses the $\pm \alpha$'s.
\end{conjecture}

\begin{theorem} \label{thm:USp2Nk}

   In the notation of Section~\ref{sec:shiftedmoments} we have
\begin{eqnarray}
J_k(USp(2N),0) =\left( 2^{k(k+1)/2} \prod_{j=1}^k \frac{j!}{(2j)!}
\right) \prod_{1\leq i\leq j\leq k} (N+\tfrac{i+j}{2}). 
\end{eqnarray}
More generally, with 
$$
G(z_1,\dots,z_k) =  \prod_{1\leq i\le j \leq k}
(1-e^{-z_{i}-z_{j}})^{-1}
$$
we have
\begin{eqnarray}
J_k(USp(2N),\alpha) & =  &
\frac{(-1)^{k(k-1)/2}2^k}{k!}\frac{1}{(2\pi i)^{k}} \nonumber \\
& &\ \ \times  \oint\cdots \oint 
\frac{ G(z_1,\dots,z_k) \Delta(z_1^2,\ldots,z_{k}^2)^2 
\prod_{j=1}^{k} z_j} { \displaystyle \prod_{i=1}^{k} \prod_{j=1}^{k}
(z_j-\alpha_i)(z_j+\alpha_i)}e^{N\sum_{j=1}^k z_j } dz_1\cdots
dz_{k}  , 
\end{eqnarray} 
where the contours of integration enclose the $\pm \alpha$'s.
\end{theorem}

Comments:
\begin{enumerate}

\item When comparing Theorem~\ref{thm:USp2Nk} with 
Conjecture~\ref{thm:Lhalfchidconjecture}, equating log conductors
(i.e., the density of zeros) gives
the equivalence 
\begin{equation}
2 N= cond(d):=\log ({|d|}/\pi) +(\Gamma'/\Gamma)(\tfrac14+a).
\end{equation}
The conductor we use here should be contrasted with the ``usual''
conductor associated with Dirichlet $L$-functions:~$\log(|d|/\pi)-\log(2)$. 
We believe this difference is significant, so we discuss it briefly.

The following
manipulations show that our conductor arises naturaly.  
In the derivation of the conjecture, one encounters
the function $X_d(\frac12+z)^{-\frac12}$, which can we rewritten in
several ways:
\begin{eqnarray}\label{eqn:conductormanipulation}
X_d(\tfrac12+z)^{-\frac12} &=& e^{\frac12 \log d \cdot z }
X(\tfrac12+z,a)^{-\frac12} \cr
&=& e^{\frac12 cond(d) \cdot z } \mathcal{G}(z) ,
\end{eqnarray}
where $\mathcal{G}(z)=1+O(z^3)$.  In the statement of the conjecture we 
used the first line of~(\ref{eqn:conductormanipulation}), 
incorporating the product over 
$X(\tfrac12+z,a)^{-\frac12}$ into the factor $G(z_1,\ldots,z_k)$.
If we chose instead to use the second line 
of~(\ref{eqn:conductormanipulation}), then the conjecture
would be written as a sum over $Q_k(cond(d))$.
One would still have that $Q_k$ is a polynomial of degree $k(k+1)/2$.
Since $\mathcal{G}(z)=1+O(z^3)$, the first $3$ leading terms
in that polynomial would not explicitly depend on the factor 
$X_d$ from the functional equation, although the lower degree terms would.
This phenomenon does not occur for moments of $L$-functions
in $t$-aspect.

\item $\mathcal{Q}_k$ in (\ref{eqn:sumchidk}) can be recovered from $Q_k$ above
by taking $g(|d|)= \chi_{[0,D]}(|d|)$, and using the estimate
$\sumstar_{-D<d<0} 1 = 3D/\pi^2 +O(D^{\frac12  +\varepsilon})$; the same
estimate holds for positive~$d$.

\item A heuristic derivation of Conjecture~\ref{thm:Lhalfchidconjecture} is given in
Section~\ref{sec:symplectic}. 

\item The leading term of $Q_k$ coincides with the leading term
conjectured by Keating and Snaith \cite{KS2}.  The calculation is
analogous to the one given in Section~\ref{sec:leadingorder}.

\end{enumerate}

\begin{conjecture} \label{thm:sumLfkconjecture}
Suppose $q$ is squarefree, let $H_n(q)$ be the set of newforms in $S_n(\Gamma_0(q))$, 
and let
\begin{equation}
X_{n,q}(s)=\(\frac{q\mathstrut}{4\pi^2}\)^{\frac12-s}
\frac{\Gamma(\frac12-s+\frac{n}{2})}{\Gamma(s-\frac12+\frac{n}{2})}
\end{equation}
be the factor in the functional equation
$L_f(s)=\varepsilon_{n,q} X_{n,q}(s) L_f(1-s)$ for the $L$-functions
associated to $f\in H_n(q)$.
Then
\begin{equation} \label{eqn:Hnqconjecture}
\sum _{  f\in H_n(q)   } {L_f(\tfrac12 )^k}\,{\langle f,f\rangle^{-1}}
 =
 \sum _{  f\in H_n(q)  } \, R_k\left(n,q\right) 
\,{\langle f,f\rangle^{-1}}
 (1 +O(nq)^{-\frac12 + \varepsilon})
\end{equation}
as $nq\to \infty$, where $R_k(n,q)$ is 
given by the $k$-fold residue
\begin{equation}\label{eqn:Rkintegral}
R_k
= \frac{(-1)^{k(k-1)/2}2^{k-1}}{k!} \frac{1}{(2\pi
i)^{k}}\oint \cdots \oint \frac{G(z_1,
\dots,z_{k})\Delta(z_1^2,\dots,z_{k}^2)^2} {\prod_{j=1}^{k}
z_j^{2k-1}}
\prod_{j=1}^k X_{n,q}(\tfrac12+z_j)^{-\frac12} dz_1\dots dz_{k},
\end{equation}
where
\begin{equation}
G(z_1,\dots,z_k)=A_k(z_1,\dots,z_k) 
\prod_{1\le i< j\le k}\zeta(1+z_i+z_j)
\end{equation}
and $A_k$ is the Euler product which is absolutely convergent for
$|\Re z_j|<\frac12 $, with $j=1,\ldots,k$, defined by
\begin{eqnarray} A_k(z_1,\dots,z_k) & =  &\prod_{p\nmid q} \prod_{1\le i < j \le k}
\left( 1-\frac{1}{p^{1+z_i + z_j}} \right) \nonumber \\&   &\quad\quad\times\frac{2}{\pi}\int_0^\pi
\sin^2\theta\prod_{j=1}^k\frac{e^{i\theta}\left(1-\frac{e^{i\theta}}
{p^{\frac12+ z_j}}\right)^{-1}-
e^{-i\theta}\left(1-\frac{e^{-i\theta}} {p^{\frac12+
z_j}}\right)^{-1}}{e^{i\theta}-e^{-i\theta}}\,d\theta .
 \end{eqnarray}
\end{conjecture}

To state the more general version of Conjecture~\ref{thm:sumLfkconjecture}, involving a
sum of products of $L_f(\frac12 + u_j)$, it is natural also to
consider the sums over even~$f$ and odd~$f$ separately. See
Conjectures~\ref{thm:Zfsumall} and~\ref{thm:Zfsumevenodd}.

\begin{theorem} \label{thm:SO2Nk}

 In the notation of Section~\ref{sec:shiftedmoments} we have
\begin{eqnarray}
J_k(SO(2N),0) =\left( 2^{k(k+1)/2} \prod_{j=1}^{k-1} \frac{j!}
{(2j)!} \right) \prod_{0\leq i<j \leq k-1} (N+\tfrac{i+j}{2}). \end{eqnarray}
More generally, with
\begin{equation}
G(z_1,\dots,z_k)= \prod_{1\leq \ell< m \leq
k} (1-e^{-z_{m}-z_{\ell}})^{-1}
\end{equation}
we have
\begin{eqnarray}
J_k(SO(2N),\alpha) & =  &\frac{(-1)^{k(k-1)/2}2^k}{k!}
\frac{1}{(2\pi i)^{k}}
\nonumber \\&   &\ \ \times 
 \oint\cdots \oint
\frac{
G(z_1,\dots,z_k)
\Delta(z_1^2,\ldots,z_{k}^2)^2 \prod_{j=1}^{k} z_j} {
\displaystyle \prod_{i=1}^{k} \prod_{j=1}^{k}
(z_j-\alpha_i)(z_j+\alpha_i)} 
e^{N \sum_{j=1}^k z_j} 
dz_1\cdots dz_{k}   \end{eqnarray}
 and
\begin{eqnarray}
J_k(O^-(2N),\alpha) & =  &i^{-k}\frac{(-1)^{k(k-1)/2}2^k}{k!}
\frac{1}{(2\pi i)^{k}}
\nonumber \\&   &\ \ \times 
 \oint\cdots \oint
\frac{
G(z_1,\dots,z_k)
\Delta(z_1^2,\ldots,z_{k}^2)^2 \prod_{j=1}^{k} \alpha_j} {
\displaystyle \prod_{i=1}^{k} \prod_{j=1}^{k}
(z_j-\alpha_i)(z_j+\alpha_i)} 
e^{N \sum_{j=1}^k z_j} 
dz_1\cdots dz_{k} . 
\end{eqnarray}
\end{theorem}

Comments:
\begin{enumerate}

\item $R_k(n,q)$ does not actually depend on $f\in H_n(q)$.  We write
(\ref{eqn:Hnqconjecture}) in this manner to stress that $R_k(n,q)$
is the expected value of $L_f(\frac12)^k$.  

\item To compare Theorem~\ref{thm:SO2Nk} and
 Conjecture~\ref{thm:sumLfkconjecture}, equating
conductors gives the equivalence
\begin{eqnarray}\label{eqn:Hnqconductor}
2N&=&cond(n,q):=\log ({q}/{4\pi^2}) + (\Gamma'/\Gamma)(n/2) \cr
&=& \log ({qn}/{8\pi^2}) +  O(n^{-1}).  
\end{eqnarray}
One can express the conjectured mean value in terms of the conductor in the
following way.  In (\ref{eqn:Rkintegral}) we can write
\begin{equation}
X_{n,q}(\tfrac12+z_j)^{-\frac12} = e^{\frac12 cond(n,q)\cdot z_j} {\mathcal G}(z_j),
\end{equation}
where ${\mathcal G}(z_j)=1+O(z^3)$.  As in Conjecture~\ref{thm:Lhalfchidconjecture},
we can express
$R_k(n,q)$ as a polynomial in the conductor, the first 3 terms of which
do not depend on the $X_{n,q}$ factor in the functional equation.

\item All of our conjectures naturally contain a factor of the form
$\prod X(\frac12\pm z_j)^{-\frac12}$, it just happens that in some
cases $X(\frac12\pm z_j)$ can be closely approximated by a simple function
of the conductor.  It is interesting that this same factor occurs in 
all of the
random matrix moments. In that case $X(s)=s^{-M}$, where $M=N$ or~$2 N$,
so in the formula for the
moments there occurs
$\prod X(e^{\pm z_j})^{- \frac12}=e^{\frac{M}{2}\sum \pm z_j}$.

\end{enumerate}


\subsection{  2nd Moment of the Riemann zeta-function
}\label{sec:secondmoment}

Now we consider   the 2nd moment of the Riemann zeta-function in
detail, putting our results in the context of the literature.

Ingham's result \cite{I} on the 2nd moment can be stated as
\begin{equation}
\int_0^T \zeta(s+\alpha)\zeta(1-s-\beta)\,dt=
\int_0^T (\zeta(1+\alpha-\beta)+\tau^{\beta -\alpha }
\zeta(1+\beta-\alpha))(1+O(t^{-\frac 12+\varepsilon}))\,dt
\end{equation}
where $s=\frac12+it$ and $\tau=\tau(t)=\frac{|t|}{2\pi}$;
this is valid for $|\alpha|, |\beta|<\frac12$. If we let $\alpha$ and
$\beta$ approach 0 here we obtain Ingham's theorem
\begin{equation}
\int_0^T|\zeta(\tfrac12+it)|^2\,dt  = \int_0^T (\log\tfrac t {2\pi}
+2\gamma)\,dt+O(T^{\frac 12 +\varepsilon}).
\end{equation}
Our conjecture is compatible with these results, because, when
$k=1$, the function $G(\alpha_1,\alpha_2)$ that appears in
Conjecture~\ref{thm:zeta2kconjecture} equals $\zeta(1+\alpha_1-\alpha_2)$. Computing
the residue, we find
\begin{eqnarray}
    P_1(x)&  = & \frac{1}{4\pi^2}
             \oint \oint 
\frac{\zeta(1+z_1-z_2)(z_2-z_1)^2 }
                              {z_1^2 z_2^2} 
e^{\frac{x}{2}(z_1-z_2)}
dz_1 dz_2 \nonumber \\&  = & x +2\gamma.
\end{eqnarray}
The second moment with a different weighting is now given; this
theorem is a slight variation of the theorem of Kober presented in
Titchmarsh \cite{T} and was inspired by the numerical calculations
described in Section~\ref{sec:riemannzetacalcs}.

\begin{theorem}
 Let
\begin{equation}I(\alpha,\beta,\delta)=\int_0^\infty  
\zeta(\tfrac12+it+\alpha)\zeta(\tfrac12-it-\beta)
e^{-\delta t}\,dt.\end{equation}
Then, for any $\eta>0$, $|  \alpha|, | \beta|\le \frac12-\eta $ and $|\arg
\delta|\le \frac \pi 2 -\eta$ , we have
\begin{equation}I(\alpha,\beta,\delta)=\int_0^\infty \left(\zeta(1+\alpha-\beta)+
\left(\frac t{2\pi}\right)^{\alpha-\beta}\zeta(1-\alpha+\beta)\right)
e^{-\delta t}\,dt+C_\delta(\alpha,\beta)+O(\delta\log 1/\delta)
\end{equation}
uniformly in $\alpha$, $\beta$, and $\delta$ where $C_\delta(\alpha,\beta)\ll
\log 1/\delta$ uniformly in $\alpha$ and $\beta$  and where
$C_\delta(\alpha,-\alpha)= -2\pi\zeta(2\alpha)$.
\end{theorem}

We restate the case $\beta=-\alpha$ as
\begin{corollary}
 For any fixed $\alpha$ with $|\alpha|< \frac12 $, we have
\begin{eqnarray}
\lim_{\delta\to 0}\left(\int_0^\infty |\zeta(\tfrac 12
+\alpha+it)|^2 e^{-\delta t} \,dt\right.
&-&\int_0^\infty
\left.\left(\zeta(1+2\alpha)+\left(\frac t{2\pi}\right)^{2\alpha}
\zeta(1-2\alpha)\right)e^{-\delta t}\,dt\right)\nonumber \\
&   & =-\pi \zeta(2\alpha)-(2\pi)^{2\alpha}\zeta(1-2\alpha)\Gamma(1-2\alpha)
\sin \pi \alpha \nonumber\\
&   & =-2\pi\zeta(2\alpha)  .
\end{eqnarray}
\end{corollary}

Note that
\begin{equation}\lim_{x\to 0}
(\zeta(1+x)+w^{-x}\Gamma(1-x)\zeta(1-x))=(2\gamma+\log w) .
\end{equation}
 Thus,  letting $\alpha \to 0$, gives
\begin{corollary} \label{thm:pimoment}

\begin{equation}
\lim_{\delta\to 0}\left( \int_0^\infty |\zeta(\tfrac 12  +it)|^2
e^{-\delta t} \,dt - \int_0^\infty
\left(2\gamma+\log\left(\frac{t}{2\pi}\right)\right)e^{-\delta t}
\,dt\right)=\pi .\end{equation}
\end{corollary}

Remark:  We discovered this Corollary after seeing the numerical
results of Section~\ref{sec:riemannzetacalcs}.  This result also follows from a result of
Hafner and Ivi\'c \cite{HI}.


\subsection{  4th Moment of the Riemann zeta-function
}\label{sec:fourthmoment}

Now we consider the 4th moment of the Riemann zeta-function in
detail.  Our discussion here builds upon work of Atkinson\cite{A},
Heath-Brown\cite{H-B1}, Conrey\cite{C}, and Motohashi\cite{Mot}.

Examining Motohashi's results in detail, 
consider
       \begin{equation}
           \label{eq:motohashi}
           \int_{-\infty}^\infty  \zeta(s+u_1) \dots
           \zeta(s+u_k)\zeta(1-s+v_1)\dots \zeta(1-s+v_k) g(t)\,dt
       \end{equation}
for a function $g(t)$ which is analytic in a horizontal strip $|\Im(t)|<c$
and decays sufficiently rapidly. Motohashi obtains an exact formula for
these moments for $k=1$ and $k=2$.
We reformulate Motohashi's theorem ($k=2$) in our context. Let
\begin{equation} C(v)=\frac{(2\pi)^v}{2 \cos \frac \pi 2 v}
\end{equation}
and let
\begin{equation}G_s(u,v)=\frac{\Gamma(s-u)}{\Gamma(s-v)}.
\end{equation}

Then, in a notation analogous to Motohashi's, the $k=2$ case
        of~(\ref{eq:motohashi}) equals
        \begin{equation} L_r+L_d+L_c+L_h  ,
        \end{equation}
where $L_r$ is the (residual) main term which we are interested in here:
\begin{equation}L_r({u},{v})= \int_{-\infty}^\infty
W(t,{u},{v})g(t)\,dt  ,
\end{equation}
with
\begin{eqnarray}
W(t,{u},{v})& =  &C(0)(G_s(0,0)+G_{1-s}(0,0))Z(u_1,u_2,v_1,v_2)\nonumber \\&   &\mbox{} + C(u_1+v_1)(G_s(u_1,v_1)+G_{1-s}(u_1,v_1))Z(-v_1,u_2,-u_1,v_2)\nonumber \\&   &\mbox{} + C(u_1+v_2)(G_s(u_1,v_2)+G_{1-s}(u_1,v_2))Z(-v_2,u_2,v_1,-u_1)\nonumber \\&   &\mbox{} + C(u_2+v_1)(G_s(u_2,v_1)+G_{1-s}(u_2,v_1))Z(u_1,-v_1,-u_2,v_2)\nonumber \\&   &\mbox{} + C(u_2+v_2)(G_s(u_2,v_2)+G_{1-s}(u_2,v_2))Z(u_1,-v_2,v_1,-u_2)\nonumber \\&   &\mbox{} +
C(u_1+u_2+v_1+v_2)(G_s(u_1,v_1)G_s(u_2,v_2)+G_{1-s}(u_1,v_1)G_{1-s}(u_2,v_2))\nonumber \\&   &\qquad \times Z(-v_1,-v_2,-u_1,-u_2),
 \end{eqnarray}
where $s=\frac12+it$, and
\begin{equation}
Z(u_1,u_2,v_1,v_2)=\frac{\zeta(1+u_1+v_1)\zeta(1+u_1+v_2)\zeta(1+u_2+v_1)
\zeta(1+u_2+v_2)}{\zeta(2+u_1+u_2+v_2+v_2)}.
\end{equation}
This formula may be obtained from Motohashi's
work (\cite{Mot} pp. 174 - 178) by a
careful analysis of his terms together with appropriate use of the
functional equation in the form
\begin{equation}\Gamma(s) \zeta(s) =\frac{(2\pi)^s}{2\cos \pi s/2}\zeta(1-s)
\end{equation}
and some trigonometric identities.

If we use the approximation
\begin{equation}\frac{\Gamma(s+\alpha)}{\Gamma(s+\beta)}=
(i|s|)^{\alpha-\beta}(1+O(1/|s|),
\end{equation}
we have, using $\tau=|t|/(2\pi)$,
\begin{equation}C(u)(G_s(u)+G_{1-s}(u))
= \tau^{-u} (1+O(1/\tau))
\end{equation} and
\begin{equation}C
(u+v)(G_s(u)G_s(v)+G_{1-s}(u)G_{1-s}(v))  = \tau^{-u-v} (1+O(1/\tau)).
\end{equation}
We then have
\begin{eqnarray}
W(t,{u},{v})& =  &\Bigl(Z(u_1,u_2,v_1,v_2)
           + \tau^{-u_1-v_1}Z(-v_1,u_2,-u_1,v_2)
           + \tau^{-u_1-v_2}Z(-v_2,u_2,v_1,-u_1) \nonumber \\
&   &\mbox{} + \tau^{-u_2-v_1}Z(u_1,-v_1,-u_2,v_2)
           + \tau^{-u_2-v_2}Z(u_1,-v_2,v_1,-u_2)\nonumber \\
&   &\mbox{} + \tau^{-u_1-u_2-v_1-v_2}Z(-v_1,-v_2,-u_1,-u_2)\Bigr)(1+O(
1/\tau)) .
\end{eqnarray}
By the formulae in Sections~\ref{sec:zetameanvalues} and~\ref{sec:concisesums}, the above agrees with the
$k=2$ case of Conjecture~\ref{thm:zeta2kconjecture}.

A residue
computation shows that our conjecture can be restated as
\begin{eqnarray} \int_0^T |\zeta(\tfrac12+it)|^4\,dt & =  &\int_0^T
\frac12\frac{1}{(2\pi i)^2} \oint \oint
\left(\frac{t}{2\pi}\right)^{x+y}
\nonumber \\&   &\ \ \ \times \frac{\zeta(1+x)^4\zeta(1+y)^4}
{\zeta(1+x-y)\zeta(1+y-x)\zeta(2+2x+2y)}\,dx\,dy\,dt
+O(T^{\frac12+\varepsilon}),
\end{eqnarray}
where we integrate around small circles centered on the origin.
This is in contrast to Conjecture~\ref{thm:zeta2kconjecture},
which when $k=2$ expresses the formula in terms of four contour integrals.
It may be that our formulae can be similarly simplified for
all~$k$, but we have not succeeded in doing so.


\section{ Moments in $t$-aspect 
}\label{sec:taspectrecipe}

The principle behind our method of conjecturing mean values is
that the Dirichlet series coefficients of $L$-functions have
an approximate orthogonality relation when averaged over a family.
These orthogonality relations are used to identify the main terms
in the mean values.

In this section we give a detailed account of the case of 
moments of a single primitive $L$-function. 
We describe the recipe for conjecturing the mean values,
applying it first to the case of the Riemann $\zeta$-function,
and then to a general primitive $L$-function.
In the remainder of this section we manipulate the formulas into
a more usable form, and also obtain a generalization of
Conjecture~\ref{thm:zeta2kconjecture}.
Later in Section~\ref{sec:families}
we recast our principles in a more general setting and consider
the averages of various families of $L$-functions.


\subsection{ The recipe
}\label{sec:therecipe}
The following is our recipe for conjecturing the $2k$th moment
of an $L$-function:

\begin{enumerate}
\item{}Start with a product of $2k$ shifted $L$-functions: 
\begin{equation}\label{eqn:firstZproduct}
Z(s,\alpha_1,\dots,\alpha_{2k}) =
Z(\tfrac 12 + \alpha_1)\cdots Z(\tfrac 12 + \alpha_{2k})
\end{equation}
(here we have written the $Z$-function, but the
examples below will show that the method applies to either
the $L$- or the $Z$-function).

\item{}Replace each $L$-function with the two terms from its
approximate functional equation, ignoring the remainder term.  
Multiply out the resulting
expression to obtain $2^{2k}$ terms.

\item{}Keep the $\binom{2k}{k}$ terms for which the product of $\chi$-factors from the
functional equation is not rapidly oscillating.
Use (\ref{eqn:twochiproduct}) to simplify the nonoscillating $\chi$-factors.

\item{}In each of those $\binom{2k}{k}$ terms, keep only the diagonal
from the sum.

\item{}Complete the resulting sums, and call the total 
$M(s,\alpha_1,\dots,\alpha_{2k})$.

\item{}The conjecture is
\begin{equation}
\int_{-\infty}^\infty Z(\tfrac12+it,\alpha_1,\dots,\alpha_{2k}) g(t) \,dt
=
\int_{-\infty}^\infty M(\tfrac12+it,\alpha_1,\dots,\alpha_{2k})
(1+O(t^{-\frac12 + \varepsilon})) g(t) \,dt ,
\end{equation}
for all $\varepsilon>0$, where $g$ is a suitable weight function.
In other words, $L(s,\alpha)$ and $M(s,\alpha)$ have the same 
expected value if averaged over a sufficiently large range.
\end{enumerate}

Notes:
\begin{enumerate}
\item In order to see the structure of these mean values, it is necessary 
to include the shifts~$\alpha_j$.  
One can obtain the moments of $L(\frac12 +it)$ by allowing the
shifts to tend to~$0$.  Because of the shifts~$\alpha_j$ we avoid
higher-order poles in our expressions.

\item The recipe applies to either the $L$-function or the $Z$-function,
and we give examples of both cases.  The $Z$-function case can be directly
obtained from the $L$-function, although the reverse is not true in general.

\item For the approximate functional equations in  the recipe, one can ignore
the range of summation because it will just be extended to infinity
in the final step.

\item We do not define what is meant by a ``suitable weight function'',
but it is acceptable to take $g(t)=g_T(t)=f(t/T)$ for a fixed
integrable function~$f$.  In particular, one can take 
$f$ to be the characteristic function of the interval~$[0,1]$,
obtaining the mean value $\int_0^T Z(\tfrac12+it,\alpha)dt$.
From this one can recover a fairly general weighted integral
by partial integration.

\item The error term $O(t^{-\frac12 + \varepsilon})$ fits
with known examples and numerical evidence.  See Section~\ref{sec:numerics}.

\item The above procedure is a recipe for conjecturing 
all of the main terms in the mean value of an $L$-function.  
It is not a heuristic, and the steps cannot be justified.
In particular, some steps can throw away terms which are the same size 
as the main term, and other steps add main terms back in.  Our 
conjecture is that all of those errors cancel.

\end{enumerate}


\subsection{  Moments of the Riemann $\zeta$-function
}\label{sec:zetameanvalues}

We illustrate our recipe in the case of the Riemann zeta-function.
In this section we consider the most familiar case of
moments of $\zeta(\frac12+it)$. 
In Section~\ref{sec:Lmeanvalues} we 
relate this to moments of
$Z(\frac12+it)$ and repeat the
calculation for the $Z$-function of an arbitrary primitive
$L$-function.

Consider
\begin{equation}\label{eqn:zetaproduct}
Z(s,\alpha )= \zeta(s+\alpha _1)\cdots \zeta(s+\alpha
_k)\zeta(1-s-\alpha _{k+1})\cdots \zeta(1-s-\alpha _{2k}),
\end{equation}
where $\alpha =(\alpha _1,\ldots, \alpha _{2k})$.
Note that this is slightly different notation than given 
in~(\ref{eqn:firstZproduct}).
Our goal is a formula for
\begin{equation}
\int_{-\infty}^\infty Z(\tfrac12+it,\alpha ) g(t)\, dt.
\end{equation}

For each $\zeta$-function we use the approximate functional equation
\begin{equation}
\zeta(s)=\sum_m \frac{1}{m^s}+\chi(s)\sum_n \frac {1}{n^{1-s}} + remainder.
\end{equation}
Recall that we ignore the remainder term and the limits on the sums.
Multiplying out the resulting expression we obtain $2^{2k}$
terms, and the recipe tells us to keep those terms in which the
product of $\chi$-factors is not oscillating rapidly.  

If $s=z+it$ with $z$ bounded (but not necessarily real) then
\begin{equation}
\chi(s)=\(\frac{t}{2\pi}\)^{\frac12 -s}e^{it+\pi i/4}(1+O(\tfrac{1}{t}))
\label{eqn:chis}
\end{equation}
and
\begin{equation}
\chi(1-s)=\(\frac{t}{2\pi}\)^{s-\frac12 }e^{-it-\pi i/4}(1+O(\tfrac{1}{t})),
\label{eqn:chi1minuss}
\end{equation}
as $t\to+\infty$.
We use the above formulas to determine which products of
$\chi(s)$ and $\chi(1-s)$ are oscillating.

One term which does not have an oscillating factor is the one where we 
use the ``first part'' of each approximate functional equation, 
for it does not have any $\chi$-factors.  With $s=\frac12+it$, 
that term~is
\begin{eqnarray} 
&& \sum _{\ontop{  m_1,\dots,m_k}{n_1,\dots,n_k }}
m_1^{-s -\alpha _1}\cdots m_k^{-s -\alpha _k}n_1^{s-1 +\alpha
_{k+1}}\cdots n_k^{s-1 +\alpha _{2k}}\cr
 && \phantom{XXXXXXXX}
=\sum _{\ontop{  m_1,\dots,m_k}{n_1,\dots,n_k }}
m_1^{-\frac12 -\alpha _1}\cdots m_k^{-\frac12 -\alpha _k}n_1^{-\frac12 +\alpha
_{k+1}}\cdots n_k^{-\frac12 +\alpha _{2k}}
\(\frac{n_1\cdots n_k }{  m_1 \cdots m_k}\)^{it}.
\end{eqnarray}
According to the recipe we keep the diagonal from the above sum, which is
\begin{equation}  
\sum_{m_1 \cdots m_k=n_1\cdots n_k}
m_1^{-\frac12 -\alpha _1}\cdots m_k^{-\frac12 -\alpha _k}n_1^{-\frac12 +\alpha
_{k+1}}\cdots n_k^{-\frac12 +\alpha _{2k}}.
\end{equation}
If we define
\begin{equation}\label{eqn:Rs}
R(s; {\alpha })=\sum_{m_1\cdots m_k =n_1 \cdots n_k  } 
\frac {1}{m_1^{s+\alpha
_1}\cdots m_k^{s+\alpha _k}n_1^{s-\alpha _{k+1}}\cdots n_k^{s-\alpha
_{2k}}},
\end{equation}
where the sum is over all positive $m_1,\dots,m_k, n_1, \dots,n_k$ such that
$m_1\cdots m_k =n_1 \cdots n_k$, 
then $R(\frac12; \alpha)$ is the first piece which we have identified as
contributing to the mean value.  (The sum in equation~(\ref{eqn:Rs})
does not converge for $s=\frac12$.  See Theorem~\ref{thm:Rkgeneral} for 
its analytic continuation.) 

Note that the variable $s$ in equation~(\ref{eqn:Rs}) should not
be viewed the same as the variable $s=\frac12+it$ from the previous
equations.
We are employing a trick of beginning with an expression involving
$s$ and $1-s$, noting that we will later be setting $s=\frac12$,
so instead we consider an expression only involving~$s$,
which later will be set equal to~$\frac12$.
This same trick will appear in Section~\ref{sec:generalrecipe}
when we consider more general mean values.

Now consider one of the other terms, say the one where we use the
second part of the approximate functional equation from 
$\zeta(s+\alpha _1)$ and the second part from 
$\zeta(1-s-\alpha _{k+1})$. 
By (\ref{eqn:chis}) and (\ref{eqn:chi1minuss}),
\begin{equation}\label{eqn:twochiproduct}
\chi(s+\alpha _1)\chi(1-s-\alpha _{k+1}) \sim
\left(\frac{t}{2\pi}\right)^{-\alpha _1+\alpha _{k+1}},
\end{equation}
which is not rapidly oscillating.
Using this and proceeding as above, the contribution from this term will
be
\begin{equation}
\(\frac{t}{2\pi}\)^{-\alpha _1+\alpha _{k+1}}
R(\tfrac{1}{2};
\alpha _{k+1},
\alpha _2,\dots,\alpha _k,
\alpha _1,\alpha
_{k+2},\dots,\alpha _{2k}).
\end{equation}
More generally, note that
\begin{equation}
\chi(s+\beta _1)\cdots \chi(s+\beta _J)
\chi(1-s-\gamma _{1}) \cdots \chi(1-s-\gamma _{K}) \sim
\(\frac{t}{2\pi e}\)^{-i(J-K)t} 
e^{i(J-K)\pi/4} 
\left(\frac{t}{2\pi}\right)^{-\sum\beta_j+ \sum\gamma_{j}},
\end{equation}
which is rapidly oscillating (because of the $\,it\,$ in the exponent) unless $J=K$.
Thus, the recipe tell us to keep those terms which involve
an equal number of $\chi(s+\alpha_j)$ and $\chi(1-s-\alpha_{k+j})$
factors.  This gives  a total of
$\binom{2k}{k}= \sum_{j=0}^k \binom{k}{j}^2$ terms in the final answer.

We now describe a typical term of the   conjectural formula. 
First note that the function
$ R(s; \alpha_1, \dots,\alpha _k, \alpha_{k+1},\dots,\alpha _{2k})$
is symmetric in $\alpha_1, \dots,\alpha _k$ and in 
$\alpha_{k+1},\dots,\alpha _{2k}$, so we can rearrange the entries so
that the first $k$ are in increasing order, as are the
last~$k$.  Thus, the final result will be a sum of terms indexed
by the $\binom{2k}{k}$ permutations
$\sigma\in S_{2k}$ such that $\sigma(1)<\cdots<\sigma(k)$ and
$\sigma(k+1)<\cdots< \sigma(2k)$.
We denote the set of such permutations by~$\Xi$.
Second, note that  
the product of an equal number of $\chi(s+\alpha _j)$ and $\chi(1-s-\alpha _{k+j})$,
as in (\ref{eqn:twochiproduct}), can be written as
\begin{equation}
\left( \frac{t}{2\pi}
\right)^{\frac{1}{2}(-\alpha _1-\cdots
 -\alpha _k+\alpha _{k+1}+\cdots +\alpha _{2k})}
\left( \frac{t}{2\pi}\right)^{\frac{1}{2}
(\alpha _{\sigma(1)}+ \cdots + \alpha _{\sigma(k)} -\alpha
      _{\sigma(k+1)}-\cdots -\alpha _{\sigma(2k)})} .
\end{equation}
For example, (\ref{eqn:twochiproduct}) is the case
$\sigma(1)=k+1$, $\sigma(k+1)=1$, and $\sigma(j)=j$ otherwise. 

If we set
\begin{equation}\label{eqn:Wdef}
W(z,\alpha ,\sigma)=\left( \frac{y}{2\pi}\right)^{\frac{1}{2}
(\alpha _{\sigma(1)}+ \cdots + \alpha _{\sigma(k)} -\alpha
      _{\sigma(k+1)}-\cdots -\alpha _{\sigma(2k)})} 
R(x; \alpha _{\sigma(1)},\ldots,\alpha _{\sigma(2k)}),
 \end{equation}
for $z=x+iy$ with $x$ and $y$ real,
then combining all terms we have 
\begin{equation}
M(z;\alpha):=
\left( \frac{y}{2\pi}
\right)^{\frac{1}{2}(-\alpha _1-\cdots
 -\alpha _k+\alpha _{k+1}+\cdots +\alpha _{2k})} \sum_{\sigma \in \Xi}
 W(z,\alpha ,\sigma). 
\end{equation}
The recipe has produced the conjecture
\begin{equation}\label{conj:M}  
\int_{-\infty}^\infty Z(\tfrac12+it,\alpha)g(t)\,dt =
\int_{-\infty}^\infty M(\tfrac12+it,\alpha)(1+O(t^{-\frac12 + \varepsilon}))g(t)\,dt ,
\end{equation}
with $Z(s,\alpha)$ given in (\ref{eqn:zetaproduct}) and $M(s;\alpha)$ given above.

Note that the exponent of $(t/2\pi)$ in (\ref{eqn:Wdef}) has half the
$\alpha_j$ with $+$~sign and the other half with $-$~sign, and the same
holds for $R(s,\alpha)$.  This allows an alternate 
interpretation of
$\Xi$ as the set of ways  of choosing $k$ elements 
from $\{\alpha_1,\dots,\alpha_{2k}\}$.

The general case of Conjecture~\ref{thm:zeta2kconjecture} is stated in terms of
the $Z$-function.   We can recover the mean value of the $Z$-function
directly from that of the $L$-function (in this case, the $\zeta$-function). 
By the functional equation and (\ref{eqn:chis}) we see that
\begin{eqnarray}
Z(s+\alpha_1)\cdots
Z(s+\alpha_{2k})
& =& \left( \frac{t}{2\pi}
\right)^{\frac{1}{2}(\alpha _1+\cdots
 +\alpha _k-\alpha _{k+1}-\cdots -\alpha _{2k})} (1+O(1/t)) \cr
&& \ \ \Bigl. \times \zeta(s+\alpha _1)\cdots \zeta(s+\alpha
_k)\zeta(1-s-\alpha _{k+1})\cdots \zeta(1-s-\alpha _{2k}).
\end{eqnarray}
The factor $\displaystyle  
\left( \frac{t}{2\pi}
\right)^{\frac{1}{2}(\alpha _1+\cdots
 +\alpha _k-\alpha _{k+1}-\cdots -\alpha _{2k})}
$
can be absorbed into the weight function $g(t)$, so we obtain the conjecture
\begin{equation}\label{eqn:firstZmeanvalue}
\int_{-\infty}^\infty Z(s+\alpha_1)\cdots
Z(s+\alpha_{2k})
g(t)\, dt =
\int_{-\infty}^\infty
\sum_{\sigma \in \Xi}
 W(s,\alpha ,\sigma)(1+O(t^{-\frac12 + \varepsilon}))g(t)\,dt ,
\end{equation}
where $s=\frac12+it$.

In the next subsection we directly obtain the above conjecture for 
the $Z$-function of a general primitive
$L$-function, and in the remainder of this section we perform various manipulations to 
put these in the form
of Conjecture~\ref{thm:zeta2kconjecture}.


\subsection{  Moments of a primitive $L$-function
}\label{sec:Lmeanvalues}

Consider the primitive $L$-function
\begin{equation}\mathcal{L}(s)=\sum_{n=1}^\infty \frac{a_n}{n^s}=
\prod_p \mathcal{L}_p\left(\frac1{p^s}\right), \qquad (\sigma>1) .
\end{equation}
  We assume a functional equation of the special form 
$\xi_L(s)=\gamma_L(s)L(s)=\varepsilon \overline{\xi}_L(1-s)$, where
\begin{equation}
\gamma_L(s)= Q^s \prod_{j=1}^w \Gamma(\tfrac12 s+ \mu_j) ,
\end{equation}
with $\{\mu_j\}$ stable under complex conjugation.
Note that we have $w_j=\frac12$, which is conjectured to hold for 
arithmetic $L$-functions.
We also assume
\begin{eqnarray}\mathcal{L}_p(x)&=&\sum_{n=0}^\infty a_{p^n}x^n \cr
&=&\prod_{j=1}^w (1-\gamma_{p,j} x)^{-1},
\end{eqnarray}
where $w$ is the degree of $\mathcal{L}$ and where
$|\gamma_{p,j}|=0$ or 1.  Again this is conjectured to hold for
arithmetic $L$-functions. 

We are going to evaluate the moments
of the $Z$-function $Z_\mathcal{L}(s) = 
\varepsilon^{-\frac12}\mathcal{X}(s)^{-\frac12} \mathcal{L}(s)$,
where
\begin{eqnarray}
\mathcal{X}(s) &=&
\frac{\overline{\gamma_{\mathcal{L}}(1-s)}}{\gamma_{\mathcal{L}}(s)} \cr
&=& Q^{1-2s} \prod_{j=1}^w \frac{ \Gamma(\tfrac12 (1-s)+ \overline{\mu_j})}
{ \Gamma(\tfrac12 s+ \mu_j)} .
\end{eqnarray}
We will have to determine when products of $\mathcal{X}(s)$ and
$\mathcal{X}(1-s)$ are not rapidly oscillating.  By Stirling's formula
\begin{equation}
\Gamma(\sigma+it) =
e^{-\frac{\pi t}{2}}t^{\sigma-\frac12}\(\frac {t}{e}\)^{it}
e^{\frac{i\pi}{2}(\sigma-\frac12)}
\( 1-\frac{i}{t}\(\frac{1}{12}-\frac{\sigma}{2}+\frac{\sigma^2}{2}\)+
O\(\frac{1}{t^2}\)\)
\end{equation}
we obtain
\begin{equation}\label{eqn:Xs}
\mathcal{X}(s)=Q^{1-2s} \(\frac{t}{2}\)^{w(\frac12 -s)}
\(\frac{t}{2\pi e}\)^{-\sum \Im(\mu_j)}
e^{w(it+i \pi /4)}(1+O(\tfrac{1}{t})),
\end{equation}
as $t\to +\infty$.
Note that the above expression can be simplified
because we have assumed $\sum \Im(\mu_j)=0$.

Now we are ready to produce a conjecture for
\begin{equation}\label{eqn:generalIk}
I_k(\mathcal{L},\alpha_1,\dots, \alpha_{2k},g)=
\int_{-\infty}^{\infty}  Z_\mathcal{L}(s+\alpha_1)\cdots
Z_\mathcal{L}(s+\alpha_{2k})g(t)\,dt.
\end{equation}
with $s=\frac12 +it$.  

By the definition of $Z$,
\begin{equation}
\prod_{j=1}^{2k}Z_\mathcal{L}(s+\alpha_j)=
\prod_{j=1}^{2k}\varepsilon^{-\frac 12}
\mathcal{X}(s+\alpha_j)^{-\frac 12 } \prod_{j=1}^{2k}
\mathcal{L}(s+\alpha_j) .
\end{equation}
According to the recipe, we replace each $\mathcal{L}(s+it)$ by its 
approximate functional
equation and multiply out the product obtaining $2^{2k}$ terms.
A typical term is a product of $2k$ sums arising from either the first
piece or the second piece of the approximate functional equation.
Consider a term where we have $\ell$ factors from the first
piece of the approximate functional equation and $2k-\ell$ factors from the
second piece.  To take one specific example, suppose it is the  first $\ell$
factors where we choose the first piece of the approximate functional equation,
and the last $2k-\ell$ factors where we take the second piece: 
\begin{eqnarray} \label{eqn:generalfirstterm}
&   &  \varepsilon^{- k }
\mathcal{X}(\tfrac 12 +\alpha_1+it)^{-\frac 12}\cdots 
      \mathcal{X}(\tfrac 12 +\alpha_{\ell}+it)^{-\frac 12}
\sum_{n_1}\frac {a_{n_1}}{n_1^{\frac 12 +\alpha_1+it}}
            \cdots \sum_{n_\ell}\frac {a_{n_\ell}}{n_\ell^{\frac 12 +\alpha_\ell+it}} \nonumber \\
&   &\phantom{xxxxxxx} \times \varepsilon^{2k-\ell} 
  \mathcal{X}(\tfrac 12 +\alpha_{\ell+1}+it)^{\frac12}
\cdots \mathcal{X}(\tfrac 12 +\alpha_{2k}+it)^{\frac12}
  \sum_{n_{\ell+1}}\frac {\overline{a_{n_{\ell+1}}}}{n_{\ell+1}^{\frac 12
-\alpha_{\ell+1}-it}} \cdots \sum_{n_{2k}}\frac
{\overline{a_{n_{2k}}}}{n_{2k}^{\frac 12 -\alpha_{2k}-it}}  \cr
& & \phantom{XX} = 
\varepsilon^{ k -\ell}
\left(\frac{\mathcal{X}(\tfrac 12 +\alpha_1+it)\cdots \mathcal{X}(\frac
12 +\alpha_\ell+it)} {\mathcal{X}(\tfrac 12 +\alpha_{\ell+1}+it)\cdots
\mathcal{X}(\tfrac 12 +\alpha_{2k}+it)}\right)^{-\frac 12}
\nonumber \\
&   &\phantom{XXxxxxxxx} \times 
\sum_{n_1}\cdots \sum_{n_{2k}} 
        \frac{a_{n_1}\cdots a_{n_\ell}\overline{a_{n_{\ell+1}}\cdots 
           a_{n_k}}}
{n_1^{\frac12 +\alpha_1}\cdots n_\ell^{\frac12 +\alpha_\ell} 
n_{\ell+1}^{\frac12 -\alpha_{\ell+1}} \cdots n_{2k}^{\frac 12-\alpha_{2k}}}
 \left(\frac{n_1 \cdots n_\ell}{n_{\ell+1}\cdots n_{2k}}\right)^{-it}.
\end{eqnarray}
The recipe tells us to retain only the expressions of this sort where the 
factor involving  $\mathcal{X}$ is not oscillating.
By (\ref{eqn:Xs}) the requirement is that $\ell=k$
(and in particular $2k$ has to be even),
and we have
\begin{equation}
\left(\frac{\mathcal{X}(\tfrac 12 +\alpha_1+it)\dots \mathcal{X}(\tfrac 12
+\alpha_k+it)} {\mathcal{X}(\tfrac 12 +\alpha_{k+1}+it)\dots \mathcal{X}(\tfrac
12 +\alpha_{2k}+it)}\right)^{-\frac 12}
=
\left( \frac{Q^{\frac 2w}t}{2}
\right)^{\frac{w}{2}(\alpha _1+\cdots
 +\alpha _k-\alpha _{k+1}-\cdots -\alpha _{2k})}(1+O(\tfrac1t))
.
\end{equation}

Now the recipe tells us to keep the diagonal from the remaining sums,
which in
(\ref{eqn:generalfirstterm}) is the terms where
$n_1 \cdots n_\ell=n_{\ell+1}\cdots n_{2k}$.  So in the same way
as the $\zeta$-function case in the previous section we let
\begin{equation}\label{eqn:Rkgeneral}
R(s,\alpha)=\sum_{n_1 \cdots n_k=n_{k+1}\cdots n_{2k}} 
        \frac{a_{n_1}\cdots a_{n_k}\overline{a_{n_{k+1}}\cdots 
           a_{n_k}}}
{n_1^{s +\alpha_1}\cdots n_k^{s +\alpha_k} 
n_{k+1}^{s -\alpha_{k+1}} \cdots n_{2k}^{s -\alpha_{2k}}},
\end{equation}
and
\begin{equation}
W(z,\alpha ,\sigma)=\left( \frac{Q^{\frac 2w}y}{2}\right)^{\frac{w}{2}
(\alpha _{\sigma(1)}+ \cdots + \alpha _{\sigma(k)} -\alpha
      _{\sigma(k+1)}-\cdots -\alpha _{\sigma(2k)})} 
R(x; \alpha _{\sigma(1)},\ldots,\alpha _{\sigma(2k)}),
 \end{equation}
for $\sigma\in\Xi$, the set of permutation of $\{1,\dots,2k\}$
with  $\sigma(1)<\cdots<\sigma(k)$ and
$\sigma(k+1)<\cdots< \sigma(2k)$.  Then
\begin{equation}
M(z;\alpha)=
\sum_{\sigma \in \Xi}
 W(z,\alpha ,\sigma), 
\end{equation}
and we arrive at the conjecture
\begin{equation}  
I_k(\mathcal{L},\alpha_1,\dots, \alpha_{2k},g)
=
\int_{-\infty}^\infty M(\tfrac12+it,\alpha)(1+O(t^{-\frac12 + \varepsilon}))g(t)\,dt ,
\end{equation}
which is the same form as we obtained in~(\ref{eqn:firstZmeanvalue}).

We will now examine the expressions in these conjectures in 
detail, rewriting them in a more explicit form.  


\subsection{ The arithmetic factor in the conjectures
} \label{sec:thearithmeticfactor}

We retain the notation of the previous subsection.  In particular,
$\mathcal{L}(s)$ is a primitive $L$-function having
the properties listed at the beginning of Section~\ref{sec:Lmeanvalues}
and
$R(s;\alpha)$ is given in~(\ref{eqn:Rkgeneral}). 

\begin{theorem} \label{thm:Rkgeneral}
Suppose $|\alpha_j|<\delta$ for $j=1,\dots,2k$.  Then
$R(s;\alpha_1,\dots,\alpha_{2k})$ converges absolutely for
$\sigma>\frac12 + \delta$ and has a meromorphic continuation
to $\sigma>\frac14 + \delta$.  Furthermore,
\begin{equation}
 R(s;\alpha_1,\dots,\alpha_{2k})  = 
\prod_{i,j=1}^k\zeta(2s+\alpha_i-\alpha_{k+j})
\,A_k(s;\alpha_1,\dots,\alpha_{2k})
\end{equation}
where 
\begin{equation} \label{eqn:bigAk}
A_k(s;\alpha_1,\dots,\alpha_{2k})=\prod_p 
\(
\prod_{i,j=1}^k
(1-p^{-2s-\alpha_i+\alpha_{k+j}})
\)
B_p(s;\alpha_1,\dots,\alpha_{2k}) 
\end{equation}
with
\begin{equation}  \label{eqn:bigBk} 
B_p(s;\alpha_1,\dots,\alpha_{2k}) =
\int_0^1 \prod_{j=1}^k
\mathcal{L}_p\left(\frac{e(\theta)}{p^{s+\alpha_j}}\right)
\overline{\mathcal{L}_p}\left(\frac{e(-\theta)}{p^{s-\alpha_{k+j}}}\right)
\,d\theta.
\end{equation}
\end{theorem}

\begin{proof}
We assumed $|\gamma_{p,j}|\le 1$, so we have the Ramanujan bound 
$a_n\le d_w(n)\ll n^\varepsilon$.  That implies absolute convergence
of $R(s;\alpha)$ for $\sigma>\frac12 +\delta+\varepsilon$.

Since the coefficients of $R(s;\alpha)$ are multiplicative, as is
the condition $n_1\cdots n_k=n_{k+1}\cdots n_{2k}$, we can write
$R(s;\alpha)$ as an Euler product:
\begin{eqnarray}
 R(s;\alpha_1,\dots,\alpha_{2k})&  = &
\sum _{  n_1\cdots n_k=n_{k+1}\cdots n_{2k} }
\frac{a_{n_1}\cdots a_{n_k}\overline{a_{n_{k+1}}\cdots a_{n_{2k}}}}
{n_1^{s+\alpha_1}\cdots n_k^{s-\alpha_{2k}}}\nonumber \\
&  = &\prod_p\sum_{  \sum_{j=1}^k e_j=\sum_{j=1}^k e_{k+j} }
\frac{a_{p^{e_1}}\cdots a_{p^{e_k}}
\overline{a_{p^{e_{k+1}}}\cdots a_{p^{e_{2k}}}}}
{p^{e_1(s+\alpha_1)}\cdots p^{e_{2k}(s-\alpha_{2k})}}\nonumber \\
&  = &\prod_p\left(1+|a_p|^2\sum_{i,j=1}^k \frac{1}{p^{2s+\alpha_i-\alpha_{k+j}}}
+\sum_{j=2}^\infty \frac{c_{p^j}(\alpha_1,\dots,\alpha_{2k})}{p^{2js}}
+\cdots\right)\nonumber \\
&  = &\prod_{i,j=1}^k\zeta(2s+\alpha_i-\alpha_{k+j}) \\
& & \phantom{X}\times
\prod_p\left(1+(|a_p|^2-1)\sum_{i,j=1}^k \frac{1}{p^{2s+\alpha_i-\alpha_{k+j}}}
+\sum_{j=2}^\infty \frac{c'_{p^j}(\alpha_1,\dots,\alpha_{2k})}{p^{2js}}
+\cdots\right)\nonumber \\
&  = &\prod_{i,j=1}^k\zeta(2s+\alpha_i-\alpha_{k+j})
A_k(s;\alpha_1,\dots,\alpha_{2k}) ,
 \end{eqnarray}
say.  Above $c_{p^j}$ and $c'_{p^j}$ are just shorthand for the (complicated)
coefficients in the Euler product. Estimating  them trivially
and using the fact that $|a_p|^2$ is 1 on average, 
(which is conjectured to hold for primitive elements of the
Selberg class)
we find that 
$A_k(s;\alpha)$ is analytic in a neighborhood of $\sigma=\frac12$.

Finally, we have
\begin{equation}
A_k(s;\alpha_1,\dots,\alpha_{2k})=\prod_p \prod_{i,j=1}^k
(1-p^{-2s-\alpha_i+\alpha_{k+j}})
B_p(s;\alpha_1,\dots,\alpha_{2k}) 
\end{equation}
where
\begin{eqnarray} B_p(s;\alpha_1,\dots,\alpha_{2k})&  = &
\sum_{  \sum_{j=1}^k e_j=\sum_{j=1}^k e_{k+j} }
\frac{a_{p^{e_1}}\cdots a_{p^{e_k}}
\overline{a_{p^{e_{k+1}}}\cdots a_{p^{e_{2k}}}}}
{p^{e_1(s+\alpha_1)}\cdots p^{e_{2k}(s-\alpha_{2k})}}\nonumber \\
&  = &
\int_0^1 \sum_{  e_1,\dots,e_{2k} }
\frac{a_{p^{e_1}}\cdots a_{p^{e_k}}
\overline{a_{p^{e_{k+1}}}\cdots a_{p^{e_{2k}}}}}
{p^{e_1(s+\alpha_1)}\cdots p^{e_{2k}(s-\alpha_{2k})}}
\, e\(\Bigl(\sum_{j=1}^k e_j-\sum_{j=1}^k e_{k+j}\Bigr)\theta\)  d\theta\nonumber \\
&=&\int_0^1 \prod_{j=1}^k \sum_{e_j=0}^\infty
\frac{a_{p^{e_j}}}{p^{e_j(s+\alpha_j)}}e(e_j \theta)
 \prod_{j=1}^k
\frac{\overline{a_{p^{e_{k+j}}}}}{p^{e_{k+j}(s-\alpha_{k+j})}}e(-e_{k+j}\theta)\,d\theta\nonumber
\\
&  = &\int_0^1 \prod_{j=1}^k
\mathcal{L}_p\left(\frac{e(\theta)}{p^{s+\alpha_j}}\right)
\overline{\mathcal{L}_p}\left(\frac{e(-\theta)}{p^{s-\alpha_{k+j}}}\right)
d\theta,
\end{eqnarray}
as claimed.
\end{proof}

To summarize, the conjecture for the general mean value
$I_k(\mathcal{L},\alpha_1,\dots, \alpha_{2k},g)$
involves the function $M(s;\alpha)$, which can be written as
\begin{equation}
M(s;\alpha)=
\sum_{\sigma \in \Xi}
 W(s;\alpha_{\sigma(1)},\dots ,\alpha_{\sigma(2k)}),
\end{equation}
where we have written $W(s;\alpha_{\sigma(1)},\dots ,\alpha_{\sigma(2k)})$
for $W(s;\alpha)$.
And by Theorem~\ref{thm:Rkgeneral},
\begin{equation}
 W(s;\alpha_1,\dots,\alpha_{2k})  = 
\left( \frac{Q^{\frac 2w}t}{2}
\right)^{\frac{w}{2} \sum_{j=1}^k \alpha _j-\alpha _{k+j}}
A_k(s;\alpha_1,\dots,\alpha_{2k})
\prod_{i,j=1}^k\zeta(2s+\alpha_i-\alpha_{k+j}) .
\end{equation}

One can see the above elements  in Conjecture~\ref{thm:zeta2kconjecture},
in particular the form of $B_p(s;\alpha)$, 
and $A_k(z)$ in that conjecture equals $A_k(\frac12,\alpha)$ given above. 
The overall structure is slightly different because
Conjecture~\ref{thm:zeta2kconjecture} is expressed as a multiple contour
integral, as opposed to a sum over permutations.  In the next subsection
we show how to write the sum over permutations in a compact form.
In the following subsection we return to the functions $A_k$ and write
them in a more explicit form. 


\subsection{  Concise form of permutation sums
}\label{sec:concisesums}

As we have seen, our methods naturally lead to
an expression involving a sum over  permutations.
In this section we describe how to write those sums in a compact
form involving contour integrals.
Similar combinatorial sums arise from our matrix ensemble
calculations, and we have previously stated our main results
and conjectures in this compact form.

Note that in both of these Lemmas, the terms in the sum on the left
side have singularities.  However, examining the right side of the
formula makes it clear that those singularities all cancel.

\begin{lemma} \label{thm:concisesumunitary}
 Suppose
$F(a;b)=F(a_1,\dots,a_k;b_1,\dots,b_{k})$ is a function of $2k$
variables, which is symmetric with respect to the first $k$
variables and also symmetric with respect to the second set of $k$
variables. Suppose also that $F$ is regular near $(0,\dots,0)$.
Suppose further that $f(s)$ has a simple pole of residue~$1$ at
$s=0$ but is otherwise analytic in a neighborhood about $s=0$.
Let
\begin{equation}
K(a_1,\dots,a_k;b_1,\dots b_k)= F(a_1,\dots;\dots,b_k)
\prod_{i=1}^k\prod_{j=1}^k f(a_i-b_j).
\end{equation}
If for all $1\leq i,j \leq k$, $\alpha_i-\alpha_{k+j}$ is
contained in the region of analyticity of $f(s)$ then
\begin{eqnarray}
 \sum_{\sigma \in \Xi}&   &K( \alpha_{\sigma(1)},\dots,  \alpha_{\sigma(k)};
    \alpha_{\sigma(k+1)}\dots  \alpha_{\sigma(2k)})=\nonumber \\&   &\ \ \ \ \ \ \
 \frac{(-1)^k}{k!^2}\frac{1}{(2\pi i)^{2k}}
\oint \cdots \oint \frac{K(z_1, \dots,z_k;z_{k+1},\dots,z_{2k})
\Delta(z_1,\dots,z_{2k})^2}{\prod_{i=1}^{2k}\prod_{j=1}^{2k}
(z_i-\alpha_j)} \,dz_1\dots dz_{2k},
\end{eqnarray}
where one integrates about small circles enclosing the
$\alpha_j$'s, and where $\Xi$ is the set of $\binom{2k}{k}$
permutations $\sigma\in S_{2k}$ such that $\sigma(1)<\cdots <
\sigma(k)$ and $\sigma(k+1)<\cdots < \sigma(2k)$.
\end{lemma}

The above Lemma applies to the Unitary case, which has been
the subject of this section.  The next Lemma
is useful in the Symplectic and Orthogonal cases,
which will be addressed beginning in Section~\ref{sec:symplectic}.

\begin{lemma} \label{thm:concisesumsymplectic}

 Suppose $F$ is a symmetric function of $k$
variables, regular near $(0,\ldots,0)$, and $f(s)$ has a simple
pole of residue~$1$ at $s=0$ and is otherwise analytic in a
neighborhood of $s=0$, and let
\begin{equation}
K(a_1,\ldots,a_k)=F(a_1,\ldots,a_k) \prod_{1\leq i \leq j\leq k}
f(a_i+a_j)
\end{equation}
or
\begin{equation}
K(a_1,\ldots,a_k)=F(a_1,\ldots,a_k) \prod_{1\leq i < j\leq k}
f(a_i+a_j).
\end{equation}
If $\alpha_i+\alpha_j$ are contained in the
region of analyticity of $f(s)$ then
\begin{eqnarray}
\sum_{\epsilon_j =\pm 1 }&   &K(\epsilon_1 \alpha_1,\ldots,
\epsilon_k\alpha_k) =\nonumber \\&   &\ \ \ \ \frac{(-1)^{k(k-1)/2} } {(2\pi
i)^k} \frac{2^k}{k!} \oint \cdots \oint K(z_1,\ldots,z_k)
\frac{\Delta(z_1^2,\ldots,z_k^2)^2 \prod_{j=1}^k z_j }
{\prod_{i=1}^k\prod_{j=1}^k (z_i-\alpha_j)(z_i+\alpha_j)}
\,dz_1\cdots dz_k , \end{eqnarray} and
\begin{eqnarray}
\sum_{\epsilon_j=\pm 1  }&   &(\prod_{j=1}^k
\epsilon_j)K(\epsilon_1 \alpha_1,\ldots, \epsilon_k\alpha_k) = \nonumber \\&   &\ \ \ \ \frac{(-1)^{k(k-1)/2} } {(2\pi i)^k} \frac{2^k}{k!} \oint
\cdots \oint K(z_1,\ldots,z_k) \frac{\Delta(z_1^2,\ldots,z_k^2)^2
\prod_{j=1}^k \alpha_j } {\prod_{i=1}^k\prod_{j=1}^k
(z_i-\alpha_j)(z_i+\alpha_j)} \,dz_1\cdots dz_k, \end{eqnarray}
where the path of integration encloses the $\pm \alpha_j$'s.
\end{lemma}

The proofs of the lemmas come from the following

\begin{lemma}Suppose that $F(a;b)=F(a_1,\dots,a_m;b_1,\dots,b_{n})$
is symmetric  in the $a$ variables and in the  $b$ variables and is
regular near $(0,\dots,0)$. Suppose  $f(s)=\frac 1 s +c +\cdots$
and let
$$
G(a_1,\dots,a_m;b_1,\dots b_n)= F(a_1,\dots;\dots,b_n)
\prod_{i=1}^m\prod_{j=1}^n f(a_i-b_j).
$$
 Let $\Xi_{m,n}$ be as defined above.  Then 
\begin{eqnarray*}
 \sum_{\sigma \in \Xi_{m,n}}&&
G( \alpha_{\sigma(1)},\dots,  \alpha_{\sigma(m)};
    \alpha_{\sigma(m+1)}\dots  \alpha_{\sigma(m+n)})=\\
&&
 \frac{(-1)^{(m+n)}}{m!n!}\sum_{\sigma\in \pi_{m+n}}
\mbox{\rm Res}_{(z_1,\dots,z_{m+n})=
 (\alpha_{\sigma(1)},\dots,\alpha_{\sigma(m+n)})}
\frac{G(z_1,\dots,z_{m+n})
\Delta(z_1,\dots,z_{m+n})^2}{\prod_{i=1}^{m+n}\prod_{j=1}^{m+n}
(z_i-\alpha_j)}  .
\end{eqnarray*}
\end{lemma}

\begin{proof}
 It suffices to prove that
\begin{eqnarray*}
\mbox{Res}_{(z_1,\dots,z_{m+n})=
 (\alpha_{\sigma(1)},\dots,\alpha_{\sigma(m+n)})}
\frac{\Delta(z_1,\dots,z_{m+n})^2}{\prod_{i=1}^{m+n}\prod_{j=1}^{m+n}
(z_i-\alpha_j)}=(-1)^{m+n}
\end{eqnarray*}
since each such term will appear $m!n!$ times. Consider the case
where $\sigma$ is the identity permutation. Then the residue is
$$
\frac{\prod_{j<k}(\alpha_k-\alpha_j)^2}
{\prod_{j\ne k}(\alpha_j-\alpha_k)}=(-1)^{m+n};$$ the answer will
be the same for any permutation $\sigma$.

The residue above can be expressed as $(2\pi i)^{-m-n}$ times an
$m+n$ fold integral, each path of which encircles all of the poles
of the integrand; note that the value of such an integral may be
calculated by summing the residues and note that there is no
singularity when $z_j=z_k$ because of the factor $(z_k-z_j)^2$ in
the numerator.

\end{proof}

To obtain the form of Conjecture~\ref{thm:zeta2kconjecture} from the formulas at the end of
Sections~\ref{sec:Lmeanvalues} and~\ref{sec:thearithmeticfactor}, 
apply  Lemma~\ref{thm:concisesumunitary}
with $K(\alpha_1,\dots,\alpha_{2k})  =  W(\frac12;\alpha_1,\dots,\alpha_{2k})$. That is,
$f(z)=\zeta(1+z)$ and 
\begin{eqnarray}
F(\alpha_1,\dots,\alpha_{2k}) & = &
\left( \frac{Q^{\frac 2w}t}{2}
\right)^{\frac{w}{2} \sum_{j=1}^k \alpha _j-\alpha _{k+j}}
A_k(\tfrac12;\alpha_1,\dots,\alpha_{2k})\cr
 & = &
\exp\left(w\log\( \frac{Q^{\frac 2w}t}{2}\)
\cdot {\frac{1}{2} \sum_{j=1}^k \alpha _j-\alpha _{k+j}}
\right)
A_k(\tfrac12;\alpha_1,\dots,\alpha_{2k}).
\end{eqnarray}
We arrive at the general case of Conjecture~\ref{thm:zeta2kconjecture}:
\begin{conjecture} \label{thm:general2kconjecture}
 Suppose 
$\mathcal{L}(s)$ is a primitive $L$-function having
the properties listed at the beginning of Section~\ref{sec:Lmeanvalues},
and the mean value  $I_k(\mathcal{L},\alpha_1,\dots, \alpha_{2k},g)$
is given in~(\ref{eqn:generalIk}).  Then
\begin{equation}
I_k(\mathcal{L},\alpha_1,\dots, \alpha_{2k},g) =
\int_{-\infty}^\infty 
P_k\left(w \log \(\frac{Q^{\frac{2}{w}} t}{2}\) , \alpha\right) 
(1+O(t^{-\frac{1}{2}+\varepsilon})) g(t)\,dt ,
\end{equation}
where $P_k(x,\alpha)$ and $G(z_1,\dots,z_{2k})$ are as stated in 
Conjecture~\ref{thm:zeta2kconjecture}, except that
$A_k$ is the Euler product
\begin{equation}
A_k(z) =\prod_p \prod_{i=1}^k\prod_{j=1}^k
\left(1-\frac{1}{p^{1+z_i-z_{k+j}}}\right) \int_0^1 \prod_{j=1}^k
\mathcal{L}_p\left(\frac{e(\theta)}{p^{\frac12+z_j}}\right)
\overline{\mathcal{L}_p}\left(\frac{e(-\theta)}{p^{\frac12-z_{k+j}}}\right) .
\end{equation}
\end{conjecture}

Note that for the Riemann $\zeta$-function, $w=1$ and $Q=1/\sqrt{\pi}$
and $\mathcal{L}_p(x)=(1-x)^{-1}$,
so Conjecture~\ref{thm:zeta2kconjecture} is a special case of the above.
Also note that $w \log \(\frac{Q^{\frac{2}{w}} t}{2}\)$ is the mean
density of zeros of $\mathcal{L}(\frac12+it)$, or equivalently
the log conductor, as expected.

It remains to express the arithmetic factor $A_k$ in a more explicit
form, which we do in the next section.


\subsection{  Explicit versions of the arithmetic factor
}\label{sec:explicitarithmeticfactor}

The factor $A_k(s,\alpha)$ in the $2k$th moment of a primitive $L$-function
can be expressed in a simple form.

Recall, see Theorem~\ref{thm:Rkgeneral}, that $A_k$ is the Euler product
\begin{equation}
A_k(s; \alpha) =\prod_p 
B_p(s;\alpha_1,\dots,\alpha_{2k})
\prod_{i=1}^k\prod_{j=1}^k
\left(1-\frac{1}{p^{2s +\alpha_i-\alpha_{k+j}}}\right) ,
\end{equation}
where
\begin{equation} 
B_p(s;\alpha_1,\dots,\alpha_{2k}) =
\int_0^1 \prod_{j=1}^k
\mathcal{L}_p\left(\frac{e(\theta)}{p^{s+\alpha_j}}\right)
\overline{\mathcal{L}_p}\left(\frac{e(-\theta)}{p^{s-\alpha_{k+j}}}\right)
\,d\theta.
\end{equation}

\begin{lemma} \label{thm:Bplemma} If $\mathcal{L}_p(x) = (1-\gamma_p x)^{-1}$ with
$|\gamma_p|=1$ then
\begin{equation}
B_p(s;\alpha_1,\dots,\alpha_{2k}) = 
\prod_{i=1}^k\prod_{j=1}^k
\left(1-\frac{1}{p^{2s +\alpha_i-\alpha_{k+j}}}\right)^{-1}
\sum_{m=1}^k
    \prod_{i\neq m}
       \frac{\displaystyle \prod_{j=1}^k \(1-\frac{1}{p^{2s+\alpha_j-\alpha_{k+i}}}\)}
                         {1-p^{\alpha_{k+i}-\alpha_{k+m}}}.
\end{equation}
\end{lemma}

\begin{corollary} \label{thm:zeta2kAk}
If $\mathcal{L}_p(x) = (1-\gamma_p x)^{-1}$ with
$|\gamma_p|=0$ when $p|N$ and $|\gamma_p|=1$ otherwise, then 
\begin{equation}
A_k(s;\alpha_1,\dots,\alpha_{2k})=  
 \prod_{p\nmid N} 
\sum_{m=1}^k
    \prod_{i\neq m}
       \frac{\displaystyle \prod_{j=1}^k \(1-\frac{1}{p^{2s+\alpha_j-\alpha_{k+i}}}\)}
                         {1-p^{\alpha_{k+i}-\alpha_{k+m}}}
\times
\prod_{p|N} \prod_{i=1}^k\prod_{j=1}^k
\left(1-\frac{1}{p^{2s +\alpha_i-\alpha_{k+j}}}\right).
\end{equation}
In particular, if $\mathcal{L}(s)=L(s,\chi)$ with $\chi$ a 
Dirichlet character of conductor~$N$, where the Riemann $\zeta$-function is the
case $N=1$, then
\begin{equation}
A_1(s;\alpha_1, \alpha_2)=\prod_{p|N} 
\left(1-\frac{1}{p^{2s +\alpha_1-\alpha_{2}}}\right),
\end{equation}
\begin{equation}
A_2(s;\alpha_1, \alpha_2,\alpha_3, \alpha_4)=
 \zeta(4s+\alpha_1+\alpha_2-\alpha_3-\alpha_4)^{-1} 
  \prod_{p|N} \prod_{i=1}^2\prod_{j=1}^2
   \left(1-\frac{1}{p^{2s +\alpha_i-\alpha_{2+j}}}\right),
\end{equation}
and
\begin{eqnarray} & &A_3(s;\alpha_1,\dots, \alpha_6)  \cr
& & \phantom{XXXX}=
  \prod_{p\nmid N}\biggl( 1 -  \mbox{}  p^{-\sum_1^3{\alpha_i-\alpha_{3+i}}}
    \left(
        p^{\alpha_1}+ p^{\alpha_2}+ p^{\alpha_3}
    \right)
    \left(
        p^{-\alpha_4}+ p^{-\alpha_5}+ p^{-\alpha_6}
    \right)
    p^{-4s} \nonumber \\
&&\phantom{XXXXxxxxxxx}\mbox{} + p^{-\sum_1^3{\alpha_i-\alpha_{3+i}}}
    \Bigl(
        \left(
            p^{\alpha_1}+ p^{\alpha_2}+ p^{\alpha_3}
        \right)
        \left(
            p^{-\alpha_1}+ p^{-\alpha_2}+ p^{-\alpha_3}
        \right)
         \nonumber \\
&&\phantom{XXXXxxxxxxx}\qquad\qquad\qquad\qquad+\left(
            p^{\alpha_4}+ p^{\alpha_5}+ p^{\alpha_6}
        \right)
        \left(
            p^{-\alpha_4}+ p^{-\alpha_5}+ p^{-\alpha_6}
        \right)
        -2
    \Bigr)
    p^{-6s} \nonumber \\
&&\phantom{XXXXxxxxxxx}\mbox{} - p^{-\sum_1^3{\alpha_i-\alpha_{3+i}}}
    \left(
        p^{-\alpha_1}+ p^{-\alpha_2}+ p^{-\alpha_3}
    \right)
    \left(
        p^{\alpha_4}+ p^{\alpha_5}+ p^{\alpha_6}
    \right)
    p^{-8s} \nonumber \\
&&\phantom{XXXXxxxxxxx}\mbox{} + p^{-2\sum_1^3{\alpha_i-\alpha_{3+i}}} p^{-12s}\biggr)\cr
& & \phantom{XXXX}
\times \prod_{p|N} \prod_{i=1}^3\prod_{j=1}^3
\left(1-\frac{1}{p^{2s +\alpha_i-\alpha_{3+j}}}\right).
\end{eqnarray}
\end{corollary}
For $k\ge 3$ it is not possible to express $A_k$ as a finite product of
$\zeta$-functions.

\begin{proof}[Proof of Lemma~\ref{thm:Bplemma}]  Using
$\mathcal{L}_p(x) = (1-\gamma_p x)^{-1}$
and setting 
\begin{equation}
 q_j = \frac{\gamma_p}{p^{s+\alpha_j}} \hbox{\ \ \  and \ \ \ }
 q_{k+j}=\frac{\overline{\gamma_p}}{p^{s -\alpha_{k+j}}} \hbox{ \ \ \ for } j=1,\ldots,k,
\end{equation}
we have
\begin{eqnarray} \label{eqn:arithint} 
B_p(s;\alpha_1,\dots,\alpha_{2k}) &=& \int_0^1 \prod_{j=1}^k
\left(1-e(\theta) q_j\right)^{-1}\left(1-e(-\theta) q_{k+j}\right)^{-1}
\,d\theta \cr
&  = & \frac{(-1)^k}{\prod_{j=1}^k  q_j}
\int_0^1 e(k\theta) \prod_{j=1}^k
\left(e(\theta)-1/ q_j\right)^{-1}
\left(e(\theta)- q_{k+j}\right)^{-1}
\,d\theta \cr
&  = & \frac{(-1)^k}{\prod_{j=1}^k  q_j} \frac{1}{2\pi i}
\oint z^{k-1} \prod_{j=1}^k
\left(z-1/ q_j\right)^{-1}
\left(z- q_{k+j}\right)^{-1}
\,dz ,
\end{eqnarray}
where the path of integration is around the unit circle.
Since $|q_j|<1$, by the residue theorem we have a contribution
from the poles at $q_{k+1},\dots,q_{2k}$, giving
\begin{eqnarray} 
B_p(s;\alpha_1,\dots,\alpha_{2k})&  = & \frac{(-1)^k}{\prod_{j=1}^k  q_j}
\sum_{m=1}^k q_{k+m}^{k-1}
 \prod_{i=1}^k
\left(q_{k+m}- q_i^{-1}\right)^{-1}
\prod_{i\not=m}\left(q_{k+m}- q_{k+i}\right)^{-1}
				 \cr
&  = & \sum_{m=1}^k 
 \prod_{i=1}^k
\left(1-q_i q_{k+m}\right)^{-1}
\prod_{i\not=m}\left(1- q_{k+i} q_{k+m}^{-1}\right)^{-1}
\end{eqnarray}

Since
\begin{equation}
\prod_{i,j=1}^k (1- q_i q_{k+j})
 \prod_{i=1}^k
\left(1-q_i q_{k+m}\right)^{-1}
=
\prod_{j\not=m}  \prod_{i=1}^k (1- q_i q_{k+j}),
\end{equation}
factoring out
\begin{equation}
\prod_{i,j=1}^k (1- q_i q_{k+j})^{-1},
\end{equation}
we have
\begin{equation}
B_p(s;\alpha)= \left(
    \prod_{i,j=1}^k (1- q_i q_{k+j})^{-1}
\right)
\sum_{m=1}^k
    \prod_{i\neq m} \frac{\prod_{j=1}^k ( 1- q_j q_{k+i})}
                         {1- q_{k+i} q_{k+m}^{-1}}.
\end{equation}
Since
\begin{equation}
q_j q_{k+i} = p^{-2s-\alpha_j+\alpha_{k+i}}
\hbox{\ \ \ and\ \ \ } 
q_{k+i} q_{k+m}^{-1} = p^{\alpha_{k+i}-\alpha_{k+m}}
,
\end{equation}
we obtain the formula in the Lemma.
\end{proof}

Notice that the special case $N=1$, i.e.~the Riemann $\zeta$ function,
reads in Corollary \ref{thm:zeta2kAk}, 
\begin{equation}
A_k(s;\alpha_1,\dots,\alpha_{2k})=  
 \prod_{p} 
\sum_{m=1}^k
    \prod_{i\neq m}
       \frac{\displaystyle \prod_{j=1}^k \(1-\frac{1}{p^{2s+\alpha_j-\alpha_{k+i}}}\)}
                         {1-p^{\alpha_{k+i}-\alpha_{k+m}}}.
\end{equation}
Each local factor
\begin{equation}
\label{eqn:local A_p,k}
A_{p,k}(s;\alpha)=  
\sum_{m=1}^k
    \prod_{i\neq m}
       \frac{\displaystyle \prod_{j=1}^k \(1-\frac{1}{p^{2s+\alpha_j-\alpha_{k+i}}}\)}
                         {1-p^{\alpha_{k+i}-\alpha_{k+m}}}.
\end{equation}
is actually a polynomial in $p^{-2s}$, $p^{-\alpha_j}$ and $p^{\alpha_{k+j}}$, for $j=1,\dots,k$.
That this is so in $p^{-2s}$ and $p^{-\alpha_j}$ is readily apparent from~(\ref{eqn:local A_p,k}).
The fact that it is also a polynomial in $p^{\alpha_{k+j}}$ follows from~(\ref{eqn:bigAk}) 
and~(\ref{eqn:bigBk}), from which
\begin{equation}
    A_{p,k}(s;\alpha_1,\dots,\alpha_k,\alpha_{k+1},\dots,\alpha_{2k})= 
    A_{p,k}(s;-\alpha_{k+1},\dots,-\alpha_{2k},-\alpha_{1},\dots,-\alpha_{k}) .
\end{equation}
Setting $\beta_1=-\alpha_{k+1},\dots,\beta_k=-\alpha_{2k}$, one has,
from the above discussion, that $A_{p,k}$ is a polynomial in
$p^{-\beta_j}$, i.e. in $p^{\alpha_{k+j}}$, for $j=1,\dots,k$.
Finally, use the fact that if an analytic function of several variables is
of polynomial growth in each variable separately, then it must
be a polynomial.


\subsection{  Recovering the leading order for moments of $\zeta$
}\label{sec:leadingorder}
Conjecture~\ref{thm:zeta2kconjecture} contains, as a special case,
 a conjecture for
the leading order term for the moments of the Riemann zeta function.
In this section we show that the leading order terms derived from
Conjecture~\ref{thm:zeta2kconjecture} agree with the 
leading order terms which have previously been conjectured by other methods.

As described in Section~\ref{sec:examplefamilies}, it is
conjectured that the mean values of the Riemann zeta function take
the form
\begin{equation}  \int_0^T |\zeta(\tfrac12 +it)|^{2k}\,dt  = T \, {\mathcal P}_{k}(\log T) +O(T^{\frac12 +\varepsilon}),
\end{equation}
where ${\mathcal P}_k(\log T)$ is a polynomial in $\log T$ of
degree $k^2$. Conrey and Ghosh conjectured that the
coefficient of the $\log ^{k^2} T$ term is of the form $g_ka_k/k^2!$,
where $a_k$ is given by~(\ref{eqn:zeta2kak}).  Keating and Snaith
used random matrix theory to conjecture that
$g_k$ is given by~(\ref{eqn:zeta2kgk}).
This leading order term $g_ka_k/k^2!$ will be re-derived
here, starting with Conjecture~\ref{thm:zeta2kconjecture}.

Conjecture~\ref{thm:zeta2kconjecture} implies
\begin{equation}
\int_{0}^T |\zeta(\tfrac12 +it)|^{2k}\,dt= \int_0^T
P_k\left(\log \tfrac{t}{2 \pi}\right) \,dt
+
O(T^{\frac{1}{2}+\varepsilon}) ,
\end{equation}
where $P_k$ is the polynomial of degree $k^2$ given by 
\begin{eqnarray}
P_k(x)&  = &\frac{(-1)^k}{k!^2(2\pi
i)^{2k}}\oint \cdots \oint A_k(z_1,\ldots,z_{2k}) \prod_{i=1}^k
\prod_{j=1}^k \zeta(1+z_i-z_{k+j})\nonumber \\&   &\ \ \ \ \times
\frac{ \Delta^2(z_1,\ldots,z_{2k})} {\prod_{j=1}^{2k} z_j^{2k}}
e^{\tfrac{x}{2} \sum_{j=1}^kz_j-z_{k+j}}\, dz_1\cdots dz_{2k} .
\end{eqnarray}
Our goal is to show that the leading order term of $P_k(x)$
is $(g_ka_k/k^2!) x^{k^2}$.

Using the fact that $A_k$ is analytic in a neighborhood of
$(0,\dots,0)$ and the $\zeta$-function has a simple pole at~$1$ with
residue~$1$, after a change of variables we have
\begin{eqnarray}
P_k(x)
&  = &\frac{(-1)^k}{k!^2(2\pi i)^{2k}} \oint \cdots
\oint A_k\(\frac{z_1}{x/2},\cdots,\frac{z_{2k}}{x/2}\)
\prod_{i=1}^k \prod_{j=1}^k
\zeta\(1+\frac{z_i- z_{k+j}}{x/2}\)\nonumber \\
&   &\ \
\ \ \times \frac{ \Delta^2(z_1,\ldots,z_{2k})} {\prod_{j=1}^{2k}
z_j^{2k}} e^{\sum_{j=1}^kz_j-z_{k+j}} \, dz_1\cdots dz_{2k} \cr
&  = & \frac{(-1)^k}{k!^2}\frac{A_k(0,\dots,0)}{(2\pi i)^{2k}} 
\(\frac{x}{2}\)^{k^2} 
(1+O(x^{-1}))
\oint \cdots \oint 
\frac{\Delta^2(z_1,\ldots,z_{2k}) } {\left(
\prod_{i=1}^k\prod_{j=1}^k (z_{i}-z_{k+j})\right) \prod_{j=1}^{2k}
z_j^{2k}} \nonumber \\&   &\ \ \ \ \ \ \ \times e^{\sum_{j=1}^k
z_j-z_{k+j}} \,dz_1\cdots dz_{2k} \nonumber \\&  =
&\frac{A_k(0,\dots,0)}{k!^2 2^{k^2} (2\pi i)^{2k} } x^{k^2} 
(1+O(x^{-1}))
\oint \cdots \oint 
\frac{\Delta(z_1,\ldots,z_{2k}) \Delta(z_1,\ldots,z_k)
\Delta(z_{k+1},\ldots,z_{2k}) } {\prod_{j=1}^{2k}
z_j^{2k}}\nonumber \\&   &\ \ \ \ \ \ \ \times
e^{\sum_{j=1}^kz_j-z_{k+j}} \, dz_1\ldots,dz_{2k}. 
\end{eqnarray}
Now we need only show that $A_k(0,\ldots0)=a_k$ and the remaining
factors give $g_k/k^2!$.

From Conjecture \ref{thm:zeta2kconjecture},
\begin{equation}
A_k(0,\ldots0)=\prod_p \left(1-\frac{1}{p}\right)^{k^2}\int_0^1
\left(1-\frac{e(\theta)}{p^{\frac12 }}\right)^{-k}
\left(1-\frac{e(-\theta)}{p^{\frac12 }}\right)^{-k}\,d\theta.
\end{equation}
For a given $p$, we concentrate on the integral in the above
expression, writing it as a contour integral around the unit
circle:
\begin{equation}
(-p^{\frac12 })^k \frac{1}{2\pi i} \oint \frac{z^{k-1}(z-p^{\frac12 })^{-k}}
{(z-p^{-\frac12 })^k} dz.
\end{equation}
  After expanding the two factors in the numerator around
$z=p^{-\frac12 }$ and calculating the residue we are left with the sum
\begin{equation}
(1-\tfrac{1}{p})^{-2k+1}\sum_{\ell=0}^{k-1}\binom{k-1}{\ell}
\binom{2k-\ell-2}{k-1} p^{-k+\ell+1}(1-\tfrac{1}{p})^{\ell}.
\end{equation}
Next one can perform a binomial expansion of
$(1-\frac{1}{p})^{\ell}$ and gather like powers of $\frac{1}{p}$
to obtain
\begin{equation}
(1-\tfrac{1}{p})^{-2k+1}\sum_{m=0}^{k-1}\left(\sum_{q=0}^m (-1)^q
\binom{k+q-m-1}{q}\binom{k-1}{k+q-m-1}\binom{k+m-q-1}{k-1} \right)
p^{-m}.
\end{equation}
A simple manipulation of the binomial coefficients and replacing
$q$ by $m-q$ gives
\begin{eqnarray}
&&\sum_{q=0}^m(-1)^q\binom{k+q-m-1}{q}\binom{k-1}{k+q-m-1}\binom{k+m-q-1}{k-1}\nonumber
\\
&&\qquad\qquad\qquad =\binom{k-1}{m}\sum_{q=0}^m
(-1)^{m-q}\binom{m}{q}\binom{k+q-1}{q},
\end{eqnarray}
and this final sum over $q$ is in fact just $\binom{k-1}{m}$ (see,
for example, \cite{Rio}).  Thus,
\begin{equation}
A_k(0,\ldots,0)=\prod_p \left(1-\frac{1}{p}\right)^{(k-1)^2}
\sum_{m=0}^{k-1}\binom{k-1}{m}^2 p^{-m} ,
\end{equation}
and this is indeed equal to $a_k$ defined in~(\ref{eqn:zeta2kak}).

Now we must identify the remaining terms as
$g_k/k^2 !$ defined in~(\ref{eqn:zeta2kgk}), as~$x\to\infty$.
The method applied below as far as (\ref{eqn:firstdeterminant})
follows closely that used for a similar purpose in \cite{BH}.
Expanding the determinants
$\Delta(z_1,\ldots,z_k)=\det[z_j^{m-1}]_{j,m=1}^k$, we obtain
\begin{eqnarray}
\lim_{x\rightarrow \infty}\frac{P_k(x)}{a_k{x}^{k^2}}&  =
&\frac{1} {(k!)^2 2^{k^2} (2\pi i)^{2k}} \oint \cdots  \oint
e^{\sum_{j=1}^k z_j-z_{k+j}} \nonumber \\&   &\ \ \ \ \times\left(
\sum_S \sgn(S) z_1^{S_0}z_2^{S_1} \cdots z_k^{S_{k-1}}
z_{k+1}^{S_k} \cdots z_{2k}^{S_{2k-1}} \right) \left( \sum_Q
\sgn(Q) z_1^{Q_0}\cdots z_k^{Q_{k-1}} \right)\nonumber \\&   &\ \
\ \ \times\left( \sum_R \sgn(R) z_{k+1}^{R_0}\cdots
z_{2k}^{R_{k-1}} \right) z_1^{-2k}\cdots z_{2k}^{-2k} dz_1\cdots
dz_{2k}.
\end{eqnarray}
Here $Q$ and $R$ are permutations of $\{0,1,\ldots,k-1\}$ and
$S$ is a permutation of $\{0,1,\ldots,2k-1\}$.

Since the integrand is symmetric amongst $z_1,\ldots,z_k$ and also
amongst $z_{k+1},\ldots,z_{2k}$, in each term of the sum over $Q$
we permute the variables $z_1,\ldots,z_k$ so that $z_j$ appears
with the exponent $j-1$, for $j=1,\ldots, k$.  In the sum over $S$
the effect is to redefine the permutations, and the additional
sign involved with this exactly cancels $\sgn(Q)$.  We do the same
with the sum over $R$, and as a result we are left with $k!^2$
copies of the sum over the permutation $S$:
\begin{eqnarray}
\label{eqn:lastintegral}
 \lim_{x\rightarrow
\infty}\frac{P_k(x)}{a_k {x}^{k^2}}& = &\frac{1} {2^{k^2} (2\pi
i)^{2k}} \oint \cdots  \oint e^{\sum_{j=1}^k z_j-z_{k+j}}\nonumber
\\&   &\ \ \ \  \times\sum_S \sgn(S) z_1^{-(2k-S_0)}
z_2^{-(2k-S_1-1)} \cdots z_k^{-(2k-S_{k-1}-(k-1))}\nonumber \\&
&\ \ \ \ \times z_{k+1}^{-(2k-S_k)} z_{k+2}^{-(2k-S_{k+1}-1)}
\cdots z_{2k}^{-(2k-S_{2k-1}-(k-1))} dz_1\cdots dz_{2k}.
\end{eqnarray}

Since
\begin{equation}
 \frac{1}{\Gamma(z)}=\frac{1}{2\pi i} \int_C
(-t)^{-z} e^{-t} (-dt),
\end{equation}
where the path of integration $C$ starts at $+\infty$ on the real
axis, circles the origin in the counterclockwise direction and
returns to the starting point, we can rewrite
(\ref{eqn:lastintegral}) as
\begin{eqnarray}
\lim_{x\rightarrow \infty}\frac{P_k(x)}{a_k x^{k^2}}
&=&\frac{(-1)^{k}}{2^{k^2}} \sum_S \sgn(S) \Big(\Gamma(2k-S_0) \Gamma(2k-S_1-1)
\cdots \Gamma(2k-S_{k-1}-(k-1)) \nonumber \\
& &\phantom{xxxxxxxxxxxxx}  \times(-1)^{S_k}
\Gamma(2k-S_k) (-1)^{S_{k+1}+1} \Gamma(2k-S_{k+1}-1) \cdots\nonumber \\
& &\phantom{xxxxxxxxxxxxx}  \times (-1)^{S_{2k-1}+k-1}
\Gamma(2k-S_{2k-1}-(k-1))\Big)^{-1} \nonumber \\
&=&\frac{(-1)^{k}}{2^{k^2}}\vmatrix \frac{1}{\Gamma(2k)} &
\frac{1}{\Gamma(2k-1)} & \cdots & \frac{1}{\Gamma(k+1)} &
\frac{1}{\Gamma(2k)}& \frac{-1}{\Gamma(2k-1)}& \cdots &
\frac{(-1)^{k-1}}{\Gamma(k+1)}
\\ \frac{1}{\Gamma(2k-1)} & \frac{1}{\Gamma(2k-2)} &\cdots &
\frac{1}{\Gamma(k)} & \frac{-1}{\Gamma(2k-1)} &
\frac{1}{\Gamma(2k-2)} & \cdots & \frac{(-1)^{k}}{\Gamma(k)} \\
\vdots& \vdots&\ddots&\vdots&\vdots &\vdots&\ddots&\vdots\\
\frac{1}{\Gamma(1)} & \frac{1}{\Gamma(0)} & \cdots &
\frac{1}{\Gamma(2-k)} & \frac{-1}{\Gamma(1)} & \frac{1}{\Gamma(0)}
& \cdots & \frac{(-1)^{3k-2}}{\Gamma(2-k)}
\endvmatrix\nonumber \\
&=&\frac{(-1)^{k}}{2^{k^2}}\left(\prod_{\ell=0}^{k-1} \frac{\ell!}{(k+\ell)!}
\right)\nonumber \\
 &&\ \ \ \ \times \vmatrix\binom{0}{0} & \binom{0}{1} &
\cdots & \binom{0}{k-1} & \binom{0}{0} & -\binom{0}{1} & \cdots &
(-1)^{k-1}\binom {0}{k-1}
\\ \vdots & \vdots & \ddots & \vdots & \vdots & \vdots & \ddots &
\vdots \\ \binom{2k-1}{0} & \binom{2k-1}{1} & \cdots &
\binom{2k-1}{k-1} & -\binom{2k-1}{0} & \binom{2k-1}{1} & \cdots &
(-1)^k \binom{2k-1}{k-1}
\endvmatrix.
 \label{eqn:firstdeterminant} \end{eqnarray}
The above is a $2k\times 2k$ determinant, the first $k$ columns of
which are identical to the first $k$ columns of the matrix
\begin{equation}
\label{eqn:magicmatrix1} \left(\begin{array}{cccccccccccc}
\binom{0}{0} & \binom{0}{1} & \cdots & \binom{0}{2k-1} \\ \vdots &
\vdots & \ddots & \vdots
\\ \binom{2k-1}{0} & \binom{2k-1}{1} & \cdots &
\binom{2k-1}{2k-1} \end{array}\right).
\end{equation}
The matrix (\ref{eqn:magicmatrix1}) is lower triangular and so can
easily be seen to have determinant equal to one.  It is also the
inverse of
 \begin{equation}
\left(\begin{array}{cccccccccccc} \binom{0}{0} & -\binom{0}{1} &
\cdots & -\binom{0}{2k-1}
\\-\binom{1}{0} & \binom{1}{1} & \cdots & \binom{1}{2k-1} \\
\vdots & \vdots & \ddots & \vdots
\\ -\binom{2k-1}{0} & \binom{2k-1}{1} & \cdots &
\binom{2k-1}{2k-1} \end{array}\right).\label{eqn:magicmatrix2}
\end{equation}
It so happens that matrix (\ref{eqn:magicmatrix2}) has its $k$
first columns identical to columns $k+1$ through $2k$ of the
matrix in expression (\ref{eqn:firstdeterminant}).  Therefore we
can multiply expression (\ref{eqn:firstdeterminant}) by the
determinant of (\ref{eqn:magicmatrix1}) (which is equal to one)
and this simplifies the final $k$ columns of the resulting
determinant significantly:
\begin{eqnarray*}
&\displaystyle{\lim_{x\rightarrow \infty}
\frac{P_k(x)}{a_k x^{k^2}} \hskip -0.5in }& \nonumber \\
&=&\frac{(-1)^{k}}{2^{k^2}} \left( \prod_{\ell=0}^{k-1} \frac{\ell!}{(k+\ell)!}
\right) \nonumber \\
& &\times
     \vmatrix \binom{0}{0} & \binom{0}{1} & \cdots & \binom{0}{2k-1} \nonumber  \\
\vdots & \vdots & \ddots & \vdots
\\ \binom{2k-1}{0} & \binom{2k-1}{1} & \cdots &
\binom{2k-1}{2k-1} \endvmatrix \cdot \vmatrix \binom{0}{0} &
\cdots & \binom{0}{k-1} & \binom{0}{0} & \cdots &
(-1)^{k-1}\binom{0}{k-1}
\\ \vdots & \ddots & \vdots & \vdots & \ddots & \vdots \\
\binom{2k-1}{0} & \cdots & \binom{2k-1}{k-1} & -\binom{2k-1}{0} &
\cdots & (-1)^k\binom{2k-1}{k-1}
\endvmatrix \nonumber  \\
&=&
\frac{(-1)^{k}}{2^{k^2}} \left( \prod_{\ell=0}^{k-1} \frac{\ell!}{(k+\ell)!}
\right)\vmatrix \binom{0}{0} & 0 & \cdots & 0 & 1& 0& \cdots & 0
\\ 2\binom{1}{0} & \binom{1}{1} & \cdots & 0 & 0 & 1 & \cdots &
0\\ \vdots & \vdots & \ddots & \vdots & \vdots & \vdots & \ddots &
\vdots \\ 2^{k-1} \binom{k-1}{0} & 2^{k-2}\binom{k-1}{1} & \cdots
& \binom{k-1}{k-1} & 0& 0& \cdots &1 \\ 2^{k}\binom{k}{0} &
2^{k-1}\binom{k}{1} & \cdots & 2\binom{k}{k-1} & 0&0&\cdots &0
\\
\vdots & \vdots & \ddots & \vdots & \vdots & \vdots & \ddots &
\vdots \\ 2^{2k-1}\binom{2k-1}{0} & 2^{2k-2}\binom{2k-1}{1} &
\cdots & 2^{k}\binom{2k-1}{k-1} & 0& 0& \cdots &0
\endvmatrix \nonumber \\
 &=&
\frac{(-1)^{k(k-1)/2}}{2^{k^2}} \left( \prod_{\ell=1}^{k-1}
\frac{\ell!}{(k+\ell)!}\right) \vmatrix 2^{2k-1}\binom{2k-1}{0} &
2^{2k-2}\binom{2k-1}{1} &\cdots & 2^{k}\binom{2k-1}{k-1} \\
\vdots & \vdots & \ddots & \vdots \\ 2^k\binom{k}{0} &
2^{k-1}\binom{k}{1} &\cdots & 2\binom{k}{k-1} \endvmatrix \nonumber \\
&=&
(-1)^{k(k-1)/2} \left( \prod_{\ell=1}^{k-1}
\frac{\ell!}{(k+\ell)!}\right)\vmatrix \binom{2k-1}{0} &
\binom{2k-1}{1} &\cdots & \binom{2k-1}{k-1} \\
\vdots & \vdots & \ddots & \vdots \\ \binom{k}{0} & \binom{k}{1}
&\cdots & \binom{k}{k-1} \endvmatrix.
 \end{eqnarray*}

The matrix above can be decomposed as
\begin{eqnarray}
&   &\left(\begin{array}{cccccccccccc} \binom{2k-1}{0}&
&\binom{2k-1}{1}&   &\cdots&   &\binom{2k-1}{k-1}
\nonumber \\ \binom{2k-2}{0}&   &\binom{2k-2}{1}&   &\cdots&   &\binom{2k-2}{k-1} \nonumber \\
\vdots&   &\vdots&   &\ddots&   &\vdots \nonumber \\ \binom{k}{0}&
&\binom{k}{1}&   &\cdots&   &\binom{k}{k-1}
\end{array}\right)\nonumber \\
&   &\ \ \ \ \ \  =\left(\begin{array}{cccccccccccc}
\binom{k-1}{0}&   &\binom{k-1}{1}&   &\cdots&   &\binom{k-1}{k-1}
\nonumber \\ \binom{k-2}{0}&   &\binom{k-2}{1}&   &\cdots&   &\binom{k-2}{k-1} \nonumber \\ \vdots&   &\vdots&   &\ddots&   &\vdots \nonumber \\
\binom{0}{0}&   &\binom{0}{1}&   &\cdots&   &\binom{0}{k-1}
\end{array}\right) \times \left(\begin{array}{cccccccccccc}
\binom{k}{0}&   &\binom{k}{1}&   &\cdots&   &\binom{k}{k-1} \nonumber \\
\binom{k}{-1}&   &\binom{k}{0}&   &\cdots&   &\binom{k}{k-2}
\nonumber \\\vdots&   &\vdots&   &\ddots&   &\vdots \nonumber \\
\binom{k}{-k+1}&   &\binom{k}{-k+2}&   &\cdots&   &\binom{k}{0}
\end{array}\right) .
\end{eqnarray}
The first matrix on the right side is zero in the lower right
triangle, and the second matrix on the right side is upper
triangular. Thus we read that the determinant of the matrix on the
left hand side is $(-1)^{k(k-1)/2}$. Therefore,
\begin{equation}
\lim_{x\rightarrow \infty}\frac{P_k(x)}{a_k x^{k^2}}=
\prod_{\ell=1}^{k-1} \frac{\ell!}{(k+\ell)!},
\end{equation}
and this is $g_k/k^2!$ from~(\ref{eqn:zeta2kgk}), as required. 

A similar method applies to the orthogonal and symplectic cases.



\section{  Families of characters and families of $L$-functions
}\label{sec:families}
We will describe a particular kind of ``family'' of primitive  $L$-functions
based on the idea of twisting a single $L$-function by a family of
``characters.''  
In the next section we provide a general recipe for conjecturing the 
critical mean value of products of $L$-functions averaged
over a family and we demonstrate the recipe in  
several examples.

Note that we use ``character'' somewhat more generally
than is usually covered by that term.


\subsection{  Families of primitive characters
}\label{sec:familiesofcharacters} We describe sets of arithmetic functions
which we call ``families of characters.''

Let $\mathcal{F} = \{f\}$ be a collection of arithmetic
functions $f(n)$,  and
assume that for each $f\in \mathcal{F}$ the associated $L$-function
$L_f(s)=\sum f(n)\;n^{-s}$ is a primitive $L$-function 
with  functional equation $L_f(s)=\varepsilon_f X_f(s)\overline{L}_f(1-s)$ and
an Euler product of the form
\begin{equation}\label{eqn:genep}
L_f(s)=\sum_{n=1}^\infty
\frac{f(n)}{n^s}=\prod_p\prod_{j=1}^v(1-\beta_{p,j}/p^s)^{-1} .
\end{equation}
The quantity 
\begin{equation}
c(f)=|(\varepsilon_f X_f)^\prime(\tfrac12)|
\end{equation}
is called the {\it log conductor} of $f$.

Note that if $f=\chi$, a primitive Dirichlet character of conductor~$q$,
then the log conductor is
\begin{equation}
c(\chi)=
\begin{cases}
\log q  -\log \pi +\frac{\Gamma'}{\Gamma}(\tfrac 14) & \chi \text{\ \ even} \\
\log q  -\log \pi +\frac{\Gamma'}{\Gamma}(\tfrac 34) & \chi \text{\ \ odd} .
\end{cases}
\end{equation}
If $f(n)=n^{-it}$ then the log conductor is
$c(n^{-it})=\log\frac{t}{2\pi}+O(t^{-1})$.
Generally the log conductor $c(f)$ 
scales as
the $\log$ of the ``usual'' conductor of~$f$.

In the case that  $\mathcal{F}$ is finite, we
require that the data $Q$, $w_j$, $\mu_j$ in the functional equation
(\ref{eqn:gammafactors}) is \emph{the same}  for all~$f\in \mathcal{F}$.
In particular,
the conductor
$c(f)$ is the same for all~$f\in \mathcal{F}$.

In the case that  $\mathcal{F}$ is infinite, we require that
the data $Q$, $w_j$, $\mu_j$ in the functional equation
(\ref{eqn:gammafactors}) are \emph{monotonic functions of the 
conductor} $c(f)$.  Furthermore, we
define the counting function 
$ M(T)=\#\{f\ :\ c(f)\le T\}$ and
require that  $M(\log(T))= F(T^A, \log T) + O(T^{\frac{A}{2}+\epsilon})$ 
for all $\epsilon>0$,
where $A>0$ and $F(\cdot,\cdot)$ is a polynomial.

If $G$ is a function on $\mathcal{F}$, then we define the
{\it expected value} of $G$ by
\begin{equation}
\langle G(f)\rangle =
\lim_{T\to \infty}
M(T)^{-1}\sum_{\ontop{f\in \mathcal{F}}{c(f)<T}}
G(f) ,
\end{equation}
assuming the limit exists.  In the case of a continuous family, the
sum is an integral.  

We require that if $m_1,\dots,m_k$ are integers then
the expected value
\begin{equation}
\delta_{\ell}(m_1,\dots,m_k)
=
\langle f(m_1)\dots
f(m_\ell)\overline{f(m_{\ell+1})\dots f(m_k)} \rangle 
\end{equation}
exists and is multiplicative.  That is,
if $(m_1 m_2 \dots m_k, n_1 n_2 \dots n_k)=1$, then
\begin{equation}\label{eqn:deltamultiplicative}
\delta_\ell(m_1n_1,m_2n_2,\dots m_kn_k)=
\delta_\ell(m_1,\dots,m_k)\delta_\ell(n_1,\dots n_k).
\end{equation}
We sometimes refer to $\delta$ as the ``orthogonality relation''
of the family.

The practical use of being multiplicative is that a multiple Dirichlet
series with $\delta_\ell$ coefficients factors has an Euler product:
\begin{eqnarray}\label{eqn:multipleEP}
\sum_{m_1,\ldots,m_\ell} 
\frac{\delta_\ell(m_1,\ldots,m_\ell)}{m_1^{s_1}\cdots m_\ell^{s_\ell}}
=
\prod_p \sum_{e_1,\ldots,e_\ell} 
\frac{\delta_\ell(p^{e_1},\ldots,p^{e_\ell})}{p^{e_1s_1+\cdots+e_\ell s_\ell}} .
\end{eqnarray}
We will use the above relation in our calculations.

To summarize, a family of characters $ \mathcal{F}=\{f\}$  is a 
collection of arithmetic
functions, each of which are the coefficients of a particular
kind of $L$-function.  The characters are partially ordered 
by conductor~$c(f)$, and the expected values
$\delta_\ell(m_1,\dots,m_k)$ are multiplicative functions.

The following are examples of families of characters, two of which are
finite and two are infinite.  The term ``finite family'' is
somewhat misleading, because those families depend on a parameter,
and the size of the family grows with the parameter.    

\noindent {\bf 1.  The family of $t$-twists.}
\begin{equation}
{\mathcal F}_{\text{t}}=\{f_t(n)=n^{-it},0<t<T\}.
\end{equation}
We have
\begin{equation}
\frac{1}{T}\int_0^T (m/n)^{it} dt =
\begin{cases}
1 & n=m \\
\frac{(m/n)^{i T} -1}{T\log(m/n)} & otherwise ,
\end{cases}
\end{equation}
leading to the
expected values
\begin{equation}
\langle f_t(n)\overline{f_t(m)} \rangle=\langle  n^{-it} m^{it} \rangle
= \langle  (m/n)^{it} \rangle 
=
\begin{cases}
1& n=m\\
0& \text{otherwise}.
\end{cases}
\end{equation}
Therefore the orthogonality relation is
\begin{equation}
\delta_\ell(n_1,\dots,n_{k})=\delta(n_1\cdots n_\ell = n_{\ell+1}\cdots n_{k}).
\end{equation}


\noindent {\bf 2. The family of primitive Dirichlet characters.}
For each positive integer $q$ we set
\begin{equation}{\mathcal F}_{\text{ch}}(q)=\{f_\chi(n)=\chi(n): \chi \text{ is a primitive
Dirichlet character mod $q$}\}.
\end{equation}
We have
\begin{equation}
\frac{1}{q^*}\sumstar_{\chi \mod q} \chi(n)\overline{\chi(m)} =
\begin{cases}
1 & n\equiv m \mod q \ \text{ and } \ (mn,q)=1 \\
0 & otherwise ,
\end{cases}
\end{equation}
where
the sum is over the primitive characters mod~$q$, and $q^*$ is the number 
terms in the sum.
This leads to the
expected values
\begin{equation}
\langle f_\chi(n)\overline{f_\chi(m)}  \rangle 
=
\langle \chi(n)\overline{\chi(m)}  \rangle 
=
\begin{cases}
1& n= m \ \text{ and } \ (mn,q)=1\\
0& otherwise.
\end{cases}
\end{equation}
Since  $\chi(m_1)\chi(m_2)=\chi(m_1m_2)$   we obtain 
\begin{equation}
\delta_\ell(m_1,\dots,m_k)=\delta\Bigl(m_1\cdots m_\ell=m_{\ell+1}\cdots m_k
\ \ \mbox{and}\ \ \ (m_1\cdots m_k, q)=1\Bigr).
\end{equation}
Note that the condition in the definition of $\delta_\ell$ is \emph{not}
$m_1\cdots m_\ell\equiv m_{\ell+1}\cdots m_k\mod q$.  We are computing the
expected value as a function of $q$, so one should think of the $m_j$ as
fixed and $q\to\infty$.  The only way to have
$m_1\cdots m_\ell\equiv m_{\ell+1}\cdots m_k\mod q$ for sufficiently
large $q$ is to have actual equality.

Note that by our definition, ${\mathcal F}_{\text{ch}}(q)$ is not a family,
but it is the union of two families consisting of the even characters
and the odd characters separately.

\noindent {\bf 3.  The family of real primitive Dirichlet characters.}
\begin{equation}
{\mathcal F}_{\text{d}} =\{f_d(n)= \chi_d(n): \chi_d \text { is a primitive
real character $\mod d$, $|d|<X$  }\},
\end{equation}
where $d$ runs over fundamental discriminants.
 We have 
expected values
\begin{equation}
\langle f_d(n)f_d(m)\rangle =\langle  \chi_d(n)\chi_d(m) \rangle
=
\langle  \chi_d(nm) \rangle
=
\begin{cases}
\prod_{p|nm} (1+\frac1p)^{-1} & nm =\square \\
0& \text{otherwise}. 
\end{cases} .
\end{equation}
The calculation in the case $nm =\square$ is nontrivial and was first done
by Jutila\cite{J}.  (If one were summing over all $d$ then the
expected value when $nm =\square$ would be $\varphi(nm)/nm$).
 
In practice one encounters more restricted families, so 
we let ${\mathcal F}_{\text{d}}(+)=\{f_d:d>0\}$ and 
${\mathcal F}_{\text{d}}(-)=\{f_d:d<0\}$,
and also
\begin{equation}{\mathcal F}_{\text{d}}(a,N,\pm)=
\{f_d\in {\mathcal F}_{\text{d}}(\pm):d \equiv a \bmod N\}.
\end{equation}
For the family ${\mathcal F}_{\text{d}}(a,N,\pm)$, 
evaluating the expected value of $\chi_d(n)$ can be tricky, so
we provide some useful asymptotics.

Below we restrict ourselves to $0 < d < X$, but the same
asymptotics hold if one restricts to $-X < d < 0$. 

\begin{theorem}\label{thm:funddisc}
Let $Q=\gcd(a,N)$ not be divisible by the square of an odd prime. Then
\begin{equation}
\label{eqn:fund discr progressions}
    \sumstar_{\ontop{0 < d < X}{d \equiv a \bmod N}} 1
    \sim \frac{1}{\phi(4N/Q)} \frac{X}{Q} \frac{6}{\pi^2}
    h_2(a,N)
    \prod_{p|2N} \frac{p}{p+1}.
\end{equation}
Next, assume further that either $N$ is odd or divisible by at least 8
(this condition is related to the fact that $\chi_d(2)$ is
periodic mod 8), and say
$n=g\square$, with  $( \square,N)=1$, and  with all prime factors of 
$g$ being prime factors of $N$. Then 
\begin{equation}
    \label{eqn:chi over fund discr}
    \sumstar_{\ontop{0 < d < X}{d \equiv a \bmod N}} \chi_d(n)
    \sim
    \chi_a(g)  \phi(\square)
    \frac{1}{\phi(4N\square/Q)} \frac{X}{Q} \frac{6}{\pi^2}
    h_2(a,N)
    \prod_{p|2N\square} \frac{p}{p+1}.
\end{equation}
Here $h_2(a,N)$ is determined according to Table~\ref{tab:h_2}.
Consequently,
for the family ${\mathcal F}_{\text{d}}(a,N,\pm)$,
\begin{equation} 
    \label{eqn:chidaverage}
    \langle  \chi_d(n) \rangle = 
    \begin{cases} \chi_a(g) \prod_{p| \square} (1+\frac1p)^{-1}
    & n= g \square \\
    0& \text{otherwise}.
    \end{cases} 
\end{equation}
\end{theorem}

\begin{table}[ht]
\begin{center}
\begin{tabular}{|c|c|c|}
$\beta$ & $a$ & $h_2(a,N)$ \cr \hline \hline 
0        & $a \in {\mathbb Z}$                   & 3/2 \cr \hline 
1        & $a \equiv 0 \bmod 2$               & 1   \cr
         & $a \equiv 1 \bmod 2$               & 2   \cr \hline 
2        & $a \equiv 0 \bmod 4$               & 2   \cr
         & $a \equiv 1 \bmod 4$               & 4   \cr
         & $a \equiv 2,3 \bmod 4$             & 0   \cr \hline
3        & $a \equiv 0,4 \bmod 8$             & 2   \cr
         & $a \equiv 1,5 \bmod 8$             & 4   \cr
         & $a \equiv 2,3,6,7 \bmod 8$         & 0   \cr \hline
$\geq 4$ & $a \equiv 1,5,8,9,12,13 \bmod 16$  & 4   \cr
         & otherwise                          & 0   \cr \hline 
\end{tabular}
\end{center}
\caption{ 
The function $h_2(a,N)$ that appears in Theorem~\ref{thm:funddisc},
where $N=2^\beta N_0$, with $N_0$ odd.
}\label{tab:h_2}
\end{table}

\begin{proof}
We first outline the proof for (\ref{eqn:fund discr progressions}). One can count 
odd fundamental discriminants $|d|< X$ by using the Dirichlet series
\begin{equation}
    \sumstar_{\text{$d$ odd}} \frac{1}{|d|^s} = 
    \prod_{\text{$p$ odd}} \left( 1+\frac{1}{p^s}\right) =
    \frac{\zeta(s)}{\zeta(2s)} \left( 1+\frac{1}{2^s}\right)^{-1}.
\end{equation}
As in the proof of the prime number theorem, the main contribution comes from the pole at $s=1$,
and one has 
\begin{equation}
    \sumstar_{\ontop{|d| < X}{\text{$d$ odd}}} 1
    \sim 
    \frac{4}{\pi^2} X.
\end{equation}
Next, assume that $N$ is odd and $(a,N)=1$. To count odd fundamental discriminants in arithmetic progression,
$d \equiv a \mod N$, one imitates Dirichlet's theorem for primes in arithmetic progression,
looking at linear combinations involving Dirichlet characters mod~$N$~of
\begin{equation}
    \sumstar_{\text{$d$ odd}} \frac{\chi(d)}{|d|^s} = 
    \prod_{\text{$p$ odd}} \left( 1+\frac{\chi(p)}{p^s}\right).
\end{equation}
If one wishes to further specify $d>0$ or $d<0$, one can restrict to 
$|d| \equiv 1 \mod 4$ or $|d| \equiv 3 \mod 4$ respectively, with $\chi$ ranging
over Dirichlet character mod $4N$. The main contribution comes from the
trivial character whose corresponding Dirichlet series is
\begin{equation}
    \label{eqn:trivial character fund discr}
    \sumstar_{\text{$d$ odd}} \frac{\chi_0(d)}{|d|^s} = 
    \prod_{\ontop{\text{$p$ odd}}{p \nmid N}} \left( 1+\frac{1}{p^s}\right)=
    \frac{\zeta(s)}{\zeta(2s)} \prod_{p \mid 2N} \left( 1+\frac{1}{p^s}\right)^{-1},
\end{equation}
and whose main pole is at $s=1$. Therefore, for $N$ odd and $(a,N)=1$, we have
\begin{equation}
    \sumstar_{\ontop{0 < d < X}{\ontop {d \equiv a \bmod N}{\text{$d$ odd}}}} 1
    \sim \frac{1}{\phi(4N)} X \frac{6}{\pi^2}
    \prod_{p|2N} \frac{p}{p+1},
\end{equation}
with the same result for $-X < d < 0$. 

Next, for $N$ odd and $(a,N)=Q>1$, one can write, for $d \equiv a \bmod N$,  $d=d_1 Q$.
Apply the above method to $d_1$ with $0 < d_1 < X/Q$, $d_1 \equiv (a/Q) \bmod (N/Q)$, $d_1$ odd,
and, because $d$ is squarefree, the extra condition that $(d_1,Q)=1$. 
Because of this last condition, the Euler product 
that we need to take in (\ref{eqn:trivial character fund discr}) is not just over
odd $p \nmid (N/Q)$ but also $p \nmid Q$, i.e. it is still 
\begin{equation}
    \prod_{\ontop{\text{$p$ odd}}{p \nmid N}} \left( 1+\frac{1}{p^s}\right).
\end{equation}
Hence, if $(a,N)= Q$,
\begin{equation}
    \sumstar_{\ontop{0 < d < X}{\ontop{d \equiv a \bmod N}{\text{$d$ odd}}}} 1
    \sim \frac{1}{\phi(4N/Q)} \frac{X}{Q} \frac{6}{\pi^2}
    \prod_{p|2N} \frac{p}{p+1},
\end{equation}
and the same result for $-X < d < 0$. 

Finally, we wish to take into account even $d$. The set of even fundamental discriminants
consists of $-4$ and $\pm 8$ times the odd fundamental discriminants. 

Again, assume $N$ is odd.
One can count discriminants, $d \equiv a \bmod N$, lying in the interval $(0,X)$
by counting odd discriminants lying in $(0,X)$, together with odd discriminants in $(-X/4,0)$,
$(0,X/8)$ and $(-X/8,0)$. Overall, this gives the same asymptotics as before, but with an extra
factor of $(1+1/4+2/8)=3/2$. This accounts for line 1 in Table~\ref{tab:h_2}. The other lines
in the table can be obtained by similar considerations.

We now apply (\ref{eqn:fund discr progressions}) to obtain (\ref{eqn:chi over fund discr})
and (\ref{eqn:chidaverage}). Consider
\begin{equation}
\label{eqn:fund discr progr kronecker ratio}
    \frac{\displaystyle
        \sumstar_{\ontop{0 < d < X}{d \equiv a \bmod N}} \chi_d(n)
    }
    {\displaystyle
        \sumstar_{\ontop{0 < d < X}{d \equiv a \bmod N}} 1
    }
\end{equation}
(the following analysis also holds for $-X < d <0$).

Write $N=N_1^{r_1}\cdot\ldots\cdot N_m^{r_m}$, the prime factorization of $N$,
and let $g=N_1^{u_1}\cdot\ldots\cdot N_m^{u_m}$. Then
\begin{equation}
    \chi_d(n) = \chi_d(N_1)^{u_1} \ldots \chi_d(N_m)^{u_m} \chi_d( \square).
\end{equation}
Now, if $N_i$ is an odd prime, $\chi_d(N_i) = \chi_a(N_i)$, since $d \equiv a \bmod N$, and so
$d \equiv a \bmod N_i$. If $N_i=2$ we need to be careful because $\chi_d(2)$ is
periodic mod 8. Now we are assuming that 
if $N$ is even it be a least divisible by 8, i.e. that $d \equiv a$ mod $8$,
and thus that $\chi_d(2) = \chi_a(2)$. 

Therefore
\begin{equation}
     \chi_d(n) = \chi_a(g) \chi_d(\square) =
     \begin{cases}
         \chi_a(g) &\text{if $(d,\square)=1$} \\
         0 &\text{otherwise,}
     \end{cases}
\end{equation}
and one gets that (\ref{eqn:fund discr progr kronecker ratio}) equals
\begin{equation}
\label{eqn:fund discr progr kronecker ratio 2}
    \chi_a(g)
    \frac{\displaystyle
        \sumstar_{\ontop{\ontop{0 < d < X}{d \equiv a \bmod N}}{(d,\square)=1}} 1
    }
    {\displaystyle
        \sumstar_{\ontop{0 < d < X}{d \equiv a \bmod N}} 1
    }.
\end{equation}
Since $(\square,N)=1$, the sum in the numerator can be split into sums $d \bmod N \square$.
Naively, one expects to have $\phi(\square)$ sums, one for each residue class $(d,\square)=1$.
However, if $\square$ is even, then only half of these residue classes, namely 
those that have $d \equiv 1 \bmod 4$, contain fundamental discriminants, so one only gets
$\phi(\square)/2$ sums. We thus consider the case that $\square$ is odd separately from
the case that it is even. Both cases end up giving the same answer.

Assume that $\square$ is odd. To apply our formula (\ref{eqn:fund discr progressions})
to each of the $\phi(\square)$ residue classes $\bmod N \square$, one needs to compute
the various components that go into the formula.

Given $d \equiv a \bmod N$ 
and $d \equiv b \bmod \square$, one has via the chinese remainder 
theorem $d \equiv \tilde{a} \bmod N \square$. Now, $Q=(a,N)=(d,N)$, and $(d,\square) = 1$,
so $(\tilde{a},N\square)=(d,N\square)=Q$.

One also needs to evaluate $h_2(\tilde{a}, N \square)$. Let $N=2^\beta N_0$, with
$N_0$ odd. Now $\square$ is odd, and so $h_2(\tilde{a},N\square)$ only depends on
$\tilde{a} \bmod 2^\beta$, but this is determined by $a \bmod N$. So 
$h_2(\tilde{a},N\square)=h_2(a,N)$.
Therefore, the numerator of (\ref{eqn:fund discr progr kronecker ratio 2}) is 
asymptotically
\begin{equation}
    \chi_a(g)  \phi(\square)
    \frac{1}{\phi(4N\square/Q)} \frac{X}{Q} \frac{6}{\pi^2}
    h_2(a,N)
    \prod_{p|2N\square} \frac{p}{p+1}.
\end{equation}
Canceling factors appearing in the asymptotics (\ref{eqn:fund discr progressions})
of the denominator of
(\ref{eqn:fund discr progr kronecker ratio 2}) we get
\begin{equation}
    \chi_a(g)  \prod_{p| \square} \frac{p}{p+1}.
\end{equation}

If $\square$ is even, write $\square = 2^\lambda \square_0$, with
$\lambda \geq 2$, and $\square_0$ odd.
Now, $(d,\square)=1$, so $d$ is odd. In all cases, according to Table~\ref{tab:h_2},
$h_2(\tilde{a},N\square)$ is therefore 4. Furthermore, as in the odd case,
$(\tilde{a},N\square)=Q$. 

Hence, one gets asymptotically 
for the numerator of (\ref{eqn:fund discr progr kronecker ratio 2})
\begin{equation}
    \chi_a(g) \frac{\phi(\square)}{2}
    \frac{1}{\phi(4N\square/Q)} \frac{X}{Q} \frac{6}{\pi^2}
    4 \prod_{p|2N\square} \frac{p}{p+1}.
\end{equation}
Since $\square$ is even, $N$ is odd.
Hence $h_2(a,N)=3/2$, and the denominator of 
(\ref{eqn:fund discr progr kronecker ratio 2}) is asymptotically
\begin{equation}
    \frac{1}{\phi(4N/Q)} \frac{X}{Q} \frac{6}{\pi^2}
    \frac{3}{2}
    \prod_{p|2N} \frac{p}{p+1}.
\end{equation}
Canceling numerator and denominator, taking special care for
powers of 2 appearing in $\square$, we get 
\begin{equation}
    \chi_a(g) \frac{2}{3}
    \prod_{p|\square_0} \frac{p}{p+1} = \chi_a(g) \prod_{p|\square} \frac{p}{p+1}.
\end{equation}

\end{proof}


\noindent {\bf 4. The family of coefficients of holomorphic newforms.} 
\begin{equation}{\mathcal F}_{\text{mod}}(k,q)=\{f(n)= \lambda_f(n):
\sum n^{(k-1)/2}\lambda_f(n) \in H_k(q)\},  
\end{equation}
where $H_k(q)$  the set of newforms in  $S_k(\Gamma_0(q))$.
A good reference for these functions is Iwaniec\cite{Iw1}.
In this family the parameter tending to infinity can be either
$k$, or $q$, or some combination.

The
orthogonality relation here is somewhat subtle, and in fact there are 
two natural ways to average over these characters.  In both cases the
starting point is the Hecke relation
\begin{equation}
\lambda_f(m)\lambda_f(n)=\sum_{ \ontop{ d|m, d|n}{(d,q)=1} } \lambda_f(mn/d^2),
\end{equation}
which imply that any product
\begin{equation}\lambda_f(m_1)\cdots \lambda_f(m_k)
\end{equation}
can be expressed as a linear combination
\begin{equation}\label{eqn:lambdajsum}
\sum_{j\ge 1} b_j \lambda_f(j)
\end{equation}
for some integers $b_j$, and in fact only for $j$ a prime power.
Thus, we need only determine the expected value of
$\lambda_f(p^j)$.

If one averages over $H_k(q)$ in the most straightforward way,
then for $p\nmid q$,
\begin{equation}
\langle \lambda_f(p^j) \rangle = 
\begin{cases}
p^{-j/2} & j \text{ \ even} \\
0        & j \text{ \ odd},
\end{cases}
\end{equation}
and more generally, if $(n,q)=1$,
\begin{equation}
\langle \lambda_f(n) \rangle = 
\begin{cases}
n^{-1/2} & n=\square \\
0        & \text{otherwise},
\end{cases}
\end{equation}
This follows from the Selberg trace formula.
However, if one averages with respect to a weighting by the
Petersson norm:
\begin{equation}
\sumhar_{f\in H_k(q)} * \ \ =\ \sum_{f\in H_k(q)} */\langle
f,f\rangle ,
\end{equation}
then
\begin{equation}
\langle \lambda_f(p^j) \rangle = 
\begin{cases}
1 & j =0 \\
0        & \text{otherwise},
\end{cases}
\end{equation}
and more generally, if $(n,q)=1$,
\begin{equation}
\langle \lambda_f(n) \rangle = 
\begin{cases}
1 & n =1 \\
0        & \text{otherwise}.
\end{cases}
\end{equation}
This follows from the Petersson formula (see \cite{Iw1}), if $(mn,q)=1$,
\begin{equation}\label{eqn:petersson}
\sumhar_{f\in H_k(q)} \lambda_f(m)\lambda_f(n)=
\delta(m,n)+2\pi i^k  \sum_{c=1}^\infty \frac{S(m,n;cq)J_{k-1}(4\pi\sqrt{mn}/cq)}{cq}  .
\end{equation}
Here $J_k$ is the Bessel function
and
\begin{equation}
S(m,n;c)=\sum_{ad=1\mod c}e\left(\frac {ma+nc}{c}\right)
\end{equation}
is the Kloosterman sum.  Since the Petersson weighting leads to a somewhat
simpler expression, we will consider that weighting in our example.
When passing from the Petersson formula 
to the expected value, using the Weil bound for the Kloosterman sum
and the fact that $J_k$ has a $k$th order zero at~$0$,
we see that for fixed $m, n$ the sum on
the right side of (\ref{eqn:petersson})
vanishes as $k\to\infty$ or~$q\to\infty$.

Let
\begin{equation}\delta(m_1,\dots,m_k) = 
\langle \lambda_f(m_1)\cdots \lambda_f(m_k) \rangle.
\end{equation} 
So in the Petersson weighting, $\delta(m_1,\dots,m_k)$ is the 
the coefficient $b_1$ of $\lambda_f(1)=1$ in~(\ref{eqn:lambdajsum}).
One can use the Hecke relations to show by induction that $\delta$ is multiplicative
in the sense of~\ref{eqn:deltamultiplicative}. 
Thus, we only need to know $\delta$ on prime powers.

 \begin{lemma}\label{thm:chebyshevweights}
With respect to the Petersson weighting, if $p\nmid q$ then 
\begin{eqnarray}\label{eqn:petersondelta} 
\delta(p^{m_1},\dots,p^{m_k})
&  = &\frac{2}{\pi}\int_0^\pi
\sin^2\theta
\prod_{j=1}^k\frac{\sin(m_j+1)\theta}{\sin \theta}\,d\theta \nonumber \\&  = & \frac{2}{\pi}\int_0^\pi \sin^2\theta \prod_{j=1}^k
\frac{e^{i(m_j+1)\theta}-e^{-i(m_j+1)\theta}}{e^{i\theta}-e^{-i\theta}}
\,d\theta.
\end{eqnarray}
For the unweighted sum we have
\begin{equation}
\delta(p^{m_1},\dots,p^{m_k})
  = \frac{4}{\pi}\int_0^\pi
\frac{\sin^2\theta}{1-\frac{2\cos\theta}{\sqrt{p}} +\frac1p}
\prod_{j=1}^k\frac{\sin(m_j+1)\theta}{\sin \theta}\,d\theta .
\end{equation}
If $p|q$ then $\delta(p^{m_1},\dots,p^{m_k})=0$ unless $m_1=\cdots = m_k=0$.
\end{lemma}

\begin{proof} We only give the details for~(\ref{eqn:petersondelta}).
Beginning from 
\begin{equation}\mathcal{L}_p(x)=\sum_{j=0}^\infty
\lambda_f(p^j)x^j=\left(1-e^{i\theta_{p,f}}x\right)^{-1}
\left(1-e^{-i\theta_{p,f}}x\right)^{-1}
\end{equation}
we have
$\lambda_f(p^j)=\frac{\sin (j+1)\theta_{f,p}}{\sin
\theta_{f,p}}=U_j(\cos \theta_{f,p})$ where $U_j$ is the usual
Tchebychev polynomial. 
Then $\delta(p^{m_1},\dots, p^{m_k})=c_0$ where 
\begin{equation}
U_{m_1}U_{m_2}\dots U_{m_k}=\sum_{e\ge 0} c_e U_e.
\end{equation}
If we evaluate both sides of this equation at $\cos \theta$ and
integrate from 0 to $\pi$ with respect to the measure $\frac 2 \pi
\sin^2 \theta \,d\theta$, then the result follows from the
orthogonality of the Tchebychev polynomials with respect to this
measure.
\end{proof}

\noindent {\bf 5. The family of coefficients of Maass newforms.}
\begin{equation}{\mathcal F}_{M}(q)=\{f(n)= \lambda_f(n):
\sqrt{y} \sum \lambda_f(n) K_{iR}(2\pi |n|y) e^{2\pi i n x} \in H(q)\},  
\end{equation}
where $H(q)$  the set of Maass newforms on  $\Gamma_0(q)$.
A good reference for these functions is Iwaniec~\cite{Iw2}.
The orthogonality relation is derived from the Kuznetsov trace
formula.  See Chapter~9 of~\cite{Iw2}.

\subsection{ Families of $L$-functions}\label{sec:familiesofLfunctions}

We use a family of characters to create a family of $L$-functions
in the following manner. 

Begin with a fixed primitive $L$-function
\begin{equation}\mathcal{L}(s)=\sum_{n=1}^\infty \frac{a_n}{n^s}=
\prod_p \mathcal{L}_p\left(\frac1{p^s}\right), \qquad (\sigma>1) .
\label{eqn:generalL}
\end{equation}
  We assume that
\begin{equation}\label{eqn:generalLp}
\mathcal{L}_p(x)=\sum_{n=0}^\infty a_{p^n}x^n 
=\prod_{j=1}^w (1-\gamma_{p,j} x)^{-1},
\end{equation}
where $w$ is the degree of $\mathcal{L}$ and where
$|\gamma_{p,j}|=0$ or 1. Assume $\mathcal{L}(s)$ satisfies
functional equation
\begin{equation}
\mathcal{L}(s)=
\varepsilon {\mathcal X}(s) \overline{\mathcal{L}}(1-s), 
\label{eqn:generalLfe}
\end{equation}
as described in Section~\ref{sec:propofLfunctions}.

We create a family of $L$-functions by twisting $\mathcal{L}$
by a family of characters.
Let $\mathcal F=\{f\}$ be a family of
characters, with the properties described 
in Section~\ref{sec:familiesofcharacters}.
The twist of $\mathcal{L}$ by $f$  is denoted by $\mathcal{L}(s,f)$ and is
given by a Rankin-Selberg convolution:
\begin{equation} \label{eqn:rankinselberg}
\mathcal{L}(s,f)=
\prod_p\prod_{i=1}^v\prod_{j=1}^w(1-\beta_{p,i}\gamma_{p,j}/p^s)^{-1}=\sum_{n=1}^\infty
\frac{a_n(f)}{n^s}. 
\end{equation}
Note that if $w=1$ or $v=1$, as will be the case in our detailed
examples,
\begin{equation}
\mathcal{L}(s,f)=\sum_{n=1}^\infty
\frac{a_n f(n)}{n^s}.
\end{equation}

We require that $\mathcal{L}(s,f)$ is an $L$-function.  That is, our family
of $L$-functions must consist of $L$-functions!  In particular, $\mathcal{L}(s,f)$
 satisfies a functional equation
\begin{equation}\mathcal{L}(s,f)=
\varepsilon_f \mathcal{X}_f(s) \overline{\mathcal{L}}(1-s,f),
\label{eqn:genfe}
\end{equation}
as described in Section~\ref{sec:propofLfunctions}.
As part of our definition of ``family'', 
we make a restrictive, but natural, assumption on $\mathcal{X}_f$.
We have $\mathcal{X}_f(s)= \overline{\gamma_f}(1-s)/\gamma_f(s)$ where
\begin{equation}
\gamma_f(s)=Q_f^s\prod_{j=1}^w \Gamma(\tfrac12 s+\mu_{j,f}) .
\end{equation}
We assume that $w$ is constant, and each of $Q_j$, $\Re(\mu_{j,f})$,
and $\Im(\mu_{j,f})$ is a monotonic function of the conductor~$c(f)$.
In practice this will mean that each of those quantities will
either be constant or will be tending to infinity with the conductor.

For example, the collection of all real primitive Dirichlet $L$-functions 
$L(s,\chi_d)$
is not a family, because for $d>0$ we have $\mu_{1,d}=0$ and 
for $d<0$ we have  $\mu_{1,d}=1$, so $\mu_{1,d}$ is not
a monotonic function of~$c(\chi_d)$.
So we consider these two families separately.

Finally, we will make use of an approximate functional equation
of shape
\begin{equation}\label{eqn:generalappfe}
\mathcal{L}(s,f)= \sum \frac{a_n(f)}{n^s}+
\varepsilon_f \mathcal{X}_f(s) \sum
\frac{\overline{a_n(f)}}{n^{1-s}} + remainder.
\end{equation}
Note: we are not claiming that the ``remainder'' in the above equation
is small, nevertheless we will ignore the remainder in our calculations.

\section{A recipe for conjecturing moments, with examples}\label{sec:recipe}

We give a general recipe for conjecturing the moments of a primitive family
of $L$-functions, and then apply the recipe to several interesting examples.

\subsection{ The general recipe
}\label{sec:generalrecipe}

Suppose $\mathcal L$ is an $L$-function and $f$ is a character with conductor $c(f)$,
as described in Section~\ref{sec:families}.  So
\begin{equation}Z_\mathcal{L}(s,f)=\varepsilon_f^{-\frac 12}\mathcal{X}_f(s)^{-\frac 12}\mathcal{L}(s,f),
\label{eqn:ZLf}
\end{equation}
which satisfies the functional equation
\begin{equation}Z_\mathcal{L}(s,f)=  \overline{Z_\mathcal{L}}(1-s,f),
\end{equation}
so $Z_\mathcal{L}(s,f)$ is real on the $\frac12$-line. 
Note that the square root of $\varepsilon_f$ involves a choice which
must be made consistently.

We consider the
moment
\begin{equation}
\sum_{f\in \mathcal F} Z_\mathcal{L}(\tfrac12 +\alpha_1,f)\dots
Z_\mathcal{L}(\tfrac12 +\alpha_k,f) g(c(f))
\end{equation}
where $g$ is a suitable test function. The recipe below also applies 
to averages of products of 
$\mathcal{L}(\tfrac12 +\alpha,f)$.  The sum is an integral
when ${\mathcal F}={\mathcal F}_t$.  

Here is a recipe for conjecturing a formula for the above moment:

\begin{enumerate}
\item{} Start with a product of $k$ shifted $L$-functions:
\begin{equation}
Z_f(s,\alpha_1,\dots,\alpha_{k}) =
Z_\mathcal{L}(s +\alpha_1,f)\dots
Z_\mathcal{L}(s +\alpha_k,f).
\end{equation}
As we will demonstrate in our examples, the recipe applies to
the $Z$-function as well as the $L$-function.

\item{} Replace each $L$-function with the two terms from its
approximate functional equation (\ref{eqn:generalappfe}), ignoring the remainder term.
Multiply out the resulting
expression to obtain $2^{k}$ terms.  Write those terms as
\begin{equation}
(\text{product of $\varepsilon_f$ factors})
(\text{product of $\mathcal{X}_f$ factors})
\sum_{n_1,\dots,n_k} (\text{summand}) .
\end{equation}

\item{} Replace each product of $\varepsilon_f$-factors 
by its expected value when averaged over the family.

\item{} Replace each summand by its expected value when averaged over the family.

\item{} Complete the resulting sums, and call the total
$M(s,\alpha_1,\dots,\alpha_{2k})$.

\item{} The conjecture is
\begin{equation}
\sum_{f\in \mathcal F} Z_f(\tfrac12,\alpha_1,\dots,\alpha_{2k}) g(c(f))
=
\sum_{f\in \mathcal F} M_f(\tfrac12,\alpha_1,\dots,\alpha_{2k})
(1+O(e^{(-\frac12 + \varepsilon)c(f)})) g(c(f)),
\end{equation}
for all $\varepsilon>0$, where $g$ is a suitable weight function.
\end{enumerate}

In other words, $Z_f(s,\alpha)$ and $M_f(s,\alpha)$ have the same value
distribution if averaged over a sufficiently large portion of the family.
Note that the dependence of $M_f$ on $f$ only occurs in the product
of $\mathcal{X}_f$ factors.

As we mentioned earlier, some of the individual steps in this 
recipe cannot be rigorously justified.  Only by using the entire
recipe does one arrive at a reasonable conjecture.  In particular,
we ignore off-diagonal terms which actually
make a contribution.  However, comparison with examples in the literature,
random matrix moments, and
numerical data, suggests that the various errors in our recipe all
cancel.  The underlying cause for this remains a mystery.

We will apply the recipe to several examples, but first we
do the initial steps of the recipe in some generality.

For each
$Z_\mathcal{L}$ substitute the expression in (\ref{eqn:ZLf}). 
After replacing
each $\mathcal{L}(s,f)$ by its approximate functional equation 
(\ref{eqn:generalappfe}),
multiply out the product. A typical term is
a product of $k$ sums arising from either the first piece or the
second piece of the approximate functional equation. Consider a
term where we have $\ell$ factors from the first piece of an
approximate functional equation and $k-\ell$ factors from the
second piece.  To take one specific example, suppose it is the
first $\ell$ factors from which we choose the first piece of the
approximate functional equation, and the last $k-\ell$ factors
from which we take the second piece of the approximate functional
equation: 
\begin{eqnarray}
&   & \varepsilon_f^{-\frac \ell 2} \mathcal{X}_f(s +\alpha_1)^{-\frac 12}\cdots 
\mathcal{X}_f(s +\alpha_\ell)^{-\frac 12}
\sum_{n_1}\frac {a_{n_1}(f)}{n_1^{s +\alpha_1}}\cdots 
\sum_{n_\ell}\frac {a_{n_\ell}(f)}{n_\ell^{s +\alpha_\ell}}\nonumber \\
&   &\phantom{xxxxxx} \times \varepsilon_f^{\frac12(k-\ell)} 
\mathcal{X}_f(s +\alpha_{\ell+1})^{\frac12} 
\cdots 
\mathcal{X}_f(s +\alpha_{k})^{\frac12}
\sum_{n_{\ell+1}}\frac
{\overline{a_{n_{\ell+1}}(f)}}{n_{\ell+1}^{1-s -\alpha_{\ell+1}}} 
\cdots \sum_{n_{k}}\frac
{\overline{a_{n_{k}}(f)}}{n_{k}^{1-s -\alpha_{k}}} .
\end{eqnarray}
Rearranging this expression and using the fact that
$\mathcal{X}_f(s)=\mathcal{X}_f(1-s)^{-1}$, we have
\begin{eqnarray}
&   & 
\varepsilon_f^{\frac{k}{2}-\ell}
\prod_{j=1}^\ell \mathcal{X}_f(s +\alpha_j)^{-\frac12}
\prod_{j=\ell+1}^k \mathcal{X}_f(1-s -\alpha_j)^{-\frac12}
\cr
&   & \phantom{XXXXXXXXXXXxxxx} \times
\sum_{n_1,\dots, n_k}
    \frac{
a_{n_1}(f)\cdots a_{n_\ell}(f)
\overline{a_{n_{\ell+1}}(f)\cdots a_{n_k}(f)}
}{n_1^{s +\alpha_1}\cdots n_\ell^{s +\alpha_\ell}
n_{\ell+1}^{1-s -\alpha_{\ell+1}} \cdots n_k^{1-s-\alpha_k}}
\, 
.
 \end{eqnarray}
A little trick:  since we will eventually set $s=\frac12$, we 
replace the above expression by
\begin{eqnarray}\label{eqn:generalfirstell}
&   & 
\varepsilon_f^{\frac{k}{2}-\ell}
\prod_{j=1}^\ell \mathcal{X}_f(s +\alpha_j)^{-\frac12}
\prod_{j=\ell+1}^k \mathcal{X}_f(s -\alpha_j)^{-\frac12}
\cr
&   & \phantom{XXXXXXXXXXXxxxx} \times
\sum_{n_1,\dots, n_k}
    \frac{
a_{n_1}(f)\cdots a_{n_\ell}(f)
\overline{a_{n_{\ell+1}}(f)\cdots a_{n_k}(f)}
}{n_1^{s +\alpha_1}\cdots n_\ell^{s +\alpha_\ell}
n_{\ell+1}^{s -\alpha_{\ell+1}} \cdots n_k^{s-\alpha_k}}
\, 
.
 \end{eqnarray}
It is expression (\ref{eqn:generalfirstell}), and the corresponding
pieces from the other terms when multiplying out the approximate functional
equation, which will appear in the final conjecture, evaluated
at~$s=\frac12$.

Now consider the product of $\varepsilon_f$ factors $\varepsilon_f^{\frac{k}{2}-\ell}$, 
which according
to the recipe should be replaced by its expected value.
An important issue is the choice of the square root.
We believe that there is a natural choice of $\varepsilon_f^{\frac{1}{2}}$
so that the following hold.
\begin{itemize}
\item{}Unitary case: the $\varepsilon_f$ are
uniformly distributed on the unit circle, and 
$\langle \varepsilon_f^{\frac{k}{2}-\ell}\rangle=0$ unless
$\frac{k}{2}-\ell=0$.  In particular, $k$ must be even.
There will be $\binom{k}{k/2}$ terms in the final answer.

\item{}Orthoconal case:  $\varepsilon_f$ is constant ($1$ or $-1$) over the family,
or $\varepsilon_f=1$ for approximately half the~$f$
and $\varepsilon_f=-1$ for the other half. We have
$\langle \varepsilon_f^{\frac{k}{2}-\ell}\rangle=0$ unless
$\frac{k}{2}-\ell$ is even.  In particular, $k$ must be even and
there will be $2^{k-1}$ terms in the final answer.

\item{}Symplectic case: $\varepsilon_f=1$ for all~$f$, 
and $\langle \varepsilon_f^{\frac{k}{2}-\ell}\rangle=1$ for all $k$ and $\ell$.
There is no restriction and there will be $2^k$ terms in the final answer.
\end{itemize}

Note that if we are considering the $L$-function, instead of the $Z$-function,
then the issue of $\varepsilon_f^{\frac{1}{2}}$ does not arise and
the calculation is somewhat easier.  See 
\eqref{eqn:cuspsumterm} and the discussion following.
Also note that in the Unitary and Orthogonal cases, odd powers of the
$Z$-function will average to zero,
while odd powers of the $L$-function will not.

The recipe now says to
replace the summand by its expected value when averaged over the family.
That is, we replace
\begin{equation}
a_{n_1}(f)\dots a_{n_\ell}(f)\overline {a_{n_{\ell+1}}(f)\dots a_{n_k}(f)} 
\end{equation}
by its expected value when averaged over the family.  
In practice, this will be of the form
\begin{equation}
c(\mathcal F) \delta_\ell(n_1,\dots,n_k)
\end{equation}
where $c(\mathcal F)$ depends only on the family, and where the $\delta_\ell$ are
multiplicative functions, i.e.
\begin{equation}
\delta_\ell(m_1n_1,\dots ,m_kn_k)=\delta_\ell(m_1,\dots,m_k)\delta_\ell(n_1,\dots,n_k)
\end{equation}
whenever $(m_1\cdots m_k,n_1\cdots n_k)=1$. 

The final step is to extend the range of summation.  This produces one 
term in the conjecture.  By considering the other terms when
multiplying out the approximate functional equations, one arrives
at a conjecture for the original mean value.

Although the above steps have produced an answer, it is not
written in a particularly usable form.  There are three more steps
to put the conjecture in the form of Conjecture~\ref{thm:zeta2kconjecture}:
writing the main terms as an Euler product, identifying the polar part,
and expressing the combinatorial sum as a multiple integral.

Since the  $\delta_\ell$  are multiplicative, we can write the
main term as an
Euler product. Specifically,
\begin{equation}
\sum_{n_1,\dots ,n_k=1}^\infty
\frac {\delta_\ell(n_1,\dots,n_k)}{n_1^{s+\alpha_1} \cdots
n_k^{s+\alpha_k}} =\prod_p \sum_{e_1,\dots ,e_k=0}^\infty
\frac{\delta_\ell(p^{e_1},\dots,p^{e_k})}{p^{e_1(s+\alpha_1)+\cdots
+e_k(s +\alpha_k)}}
\end{equation}
assuming that $\Re s$ is sufficiently large. 

Next we determine the leading order poles.  It usually turns out
that 
$\delta_\ell(p^{e_1},\dots,p^{e_k})=0$ when 
$\sum e_j=1$. Thus, the first poles come from those terms
where $\sum e_j=2$.  This happens in two ways, so the rightmost poles
are the same as the poles of
\begin{equation}
\prod_{1\le i<j\le k} \prod_p 
\(1+\frac{\delta_{\ell,i,j}(p,p)}{p^{2s+\alpha_i+\alpha_j}} \)
\times
\prod_{j=1}^k \prod_p \(1+\frac{\delta_{\ell,j}(p^2)}{p^{2s+2\alpha_j}}\)
.
\end{equation}
In practice the first factor has simple poles at 
$\frac12-\frac12(\alpha_i+\alpha_j)$,
and the second factor is either regular in a neighborhood of $\sigma=\frac12$,
or else it has a simple pole at $s=\frac12-\alpha_j$.  
Accordingly, we factor out either
\begin{equation}
\prod_{1\le i<j\le k} \zeta(2s+\alpha_i+\alpha_j)
\ \ \ \ \ \ \text{ or }
\ \ \ \ \ \ 
\prod_{1\le i\le j\le k} \zeta(2s+\alpha_i+\alpha_j)
.
\end{equation}
The remainder is the $A_k$ in our conjectures, and it is 
regular in a neighborhood of~$\sigma=\frac12$.

Having identified the polar part of our main terms, we can apply the lemmas
in Section~\ref{sec:concisesums} to express the sum of terms
as a contour integral.  The result is an expression similar
to Conjecture~\ref{thm:zeta2kconjecture}.

We have already seen this procedure in Section~\ref{sec:zetameanvalues} for the case of
mean values of the zeta-function. In the following sections we
carry out example calculations for families of each of the three
symmetry types.

\subsection{ Unitary: moments of a primitive $L$-function
}\label{sec:unitarymean}

The recipe for mean values in Section~\ref{sec:therecipe} is a special case
of the general recipe.  To see this, note that if
$f_t\in \mathcal{F}_t$ then $f_t(n)=n^{-it}$,
so $\mathcal{L}(s,f_t)=\mathcal{L}(s+it)$.  From the
functional equation 
$\mathcal{L}(s) = \varepsilon \mathcal{X}(s) \overline{\mathcal{L}}(1-s)$ 
we obtain the functional equation
\begin{equation}
\mathcal{L}(s,f_t) = \varepsilon_t \mathcal{X}_t(s) \overline{\mathcal{L}}(1-s,f_t),
\end{equation}
where
\begin{equation}
\varepsilon_t = {\varepsilon}{ \mathcal{X}(\tfrac12+it)},
\text{\ \ \ \ \ \ \ \ and \ \ \ \ \ \ \ \ \ }
\mathcal{X}_t(s) = \frac{\mathcal{X}(s+it)}{ \mathcal{X}(\tfrac12+it)} .
\end{equation}
Note that these satisfy the requirements 
$|\mathcal{X}_t(\frac12 + iy)| = 1$ for $y$ real, with 
$\mathcal{X}_t(\frac12) = 1$ and 
$|\varepsilon_t|=1$.  
Also note that the log conductor of $\mathcal{L}(s,f_t)$, defined
as $|(\varepsilon_t\mathcal{X}_t)^\prime(\frac12)|$, equals $|\mathcal{X}^\prime(\frac12 +it)|$,
in agreement with the usual notion of conductor in $t$-aspect.

Replacing the product of 
$\varepsilon_t$-factors by their expected value is the same
``keep the terms where the product of the $\chi$-factors is not oscillating.''
Thus, after multiplying out
the approximate functional equations there will be
$\binom{2k}{k}$ terms which contribute.  
In each of those terms replacing the
summand by its expected value is the same as ``keeping the diagonal.''
Thus, we arrive at the same conjecture as before.


\subsection{ Unitary: all Dirichlet $L$-functions
}\label{sec:unitary}

We apply our recipe to conjecture the average
\begin{equation}
 \sumstar_{\ontop{\chi\mod q}{\chi\ \text{even or odd}}}
Z_\chi(\tfrac12;\alpha_1,\dots,\alpha_{2k}) ,
\end{equation}
where the sum is over either the even or the odd
 primitive Dirichlet characters $\mathstrut\mod q$
and
\begin{equation}
Z_\chi(s;\alpha_1,\dots,\alpha_{2k})=
 Z(s +\alpha_1,\chi)\dots
Z(s +\alpha_{2k},\chi)
.
\end{equation}
Here $Z(s +\alpha_k,\chi) = (\varepsilon_\chi X_\chi (s) )^{-\frac12} L(s, \chi)$
where $L(s, \chi)=\varepsilon_\chi X_\chi(s) L(1-s,\chi)$.
Note that $\varepsilon_\chi=\tau(\chi)q^{-\frac12}$,
which is uniformly distributed on the unit circle.

Following the general discussion in Section~\ref{sec:generalrecipe},
equation (\ref{eqn:generalfirstell}) specializes in this case~to
\begin{eqnarray}
&   & 
\varepsilon_\chi^{{k}-\ell} 
\prod_{j=1}^\ell {X}_\chi(\tfrac12 +\alpha_j)^{-\frac12}
\prod_{j=\ell+1}^{2k} {X}_\chi(\tfrac12 -\alpha_j)^{-\frac12}
\sum_{n_1,\dots, n_{2k}}
    \frac{1}{n_1^{\tfrac12 +\alpha_1}\cdots n_{2k}^{\tfrac12-\alpha_{2k}}}
\, 
\chi(n_1)\cdots \overline{\chi}(n_{2k})
.
 \end{eqnarray}
According to the recipe, we replace $\varepsilon_\chi^{{k}-\ell}$
by its expected value.  
Since the $\varepsilon_\chi$ are uniformly distributed on the unit circle,
the expected value is 1 if $\ell=k$ and 0 otherwise, so we keep
$\binom{2k}{k}$ terms.

Next we replace the summand by its expected value, which is
\begin{eqnarray}
\delta(n_1,\dots,n_{2k})
&=&
\langle \chi(n_1)\cdots\chi(n_k) 
         \overline{\chi}(n_{k+1})\cdots \overline{\chi}(n_{2k})\rangle \cr
&=&
\begin{cases}
1 & n_1\cdots n_k = n_{k+1} \cdots  n_{2k} \text{\ \ and\ \ } (n_1\cdots n_{2k},q)=1 \\
0 & otherwise .
\end{cases}
\end{eqnarray}

The above is almost identical to the conjectures obtained for the mean values, 
in $t$-aspect, for a primitive $L$-function.  So one obtains the
same formulas as appear in Conjectures~\ref{thm:zeta2kconjecture} 
and~\ref{thm:general2kconjecture}, the only changes
being that one omits the factors $p|q$ in the Euler product~$A_k$,
and one must use the factors $X_\chi(\frac12\pm \alpha_j)^{-\frac12}$.
Specifically, in  Conjectures~\ref{thm:zeta2kconjecture} and~\ref{thm:general2kconjecture}
a simplification occurred by use of equations (\ref{eqn:chis}) and (\ref{eqn:Xs}).
If those conjectures are written in terms of $\prod X(\frac12\pm z_j)^{-\frac12}$, then
the Dirichlet $L$-function moment conjecture would be obtained by substituting
with $\prod X_\chi(\frac12\pm z_j)^{-\frac12}$
Note that we are considering the averages over the even and odd 
primitive characters separately, so  in the sum $X_\chi$ only depends 
on the conductor of~$\chi$.
See the comments following the Theorems in Section~\ref{sec:mainconjectures}
for more discussion on these $X$-factors and conductors.


\subsection{  Symplectic and Orthogonal: quadratic twists of  a
real $L$-function
}\label{sec:symplectic}

Next we consider what happens when we average the shifts of
central values of  $\mathcal{L}(s)$ twisted by the family of quadratic
characters $\displaystyle \chi_d(n)=\left(\frac{d}{n}\right)$,
$d<0$ a fundamental discriminant.  Here $\displaystyle
\chi_d(n)=\left(\frac{d}{n}\right)$ is the Kronecker symbol which
is a primitive Dirichlet character of conductor $|d|$.  
We will see that the family can be either Symplectic or Orthogonal,
depending on the particular $L$-function we start with.

Again, let
\begin{equation}
\mathcal{L}(s)=\sum_{n=1}^\infty \frac{a_n}{n^s}=\prod_p
\mathcal{L}_p(1/p^s)
\end{equation}
be a primitive $L$-function and note that
\begin{equation}
\mathcal{L}(s,\chi_d)= \sum_{n=1}^\infty \frac {a_n
\left(\frac{d}{n}\right)}{n^s}.
\end{equation}

We assume that $\mathcal{L}$ is real,
i.e. $\mathcal{L}=\overline{\mathcal{L}}$,
as this case is relatively easy to deal with from a
fairly general perspective.
Thus,
\begin{equation}
\mathcal{L}(s)= \varepsilon \mathcal{X}(s)\mathcal{L}(1-s),
\end{equation}
where $\varepsilon=\pm1$.
The twisted $L$-function is expected to satisfy a functional equation of
the form
\begin{equation} \label{eqn:Lschidfe}
\mathcal{L}(s,\chi_d)= \varepsilon_d \mathcal{X}_d(s)
				{\mathcal{L}}(1-s,\chi_d).
\end{equation}
It is further expected that 
\begin{equation} \label{eqn:Xdform}
\mathcal{X}_d(s) =|d|^{w(\frac12 -s)} X(s, d), 
\end{equation}
where there are only finitely 
many possibilities for $X(s, d)$. 
By our definition of ``family'' we require that the parameters in the
functional equation be monotonic functions of the conductor.
Since there are only finitely many choices for $X(s, d)$, we must
restrict to averages over sets of $d$ for which $X(s, d)$
is constant (as a function of $d$).
In the situation described here, 
it is believed that there exists an integer~$N$, depending on 
$\mathcal{L}$, such that $\varepsilon_d$ and  $X(s,d)$
only depend on the sign of $d$ and on $(d\mod N)$.  Thus, we will
consider the averages
\begin{equation}
\sumstar_{  \ontop{d<0}{d\equiv a \mod N} } 
L_d(\tfrac12;\alpha_1,\dots,\alpha_k) g(|d|),
\end{equation}
(the following analysis holds also for $d>0$)
where $\sumstar$ denotes a sum over fundamental discriminants~$d$,
and
\begin{equation}
L_d(s;\alpha_1,\dots,\alpha_k)=
 Z_\mathcal{L}(s +\alpha_1,\chi_d)\dots
Z_\mathcal{L}(s +\alpha_k,\chi_d)
.
\end{equation}
Note that $\varepsilon_d=\varepsilon_a$, which may depend on the sign of~$d$.
If $N$ is even we are insisting further that it be divisible by at least $8$.

Following the general discussion in Section~\ref{sec:generalrecipe}, 
equation (\ref{eqn:generalfirstell}) specializes in this case~to
\begin{eqnarray}
&   & 
\varepsilon_f^{\frac{k}{2}-\ell}
\prod_{j=1}^\ell \mathcal{X}_d(s +\alpha_j)^{-\frac12}
\prod_{j=\ell+1}^k \mathcal{X}_d(s -\alpha_j)^{-\frac12}
\sum_{n_1,\dots, n_k}
    \frac{a_{n_1}\cdots a_{n_k}
}{n_1^{s +\alpha_1}\cdots n_k^{s-\alpha_k}}
\, 
\chi_d(n_1)\cdots \chi_d(n_k)
.
 \end{eqnarray}
According to the recipe, we replace $\varepsilon_f^{\frac{k}{2}-\ell}$
by its expected value.  
We have assumed (by our choice of $a\mod N$) that $\varepsilon_d=\varepsilon_a$
for all $d$, so the expected value is $\varepsilon_a^{\frac{k}{2}-\ell}$ and
we will have a contribution from all $2^k$ terms.
(That expression is more transparent is you separately consider the cases
$\varepsilon_a=1$ and $\varepsilon_a=-1$).

The next step in the recipe is to replace the
summand by its expected value.  Since 
$\chi_d(n_1)\cdots \chi_d(n_k) = \chi_d(n_1\cdots n_k)$,
from equation~(\ref{eqn:chidaverage}) we have
the expected value
\begin{equation}
\langle\chi_d(n_1)\cdots \chi_d(n_k)\rangle=
\begin{cases}
    \chi_a(g) \prod_{p|\square}  (1+\frac1p)^{-1} & n_1\cdots n_k=g \square  \\
    0 & \text{otherwise},
\end{cases}
\end{equation}
where $(N,\square)=1$, and all the prime factors of $g$ also being
prime factors of $N$.
So the contribution from the term where we use the first part of the
approximate functional equation for the first $\ell$ factors,
and the second part for the rest, is
\begin{equation}
\varepsilon_d^{\frac{k}{2}-\ell}
\prod_{j=1}^\ell \mathcal{X}_d(s +\alpha_j)^{-\frac12}
\prod_{j=\ell+1}^k \mathcal{X}_d(s -\alpha_j)^{-\frac12}
R_{k,N}(s ;\alpha_1,\dots,\alpha_\ell,-\alpha_{\ell+1},\dots,-\alpha_k)
 \end{equation}
where
\begin{equation}
    R_{k,N}(s; \alpha_1,\dots,\alpha_k) =
    \sum_{g \square}
    \chi_a(g)
    \sum_{n_1\cdots n_k=g \square}
    \frac{
        a_{n_1}\cdots a_{n_k}
    }
    {n_1^{s +\alpha_1}\cdots n_k^{s+\alpha_k}}
    \prod_{p|\square} (1+\tfrac1p)^{-1}
\end{equation}

Adding up all $2^k$ terms we obtain
\begin{equation}
M(s;\alpha_1,\dots,\alpha_k) = \sum_{\epsilon_i= \pm 1}
sign(\{\epsilon_j\})
\prod_{j=1}^k \mathcal{X}_d(\tfrac12 +\epsilon_j\alpha_j)^{-\frac
12}R_{k,N}(s; \epsilon_1\alpha_1, \dots, \epsilon_k\alpha_k),
\end{equation}
where
\begin{equation}\label{eqn:sign}
sign(\{\epsilon_j\})=\begin{cases}
1 & \varepsilon_a=1 \cr
(-1)^{\frac12 \sum  \epsilon_i} & \varepsilon_a=-1 .
\end{cases}
\end{equation}
So the recipe has produced the conjecture
\begin{equation}
\sumstar_{  \ontop{d<0}{d\equiv a \bmod N}}
L_d(\tfrac12 ,\alpha_1,\dots,\alpha_k)
        g(d)
=\sumstar_{ \ontop{d<0}{d\equiv a \bmod N}}
M(\tfrac12;\alpha_1,\dots,\alpha_k)
(1+ O(|d|^{-\frac 12+\varepsilon})) g(d).
\end{equation}

To put the conjecture in a more useful form, we now write $R_{k,N}$ as an 
Euler product, and then express the main term as a contour
integral.

We have $R_{k,N}=\prod_p R_{k,N,p}$,
which naturally separates into a product  over the primes which
divide $N$ and a product over the primes which do not divide~$N$.
The $p$-factor when $p\nmid N$ is 
\begin{eqnarray}
R_{k,N,p}(s)&  = &
 \left(1+ \left(1+\frac 1 p \right)^{-1}\sum_{j=1}^\infty
\sum_{e_1+\dots e_k=2j}\prod_{i=1}^k
\frac{a_{p^{e_i}}}{p^{e_i(s +\alpha_i)}} \right)\nonumber \\
&  = &
\left(1+\frac 1 p\right)^{-1} \left(\frac{1}{p}+
\sum_{j=0}^\infty\sum_{e_1+\dots e_k=2j}
\prod_{i=1}^k \frac{a_{p^{e_i}}}{p^{e_i(s +\alpha_i)}} \right)\nonumber \\
&  = &\left(1+\frac 1 p\right)^{-1}\left(\frac 1 p +\frac{1}{2}
\left(\prod_{j=1}^k
\mathcal{L}_p \left(\frac{1}{p^{s +\alpha_j}} \right)+ 
\prod_{j=1}^k\mathcal{L}_p
\left(\frac{-1}{p^{s +\alpha_j}}\right) \right)\right)
. \end{eqnarray}
Similarly, the $p$-factor when $p| N$ is
\begin{equation}
    R_{k,N,p}=
        \prod_{j=1}^k 
        \mathcal{L}_p \left(\frac{\chi_a(p)}{p^{s +\alpha_j}} \right).
\label{eqn:RNp}
\end{equation}


The above expression will enable us to locate the leading poles of $R_{k,N}$.
Consider the expansion of $R_{k,N,p}$ (for $p\nmid N$) in powers of $1/p$.
The expansion is of the form
\begin{eqnarray}\label{eqn:RNp1overp}
&&1+  \sum_{j=1}^k \frac{a_{p^2}}{p^{2s+2\alpha_j}}+
\sum_{1\le i < j \le k} \frac {(a_p)^2}{p^{2s+\alpha_i+\alpha_j}} + 
O(p^{-1-2s+\varepsilon})+O(p^{-3s+\varepsilon} )\cr
&&\phantom{XXXX}= 
\prod_{j=1}^k
\left(1+\frac{a_{p^2}}{p^{2s+\alpha_j}}\right)
\times
\prod_{1\le i < j \le k}
\left(1+\frac{(a_p)^2}{p^{2s+\alpha_i+\alpha_j}}\right)
\times
\left(1+O(p^{-1-2s+\varepsilon})+O(p^{-3s+\varepsilon} )\right)
\end{eqnarray}
We assume that
\begin{equation}\prod_p\left(1+\frac{(a_p)^2}{p^s}\right)
\label{eqn:squareofap}
\end{equation}
has a simple pole at $s=1$.  This is conjectured to be equivalent to
$\mathcal{L}(s)$ being a primitive $L$-function, and this is the
key place where the assumption of primitivity enters the calculation.
We also assume that
\begin{equation}\prod_p\left(1+\frac{a_{p^2}}{p^s}\right)
\end{equation}
has a pole of order $\delta=0$ or~1 at $s=1$.

In general, $\delta$ is expected to be 0 or 1
according to whether the symmetric square $L$-function
of $\mathcal{L}(s)$ is analytic at $s=1$ 
or has a simple pole at $s=1$. 
If $\mathcal{L}(s)$ is a degree~1 $L$-function, (that is, the 
Riemann $\zeta$-function or a Dirichlet $L$-function), then~$\delta=1$.
If $\mathcal{L}(s)$ is associated to a $GL(2)$ automorphic
form, then $\delta=0$ in general (except possibly when $\mathcal{L}$ is a
dihedral  Artin
$L$-function associated to an weight~1 modular form). 

Note that $\prod_p (1+O(p^{-1-2s})+O(p^{-3s} ))$ is regular 
for $\sigma>\frac13$.  
Thus the total order of the pole of the above product at
$s=\frac12$ when 
$\alpha_1=\dots =\alpha_k=0$
is $\frac12 k(k-1)+\delta k$.
Accordingly, we factor out appropriate zeta-factors
and write the above product as
\begin{eqnarray}
R_{k,N}(s)&=& \prod_{1\le i < j \le k} \zeta (2s+\alpha_i+\alpha_j) \prod_p
R_{k,N,p}(s) \prod_{1\le i< j\le k}
\left(1-\frac{1}{p^{2s+\alpha_i+\alpha_j}}\right) \cr
\end{eqnarray}
if $\delta = 0$, and as
\begin{eqnarray}
R_{k,N}(s)&=& \prod_{1\le i \le j \le k} \zeta (2s+\alpha_i+\alpha_j)
\prod_p R_{k,N,p}(s) \prod_{1\le i \le j\le k}
\left(1-\frac{1}{p^{2s+\alpha_i+\alpha_j}}\right) \cr
\end{eqnarray}
if $\delta=1$.  
In the first case above the family is Orthogonal, and in the second
case it is Symplectic.

In summary, we are led to conjecture
 \begin{eqnarray}\label{thm:Lhalfchidsum}&   &\sumstar_{  \ontop{d<0}{d\equiv a \bmod N}}
Z_\mathcal{L}(\tfrac12 +\alpha_1,\chi_d)\dots
        Z_\mathcal{L}(\tfrac12 +\alpha_k,\chi_d)g(|d|)\nonumber \\
&   &\phantom{xxxx} =
 \sum_{\epsilon_i=\pm 1}
sign(\{\epsilon_i\})
\prod_{j=1}^k X(\tfrac12 +\epsilon_j\alpha_j,a)^{-\frac 12}
\sumstar_{ \ontop{d<0}{d \equiv a \bmod N}}
R_{k,N}(\tfrac12 ,\epsilon_1\alpha_1, \dots, \epsilon_k\alpha_k)
  |d|^{\frac w 2 \sum_{j=1}^k\epsilon_j \alpha_j}  \nonumber \\
&   &\phantom{xxxxxxxxxxxxxxxxxxxxxXXXXXXXXXX}
\times(1+ O(|d|^{-\frac 12+\varepsilon}))\,g(|d|) . \nonumber \\
\end{eqnarray}
The analogous sum over $d>0$ leads to a similar conjecture.
Here $sign(\{\epsilon_i\})$ is given in (\ref{eqn:sign}) and in 
either case we can use Lemma~\ref{thm:concisesumsymplectic} to write the sum as a 
contour integral.

In the case that $\mathcal{L}(s)$ is the Riemann
zeta-function, the above reduces to Conjecture~\ref{thm:Lhalfchidconjecture}.


\subsection{  Orthogonal: $L$-functions associated with cusp
forms }

Recall that the set of primitive newforms $f\in S_n(\Gamma_0(q))$
is denoted by $H_n(q)$. In this section we consider the shifted
moments of $L_f(s)=\sum_{n=1}^\infty \lambda_f(n)n^{-s}$ near the
critical point $s=\frac12 $ averaged over $f\in H_n(q)$. 
Note that in the language of
Section~\ref{sec:familiesofLfunctions} these $L$-functions are the
twists of $\zeta(s)$ by the family of characters $H_k(q)$ and
would be denoted as $\zeta(s,f)$.  However, we use the more common
notation here.

The functional equation is
\begin{equation}
L_f(s)=\varepsilon_f X(s)L_f(1-s),
\end{equation}
where $\varepsilon_f=-\sqrt{q}\lambda_f(q)=\pm 1$.

We consider the ``harmonic average''
\begin{equation}
\sumhar_{f\in H_n(q)} L_f(\tfrac12 +\alpha_1)\dots L_f(\tfrac12
+\alpha_k)
\end{equation}
which attaches a weight $\langle f,f\rangle ^{-1}$ to each
summand. That is,
\begin{equation}
\sumhar_{f\in H_n(q)} * \ \ =\ \sum_{f\in H_n(q)} */\langle
f,f\rangle .
\end{equation}
Following the general discussion in
Section~\ref{sec:generalrecipe}, equation
(\ref{eqn:generalfirstell}) specializes in this case~to
\begin{eqnarray}\label{eqn:cuspsumterm}
&   & \varepsilon_f^{k-\ell} \prod_{j=\ell+1}^k X(s
-\alpha_j)^{-1} \sum_{n_1,\dots, n_k}
    \frac{\lambda_f(n_1)\cdots \lambda_f(n_k)
}{n_1^{s +\alpha_1}\cdots n_k^{s-\alpha_k}} .
 \end{eqnarray}
According to the recipe, we replace $\varepsilon_f^{k-\ell}$ by
its expected value. Since $\varepsilon_f$ is randomly $\pm1$, the
expected value is 0 unless $k-\ell$ is even.  Thus, we will have
$2^{k-1}$ terms in the final answer.

Next we replace $\lambda_f(n_1)\cdots \lambda_f(n_k)$ by its
expected value. This is given in Lemma~\ref{thm:chebyshevweights}.
After factoring into an Euler product and summing the relevant
geometric series we see that (\ref{eqn:cuspsumterm}) contributes
$\prod_{j=1}^k
X(s-\alpha_j)^{-\frac12}R(\alpha_1,\dots,\alpha_\ell,-\alpha_{\ell+1},\dots,-\alpha_{k})$,
where
\begin{eqnarray}
&&R(s,\alpha_1,\dots,\alpha_{k})   \cr &&\phantom{XXXX}=
\prod_{j=1}^k X(s +\alpha_j)^{-\frac12 } \prod_p
\frac{2}{\pi}\int_0^\pi \sin^2\theta\prod_{j=1}^k
\frac{e^{i\theta}\left(1-\frac{e^{i\theta}}{p^{s+\alpha_j}}
\right)^{-1} -e^{-i\theta}
\left(1-\frac{e^{-i\theta}}{p^{s+\alpha_j}}\right)^{-1}}{e^{i\theta}-e^{-i\theta}}\,d\theta
\end{eqnarray}
Here remember that $s$ will eventually be set to $\frac12$ and
$X(s)=X(1-s)^{-1}$.  Adding up all $2^{k-1}$ terms we obtain
\begin{equation}
M(s;\alpha_1,\dots,\alpha_k) =
\prod_{j=1}^kX(s-\alpha_j)^{-\frac12}\sum_{\ontop{\epsilon_i=\pm
1}{\prod_{j=1}^k \epsilon_j=1}} R(s; \epsilon_1\alpha_1, \dots,
\epsilon_k\alpha_k),
\end{equation}
so the recipe has produced the conjecture
\begin{eqnarray}
\sumhar_{f\in H_n(q)} L_f(\tfrac12 +\alpha_1)\dots L_f(\tfrac12
+\alpha_k) &=&\sumhar_{f\in H_n(q)}
M(\tfrac12;\alpha_1,\dots,\alpha_k) (1+ O(nq)^{-\frac
12+\varepsilon}) \cr &=& 
\(1+O(nq)^{-\frac12+\varepsilon}\)
M(\tfrac12;\alpha_1,\dots,\alpha_k)
\end{eqnarray}

Summarizing, we have
\begin{conjecture} \label{thm:Zfsumall}
 With
$A_k(\alpha_1,\dots, \alpha_k)$ as in
Conjecture~\ref{thm:sumLfkconjecture}, we have
\begin{eqnarray}&   &\sumhar_{f\in H_n(q)}L_f(\tfrac12 +\alpha_1)\dots L_f(\tfrac12 +\alpha_k)=\nonumber \\
&   &\phantom{xxxxxxxxxx}
    \prod_{j=1}^k X(\tfrac12-\alpha_j)^{-\frac12}\sum _{ \ontop{\epsilon_j=\pm 1}{\prod_{j=1}^k\epsilon_j=1 }}
\prod_{j=1}^k X(\tfrac12 +\epsilon_j\alpha_j)^{-\frac12 }
\nonumber \\
&   &\phantom{xxxxxxxxxxxxxxxxxxxx}\times\prod_{1\le i<j\le
k}\zeta(1+\epsilon_i\alpha_i +\epsilon_j\alpha_j)
 A_k(\epsilon_1\alpha_1,\dots,\epsilon_k\alpha_k)(1 +O(nq)^{-\frac
12+\varepsilon} ).
\end{eqnarray}
\end{conjecture}

For the purpose of considering averages of even forms or odd forms
separately, we  note that
\begin{eqnarray}
&&\sumhar_{f\in H_n(q)}\varepsilon_f L_f(\tfrac12 +\alpha_1)\dots
L_f(\tfrac12 +\alpha_k) \\
&&\phantom{xxxxxxxxx}=X(\tfrac12+\alpha_1)\sumhar_{f\in H_2(q)}
L_f(\tfrac12 -\alpha_1)L_f(\tfrac12+\alpha_2)\dots L_f(\tfrac12
+\alpha_k)\nonumber\\
&&\phantom{xxxxxxxxxx}=X(\tfrac12+\alpha_1)M(\tfrac12;-\alpha_1,\alpha_2,\dots,\alpha_k)
(1+O(nq)^{-\frac12+\varepsilon}).\nonumber
\end{eqnarray}
By looking at the combinations \begin{equation} \sumhar_{f\in
H_n(q)} L_f(\tfrac12 +\alpha_1)\dots L_f(\tfrac12 +\alpha_k)\pm
\sumhar_{f\in H_n(q)}\varepsilon_f L_f(\tfrac12 +\alpha_1)\dots
L_f(\tfrac12 +\alpha_k) \end{equation} this leads to

 \begin{conjecture} \label{thm:Zfsumevenodd}
With $A_k(\alpha_1,\dots, \alpha_k)$ as in
Conjecture~\ref{thm:sumLfkconjecture}, we have
\begin{eqnarray}\label{eqn:evensum}
 & & \sumhar _{ \ontop{f\in H_n(q)}{f \text{ even} }}L_f(\tfrac12 +\alpha_1)\dots L_f(\tfrac12 +\alpha_k)=\nonumber \\
& &\phantom{xxxxxxxx}\frac{1}{2}\prod_{j=1}^k X(\tfrac12
-\alpha_j)^{-\frac12}\sum_{\epsilon_j=\pm 1 } \prod_{j=1}^k
X(\tfrac12 +\epsilon_j\alpha_j)^{-\frac12 } \nonumber \\
&   &\phantom{xxxxxxxxxxxxxxxxxxx}\times\prod_{1\le i<j\le k}
\zeta(1+\epsilon_i\alpha_i +\epsilon_j\alpha_j)
 A_k(\epsilon_1\alpha_1,\dots,\epsilon_k\alpha_k)(1 +O(nq)^{-\frac
12+\varepsilon}),
\end{eqnarray}
and
\begin{eqnarray} & &  \sumhar _{ \ontop{f\in H_n(q)}{f \text{ odd} }}L_f(\tfrac12 +\alpha_1)\dots L_f(\tfrac12 +\alpha_k)=\nonumber \\
&
&\phantom{xxxxxxxx}\frac{1}{2}\prod_{j=1}^kX(\tfrac12-\alpha_j)^{-\frac12}\sum_{\epsilon_j=\pm
1 } \prod_{j=1}^k\epsilon_j X(\tfrac12
+\epsilon_j\alpha_j)^{-\frac12
} \nonumber \\
&   &\phantom{xxxxxxxxxxxxxxxxxxx}\times\prod_{1\le i<j\le k}
\zeta(1+\epsilon_i\alpha_i
+\epsilon_j\alpha_j)A_k(\epsilon_1\alpha_1,\dots,\epsilon_k\alpha_k)
 (1 +O(nq)^{-\frac 12+\varepsilon}) .
\end{eqnarray}
\end{conjecture}
The above formulae can be written as contour integrals using
Lemma~\ref{thm:concisesumsymplectic}, giving expressions analogous
to those in Conjecture~\ref{thm:sumLfkconjecture}.  In particular,
expressing (\ref{eqn:evensum}) as a contour integral and then
letting $\alpha_j\to 0$ gives
Conjecture~\ref{thm:sumLfkconjecture}. 


\section{  Numerical calculations
}\label{sec:numerics}

We compare our conjectures with some numerical calculations. The
agreement is very good. These calculations involve numerically
approximating the coefficients in the conjectured formulae, and
numerically evaluating the mean value. Both of those calculations
are interesting and we will give more details and examples in a
subsequent paper.

\subsection{  Unitary: Riemann zeta-function
}\label{sec:riemannzetacalcs}

The coefficients of $P_2(x)$ in Conjecture~\ref{thm:zeta2kconjecture} 
can be written explicitly in terms of
known constants. When $k=2$ the function
$G(\alpha_1,\alpha_2,\alpha_3,\alpha_4)$ that appears in
Conjecture~\ref{thm:zeta2kconjecture} equals
\begin{equation}
\zeta(2+\alpha_1+\alpha_2- \alpha_3-\alpha_4)^{-1} \prod_{i,j=1}^2
\zeta(1+\alpha_i-\alpha_{k+j}),
\end{equation}
which is given in Section~\ref{sec:zetameanvalues}.

But
\begin{equation}
\zeta(1+s)= {s}^{-1}+\gamma-\gamma_1 s +
\frac{\gamma_2}{2!}{s}^{2}- \frac{\gamma_3}{3!}{s}^{3} +\ldots
\end{equation}
and
\begin{equation}
\zeta(2+s)^{-1} = \frac{6}{\pi^2}-\frac{36 \zeta'(2)}{\pi^4} s +
\frac{-3\pi^2\zeta''(2) +36 \zeta'(2)^2 }{\pi^6} s^2 + \ldots.
\end{equation}
In the latter, the terms up to $s^4$ were evaluated using Maple.
Computing the residue in Conjecture~\ref{thm:zeta2kconjecture} gives

\begin{eqnarray}
P_2(x) & =  &\frac{1}{2\pi^2} x^4 +
\frac{8}{\pi^4}
    \left(\gamma{\pi }^{2}-3\zeta'(2)\right) x^3  \nonumber \\&   &\mbox{} +\frac{6}{\pi^6}
    \left(
    -48\gamma\zeta'(2){\pi }^{2}-12\zeta''(2){\pi }^{2}+7\gamma^{2}{\pi }^{4}+
    144\zeta'(2)^{2}-2\,\gamma_1{\pi }^{4}\right) x^2  \nonumber \\&   &\mbox{} +\frac{12}{\pi^8}
    \biggl(
    6\gamma^{3}\pi^{6}-84\gamma^{2}\zeta'(2)\pi^{4}+24\gamma_1\zeta'(2)
    \pi^{4}-1728\zeta'(2)^{3}+
    576\gamma\zeta'(2)^{2}\pi^{2} \nonumber \\&   &\phantom{TTTTT}+288\zeta'(2)\zeta''(2)\pi^{2}-
    8\zeta'''(2)\pi^{4}-10\gamma_1\gamma\pi^{6}-
    \gamma_2\pi^{6}-48\gamma\zeta''(2)\pi^{4}
    \biggr)
    x  \nonumber \\&   &\mbox{} +\frac{4}{\pi^{10}}
    \biggl(
    -12\zeta''''(2){\pi }^{6}+36\gamma_2\zeta'(2){\pi }^{6}+9{\gamma}^{4}{\pi }^{8}+
    21\gamma_1^{2}{\pi }^{8}+432\zeta''(2)^{2}{\pi }^{4} \nonumber \\&   &\phantom{TTTTT}+3456
    \gamma\zeta'(2)\zeta''(2){\pi }^{4}+3024{\gamma}^{2}\zeta'(2)^{2}{\pi }^{4}-
    36{\gamma}^{2}\gamma_1{\pi }^{8}-252{\gamma}^{2}\zeta''(2){\pi }^{6}\nonumber \\&   &\phantom{TTTTT}+3\gamma\gamma_2{\pi}^{8}+
    72\gamma_1\zeta''(2){\pi }^{6}+360\gamma_1\gamma\zeta'(2){\pi }^{6}-216{\gamma}^{
    3}\zeta'(2){\pi }^{6}\nonumber \\&   &\phantom{TTTTT}-864\gamma_1\zeta'(2)^{2}{\pi }^{4}+5\gamma_3{\pi }^{8} +
    576\zeta'(2)\zeta'''(2){\pi }^{4}-20736\gamma\zeta'(2)^{3}{\pi }^{2} \nonumber \\&   &\phantom{TTTTT}-
    15552\zeta''(2)\zeta'(2)^{2}{\pi }^{2}-96\gamma\zeta'''(2){\pi }^{6}+62208
    \zeta'(2)^{4}
    \biggr),
\end{eqnarray}
in agreement with a result implied in the work of
Heath-Brown\cite{H-B1} (see \cite{C} where, using \cite{H-B1}, the
same polynomial is worked out, although there are some slight
errors). Numerically,
\begin{eqnarray}
 P_2(x) & =  &0.0506605918211688857219397316048638\,x^4 \nonumber \\&   &\mbox{} + 0.69886988487897996984709628427658502\,x^3  \nonumber \\&   &\mbox{} + 2.425962198846682004756575310160663\,x^2 \nonumber \\&   &\mbox{} + 3.227907964901254764380689851274668\,x\, \nonumber \\&   &\mbox{} + 1.312424385961669226168440066229978
\end{eqnarray}

There are several ways that one can numerically compute the
coefficients of $P_3(x)$, and these will be described in a future paper.
We found
\begin{eqnarray}
P_3(x)& =  &0.000005708527034652788398376841445252313\,x^9 \nonumber \\
&   &\mbox{} + 0.00040502133088411440331215332025984\,x^8 \nonumber \\
&   &\mbox{} +  0.011072455215246998350410400826667\,x^7 \nonumber \\
&   &\mbox{} + 0.14840073080150272680851401518774\,x^6 \nonumber \\
&   &\mbox{} +  1.0459251779054883439385323798059\,x^5 \nonumber \\
&   &\mbox{} + 3.984385094823534724747964073429\,x^4 \nonumber \\
&   &\mbox{} +  8.60731914578120675614834763629\,x^3 + 10.274330830703446134183009522\,x^2 \nonumber \\
&   &\mbox{} +  6.59391302064975810465713392\,x + 0.9165155076378930590178543.
 \end{eqnarray}
One notices that the leading coefficient is much smaller than the
lower order coefficients, which means that, in numerical
calculations, the lower order terms will contribute significantly.
One might suppose that the coefficients of $P_k(x)$ are always
positive. Unfortunately, while this is true for $P_1,\ldots,P_4$,
by $k=5$, negative coefficients begin to appear (see Table~\ref{tab:cr}).

\begin{table}[h!tb]     
\centerline{\scriptsize
\begin{tabular}{|c|c|c|c|c|}
\hline
$r$&    $c_r(4)$ &$c_r(5)$ &$c_r(6)$ &$c_r(7)$ \cr \hline
0&      .24650183919342276e-12&
.141600102062273e-23&.512947340914913e-39&.658228478760010e-59\cr
1&      .54501405731171861e-10&
.738041275649445e-21&.530673280992642e-36&.120414305554514e-55\cr
2&      .52877296347912035e-8&
.177977962351965e-18&.260792077114835e-33&.106213557174925e-52\cr
3&      .29641143179993979e-6&
.263588660966072e-16&.810161321577902e-31&.601726537601586e-50\cr
4&      .1064595006812847e-4&
.268405453499975e-14&.178612973800931e-28&.246062876732400e-47\cr
5&      .25702983342426343e-3&
.199364130924990e-12&.297431671086361e-26&.773901216652114e-45\cr
6&      .42639216163116947e-2&
.111848551249336e-10&.388770829115587e-24&.194786494949524e-42\cr
7&      .48941424514215989e-1&
.484279755304480e-9&.409224261406863e-22&.403076849263637e-40\cr
8&      .38785267&.16398013e-7&.35314664e-20&.69917763e-38\cr 9&
2.1091338&.43749351e-6&.25306377e-18&.1031402e-35\cr 10&
7.8325356& .92263335e-5 &.15198191e-16  &.13082869e-33\cr 11&
19.828068& .00015376778 &.77001514e-15  &.14392681e-31\cr 12&
33.888932& .0020190278  &.3306121e-13   &.13825312e-29\cr 13&
38.203306& .020772707   &.12064042e-11  &.11657759e-27\cr 14&
25.604415& .16625059    &.37467193e-10  &.86652477e-26\cr 15&
10.618974&1.0264668     &.99056943e-9   &.56962227e-24\cr 16&
.708941&4.8485893       &.22273886e-7   &.33197649e-22\cr 17&
&17.390876& .42513729e-6        &.1718397e-20\cr 18& &47.040877&
.68674336e-5        &.79096789e-19\cr 19& &95.116618& .9351583e-4
&.32396929e-17\cr 20&     &141.44446& .0010683164
&.11809579e-15\cr 21&     &149.35697& .010180702 &.3830227e-14\cr
22&     &105.88716& .080418679  &.11044706e-12\cr 23& &44.1356&
.52296142     &.28282258e-11\cr 24& &20.108&2.7802018
&.64210662e-10\cr 25&     &-1.27&12.001114 &.12898756e-8\cr 26&
&&41.796708     &.22869667e-7\cr 27& &&116.72309 &.35683995e-6\cr
28&     &&259.39898 &.48834071e-5\cr 29& &&452.491&
.58391045e-4\cr 30& &&601.17& .00060742037\cr 31& &&573.54&
.0054716438\cr 32& &&374.8& .042465904\cr 33& &&246.5&
.28245494\cr 34& &&248.&1.6013331\cr 35&     &&1.6e+02 ?&7.6966995\cr
36&     &&-4.e+01 ?&31.20352\cr 37& &&&106.19714\cr 38&
&&&301.91363\cr 39& &&&711.742\cr 40& &&&1370.10\cr 41&
&&&2083.\cr 42& &&&2356.\cr 43& &&&1.9e+03 \cr 44&     &&&1.8e+03
\cr 45& &&&3.e+03  \cr 46& &&&3.e+03 \cr 47&     &&&8.e+01 ?\cr 48&
&&&-1.e+03 ?\cr 49& &&&-2.e+02 ? \cr \hline
\end{tabular}
}
\caption{
Coefficients of $P_k(x) = c_0(k) x^{k^2} + c_1(k) x^{k^2-1} +
\cdots + c_{k^2}(k)$, for $k=4,5,6,7$. Notice the relatively small
size of $c_0(k)$. We believe the coefficients to be correct to the
number of places listed, except in the cases indicated by question
marks, where the numerics have not quite stabilized. Two different
methods were used to compute the coefficients. The former, for $0
\leq r \leq 7$, gave us higher precision but was less efficient,
while the latter for $r \leq 49$, was more efficient but required
using less precision.
}\label{tab:cr}\end{table}

Table~\ref{tab:sixthmoment} depicts
\begin{equation}
\int_C^D |\zeta(\tfrac12 +it)|^6 dt \label{eqn:zeta6smooth}
\end{equation}
as compared to
\begin{equation}
\int_C^D P_3(\log(t/2\pi)) dt,
 \label{eqn:p3smooth}
\end{equation}
along with their ratio, for various blocks $[C,D]$ of length
50000, as well as a larger block of length 2,350,000. The data
agree with our conjecture and are consistent with a remainder of
size $|D-C|^{\frac12 }D^\varepsilon$.

\begin{table}[h!tb]
\centerline{
\begin{tabular}{|c|c|c|c|}
\hline
$[C,D]$ & conjecture (\ref{eqn:p3smooth}) & reality (\ref{eqn:zeta6smooth})& ratio \cr \hline
[0,50000] & 7236872972.7 & 7231005642.3    & .999189\cr
[50000,100000] & 15696470555.3 & 15723919113.6 & 1.001749\cr
[100000,150000]& 21568672884.1&21536840937.9   &     .998524\cr
[150000,200000]& 26381397608.2&26246250354.1   &     .994877\cr
[200000,250000]& 30556177136.5&30692229217.8   &    1.004453\cr
[250000,300000]& 34290291841.0&34414329738.9   &    1.003617\cr
[300000,350000]& 37695829854.3&37683495193.0   &     .999673\cr
[350000,400000]& 40843941365.7&40566252008.5   &     .993201\cr
[400000,450000]& 43783216365.2&43907511751.1   &    1.002839\cr
[450000,500000]& 46548617846.7&46531247056.9   &     .999627\cr
[500000,550000]& 49166313161.9&49136264678.2  &      .999389\cr
[550000,600000]& 51656498739.2&51744796875.0  &     1.001709\cr
[600000,650000]& 54035153255.1&53962410634.2  &      .998654\cr
[650000,700000]& 56315178564.8&56541799179.3  &     1.004024\cr
[700000,750000]& 58507171421.6&58365383245.2  &      .997577\cr
[750000,800000]& 60619962488.2&60870809317.1  &     1.004138\cr
[800000,850000]& 62661003164.6&62765220708.6  &     1.001663\cr
[850000,900000]& 64636649728.0&64227164326.1  &      .993665\cr
[900000,950000]& 66552376294.2&65994874052.2  &      .991623\cr
[950000,1000000]& 68412937271.4&68961125079.8  &     1.008013\cr
[1000000,1050000]& 70222493232.7&70233393177.0  &
1.000155\cr [1050000,1100000]& 71984709805.4&72919426905.7  &
1.012985\cr [1100000,1150000]& 73702836332.4&72567024812.4  &
.984589\cr [1150000,1200000]& 75379769148.4&76267763314.7  &
1.011780\cr [1200000,1250000]& 77018102997.5&76750297112.6  &
.996523\cr [1250000,1300000]& 78620173202.6&78315210623.9  &
.996121\cr [1300000,1350000]& 80188090542.5&80320710380.9  &
1.001654\cr [1350000,1400000]& 81723770322.2&80767881132.6  &
.988303\cr [1400000,1450000]& 83228956776.3&83782957374.3  &
1.006656\cr [0,2350000]& 3317437762612.4&3317496016044.9  &
1.000017 \cr \hline
\end{tabular}
}
\caption{ Sixth moment of $\zeta$
versus Conjecture~\ref{thm:zeta2kconjecture}. The `reality' column, i.e. integrals
involving $\zeta$, were computed using Mathematica.
}\label{tab:sixthmoment}\end{table}

One can also look at smoothed moments, for example,
\begin{equation}
\int_0^\infty |\zeta(\tfrac12 +it)|^{2k} \exp(-t/T) dt \label{eqn:zeta2ksmooth}
\end{equation}
as compared to
\begin{equation}
\int_0^\infty P_k(\log(t/2\pi)) \exp(-t/T) dt. \label{eqn:p2ksmooth}
\end{equation}

Table~\ref{tab:zeta2ksmooth} compares these with $T=10000$, for $k=4,3,2,1$.

\begin{table}[h!tb]
\centerline{
\begin{tabular}{|c|c|c|c|c|}
\hline
$k$ & (\ref{eqn:zeta2ksmooth}) & (\ref{eqn:p2ksmooth})& difference&
\ \ \  relative  \cr &&&&\ \ \ difference \cr \hline
 1 & 79499.9312635 &79496.7897047 & 3.14156&\ \ $3.952\times
10^{-5}$\cr 2 & 5088332.55512 &5088336.43654 &
-3.8814&$-7.628\times 10^{-7}$\cr 3 & 708967359.4 &708965694.5 &
1664.9&\ \ $2.348\times 10^{-6}$\cr 4 & 143638308513.0 &
143628911646.6 & 9396866.4 &\ \ $6.542\times 10^{-5} $\cr \hline
\end{tabular}
}
\caption{ Smoothed moment of $\zeta$ versus
Conjecture~\ref{thm:zeta2kconjecture}.
}\label{tab:zeta2ksmooth}\end{table}

For $k=3,4$ the data agrees to roughly half the decimal places.
This supports our conjecture that the error term in the
conjectured mean values is $O(T^{\frac{1}{2}+\varepsilon})$.  For
$k=1$ the numerics suggest Corollary~\ref{thm:pimoment}.


\subsection{  Symplectic: quadratic Dirichlet $L$-functions
}

We have computed the polynomials $Q_k$ of Conjecture~\ref{thm:Lhalfchidconjecture} for
$k=1,2,\ldots,8$, separately for $d<0$ and $d>0$. Table~\ref{tab:Qkcoeffn}
lists these polynomials for $d<0$, while in Table~\ref{tab:Qkcoeffp} we
consider $d>0$. Again notice the small size of the
leading coefficients.

\begin{table}[h!tb]
\centerline{\small
\begin{tabular}{|c|c|c|c|c|}
\hline
$r$&    $d_r(1)$ &$d_r(2)$ &$d_r(3)$ &$d_r(4)$ \cr \hline
0 & .3522211004995828 & .1238375103096e-1 & .1528376099282e-4 & .31582683324433e-9 \cr
1 & .61755003361406 & .18074683511868 & .89682763979959e-3 & .50622013406082e-7 \cr
2 & & .3658991414081 & .17014201759477e-1 & .32520704779144e-5 \cr
3 & & -.13989539029 & .10932818306819 & .10650782552992e-3 \cr
4 & & & .13585569409025 & .18657913487212e-2 \cr
5 & & & -.23295091113684 & .16586741288851e-1 \cr
6 & & & .47353038377966 & .59859999105052e-1 \cr
7 & & & & .52311798496e-2 \cr
8 & & & & -.1097356195 \cr
9 & & & & .55812532 \cr
10 & & & & .19185945 \cr
\hline\hline
$r$&    $d_r(5)$ &$d_r(6)$ &$d_r(7)$ &$d_r(8)$ \cr \hline
0 & .671251761107e-16 & .1036004645427e-24 & .886492719e-36 & .337201e-49 \cr
1 & .23412332535824e-13 & .67968140667178e-22 & .98944375081241e-33 & .59511917e-46 \cr
2 & .35711692341033e-11 & .20378083365099e-19 & .51762930260135e-30 & .500204322e-43 \cr
3 & .31271184907852e-9 & .36980514080794e-17 & .16867245856115e-27 & .2664702284e-40 \cr
4 & .17346173129392e-7 & .45348387982697e-15 & .38372675160809e-25 & .1010164552e-37 \cr
5 & .63429411057027e-6 & .39728668850800e-13 & .64746354773372e-23 & .29004988867e-35 \cr
6 & .15410644373832e-4 & .2563279107877e-11 & .84021141030379e-21 & .65555882460e-33 \cr
7 & .2441498848698e-3 & .12372292296e-9 & .85817644593981e-19 & .11966099802e-30 \cr
8 & .2390928284571e-2 & .44915158297e-8 & .70024645896e-17 & .17958286298e-28 \cr
9 & .127561073626e-1 & .1222154548e-6 & .4607034349989e-15 & .22443685425e-26 \cr
10 & .24303820161e-1 & .2461203700e-5 & .2455973970377e-13 & .2357312577e-24 \cr
11 & -.333141763e-1 & .3579140509e-4 & .106223013225e-11 & .20942850060e-22 \cr
12 & .25775611e-1 & .3597968761e-3 & .3719625461492e-10 & .15805997923e-20 \cr
13 & .531596583 & .230207769e-2 & .1048661496741e-8 & .10159435845e-18 \cr
14 & -.325832 & .7699469185e-2 & .2357398870407e-7 & .55665248752e-17 \cr
15 & -1.34187 & .4281359929e-2 & .416315210727e-6 & .25985097519e-15 \cr
16 & & -.2312387714e-1 & .564739434674e-5 & .103134457e-13 \cr
17 & & .109503 & .56831273239e-4 & .346778002e-12 \cr
18 & & .2900464 & .40016131254e-3 & .982481680e-11 \cr
19 & & -.9016 & .1755324808e-2 & .232784142e-9 \cr
20 & & -.89361 & .340409901e-2 & .456549799e-8 \cr
21 & & -.181 & -.2741804e-2 & .7309216472e-7 \cr
22 & & & .353555e-3 & .9368893764e-6 \cr
23 & & & .117734 & .9348804928e-5 \cr
24 & & & .20714e-1 & .69517414e-4 \cr
25 & & & -.9671 & .356576507e-3 \cr
26 & & & -.284 & .1059852e-2 \cr
27 & & & 1.3 & .8242527e-3 \cr
28 & & & -1. & -.206921e-2 \cr
29 & & & & .181031e-1 \cr
30 & & & & .862815e-1 \cr
31 & & & & -.14025 \cr
32 & & & & -.91619 \cr
33 & & & & -.942 \cr
34 & & & & -.153e-1 \cr
35 & & & & -.3 ?\cr
36 & & & & ? \cr
\hline
\end{tabular}
}
\caption{ Coefficients of
$Q_k(x) = d_0(k) x^{k(k+1)/2} + d_1(k) x^{k(k+1)/2} + \cdots $, for $k=1,\ldots,8$,
odd twists, $d<0$.
}\label{tab:Qkcoeffn}\end{table}

\begin{table}[h!tb]
\centerline{\small
\begin{tabular}{|c|c|c|c|c|}
\hline
$r$&    $e_r(1)$ &$e_r(2)$ &$e_r(3)$ &$e_r(4)$ \cr \hline
0 & .3522211004995828 & .1238375103096e-1 & .1528376099282e-4 & .31582683324433e-9 \cr
1 & -.4889851881547 & .6403273133043e-1 & .60873553227400e-3 & .40700020814812e-7 \cr
2 & & -.403098546303 & .51895362572218e-2 & .19610356347280e-5 \cr
3 & & .878472325297 & -.20704166961612e-1 & .4187933734219e-4 \cr
4 & & & -.4836560144296e-1 & .32338329823195e-3 \cr
5 & & & .6305676273171 & -.7264209058150e-3 \cr
6 & & & -1.23114954368 & -.97413031149e-2 \cr
7 & & & & .6254058547e-1 \cr
8 & & & & .533803934e-1 \cr
9 & & & & -1.125788 \cr
10 & & & & 2.125417 \cr
\hline\hline
$r$&    $e_r(5)$ &$e_r(6)$ &$e_r(7)$ &$e_r(8)$ \cr \hline
0 & .671251761107e-16 & .1036004645427e-24 & .886492719e-36 & .337201e-49 \cr
1 & .2024913313373e-13 & .6113326104277e-22 & .91146378e-33 & .556982629e-46 \cr
2 & .261100345555e-11 & .16322243213252e-19 & .437008961e-30 & .43686422e-43 \cr
3 & .187088892376e-9 & .2605311255687e-17 & .1297363095e-27 & .216465856e-40 \cr
4 & .8086250862418e-8 & .2766415183453e-15 & .2670392090e-25 & .7604817313e-38 \cr
5 & .2126496335545e-6 & .2056437432502e-13 & .404346681e-23 & .201532781e-35 \cr
6 & .319415704903e-5 & .10957094998959e-11 & .46631481394e-21 & .418459324e-33 \cr
7 & .21201987479e-4 & .42061728711797e-10 & .41831543311e-19 & .698046515e-31 \cr
8 & -.33900555230e-4 & .11491097182922e-8 & .29548572643e-17 & .951665168e-29 \cr
9 & -.775061385e-3 & .21545094604323e-7 & .1652770327e-15 & .1073015400e-26 \cr
10 & .333997849e-2 & .25433712247032e-6 & .73192383650e-14 & .1008662234e-24 \cr
11 & .22204682e-1 & .1448397731463e-5 & .25506469557e-12 & .7945270901e-23 \cr
12 & -.1538433 & -.2179868777201e-5 & .6901276286e-11 & .5257922143e-21 \cr
13 & -.19794e-1 & -.54298634893e-4 & .141485467e-9 & .2924082555e-19 \cr
14 & 2.01541 & .1698771341e-3 & .210241720e-8 & .1363867915e-17 \cr
15 & -4.451 & .22887524e-2 & .20651382e-7 & .5311448709e-16 \cr
16 & & -.1042e-1 & .101650951e-6 & .1714154659e-14 \cr
17 & & -.4339429e-1 & -.16979129e-6 & .453180963e-13 \cr
18 & & .343054 & -.37367e-5 & .9644403068e-12 \cr
19 & & -.1947171 & .97069e-5 & .160742335e-10 \cr
20 & & -3.16910 & .18351e-3 & .200188929e-9 \cr
21 & & 7.31266 & -.54878e-3 & .16931900e-8 \cr
22 & & & -.5621e-2 & .7257434e-8 \cr
23 & & & .284e-1 & -.14329111e-7 \cr
24 & & & .639e-1 & -.25913136e-6 \cr
25 & & & -.7 & .6473933e-6 \cr
26 & & & .86 & .138673e-4 \cr
27 & & & 5. & -.2339e-4 \cr
28 & & & -.1e2 & -.48124e-3 \cr
29 & & & & .162e-2 \cr
30 & & & & .976e-2 \cr
31 & & & & -.83e-1 \cr
32 & & & & -.62e-1 \cr
33 & & & & 2. \cr
34 & & & & -2. \cr
35 & & & & -9. \cr
36 & & & & 30.? \cr \hline
\end{tabular}
}
\caption{ Coefficients of
$Q_k(x) = e_0(k) x^{k(k+1)/2} + e_1(k) x^{k(k+1)/2} + \cdots $, for $k=1,\ldots,8$,
even twists, $d>0$.
}\label{tab:Qkcoeffp}\end{table}

Table~\ref{tab:Lhalfchin} compares, for $d<0$,
conjectured moments for $k=1,\ldots,8$ against numerically computed moments,
\begin{equation}
    \sumstar_{d<0} L(\tfrac12 ,\chi_d)^k g(|d|)
\label{eqn:lhchidn}
\end{equation}
versus
\begin{equation}
   \sumstar_{d<0}\,Q_k(\log{|d|}) g(|d|)
\label{eqn:Qkn}
\end{equation}
where $g$ is the smooth test function
\begin{equation}
g(t) =
   \begin{cases}
        1,   &\text{if $0\leq t < 850000$;} \\
        \exp\left(1-\frac{1}{1-\frac{(t-850000)^2}{(150000)^2}}\right),
             &\text{if $850000 \leq t \leq 1000000$;} \\
        0,   &\text{if $1000000 < t $.}
    \end{cases}
\end{equation}
Table~\ref{tab:Lhalfchip} compares the same quantities, but for $d>0$.

\begin{table}[h!tb]
\centerline{\small
\begin{tabular}{|c|c|c|c|}
\hline
$k$ & reality (\ref{eqn:lhchidn}) & conjecture (\ref{eqn:Qkn})& ratio \cr \hline
1 &1460861.8       &1460891.          &0.99998 \cr 2 &17225813.8
&17226897.5       &0.999937 \cr 3 &316065502.1      &316107868.6
&0.999866 \cr 4 &7378585496.      &7380357447.1     &0.99976  \cr
5 &198754711593.6   &198809762196.4   &0.999723 \cr 6
&5876732216291.7 &5877354317291.3   &0.999894 \cr 7
&185524225881950. &185451557119001. &1.000392 \cr 8
&6149876164696600 &6141908614344770 &1.0013 \cr \hline
\end{tabular}
}
\caption{ Smoothed moment of $L(\frac12 ,\chi_d)$ versus
Conjecture~\ref{thm:Lhalfchidconjecture}, for fundamental discriminants $-1000000 < d <
0$, and $k=1,\ldots,8$.
}\label{tab:Lhalfchin}\end{table}


Figure \ref{fig:conj vs reality, odd chi} depicts, for $k=1,\ldots,8$
and $X=10000,20000,\ldots,10^7$,
\begin{equation}
    \sumstar_{-X<d<0} L(\frac12 ,\chi_d)^k 
\label{eqn:chi odd sharp cutoff}
\end{equation}
divided by
\begin{equation}
   \sumstar_{-X<d<0}\,Q_k(\log{|d|}).
\label{eqn:conj chi odd sharp cutoff}
\end{equation}
One sees the graphs fluctuating above and below one. Interestingly, the
graphs have a similar shape as $k$ varies. This is explained by the fact that
large values of $L(1/2,\chi_d)$ tend to skew the moments, and this gets
amplified as $k$ increases.

Figure \ref{fig:conj vs reality, even chi} depicts the same but
for $0 < d \leq X$.

The values $L(\frac12 ,\chi_d)$ were computed using a smoothed form of the
approximate functional equation which expresses the $L$-function in terms of the
incomplete Gamma-function (see for example~\cite{L}).

\begin{figure}[h!tb]
    \centerline{
            \psfig{figure=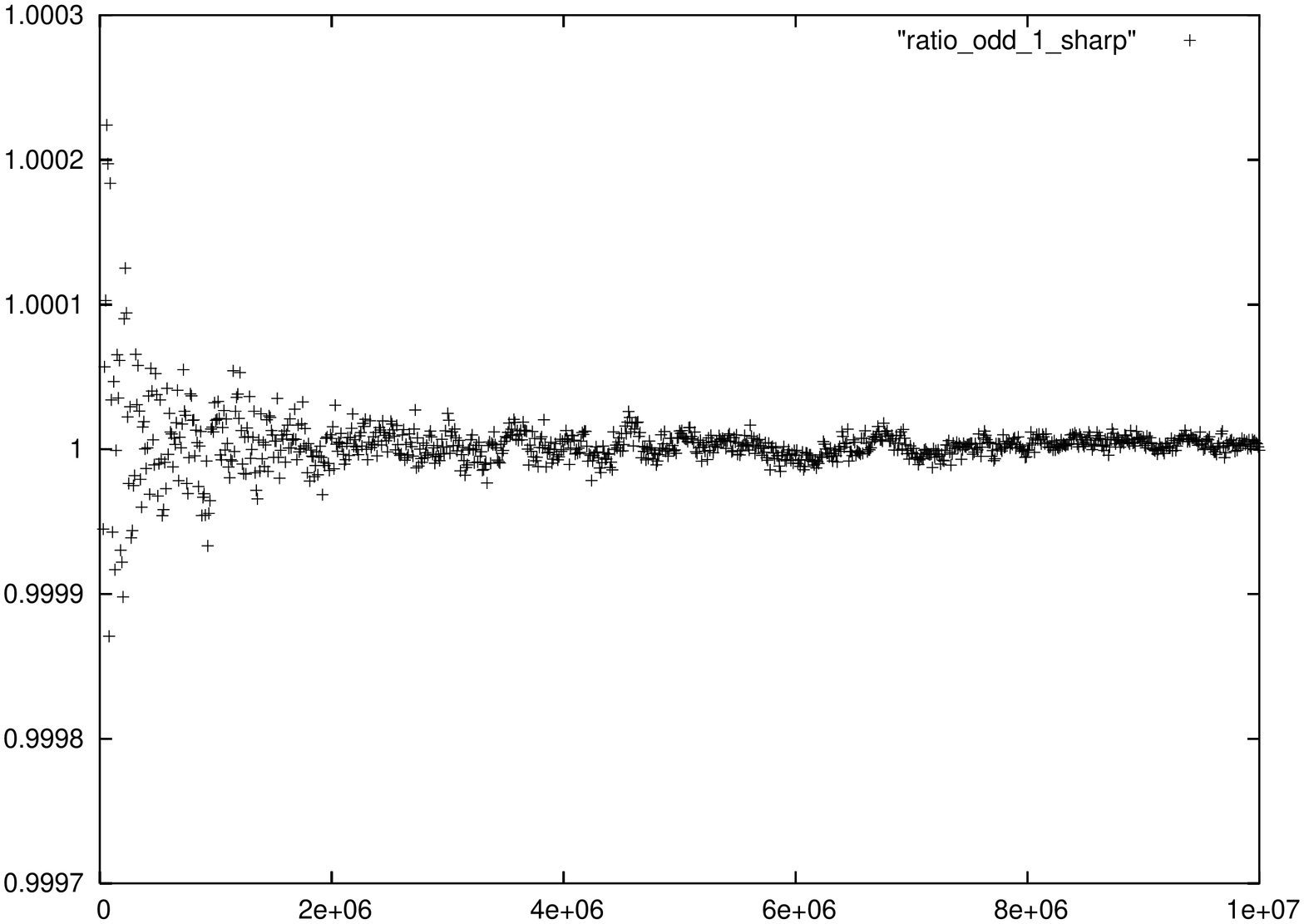,width=3in,angle=-90}
            \psfig{figure=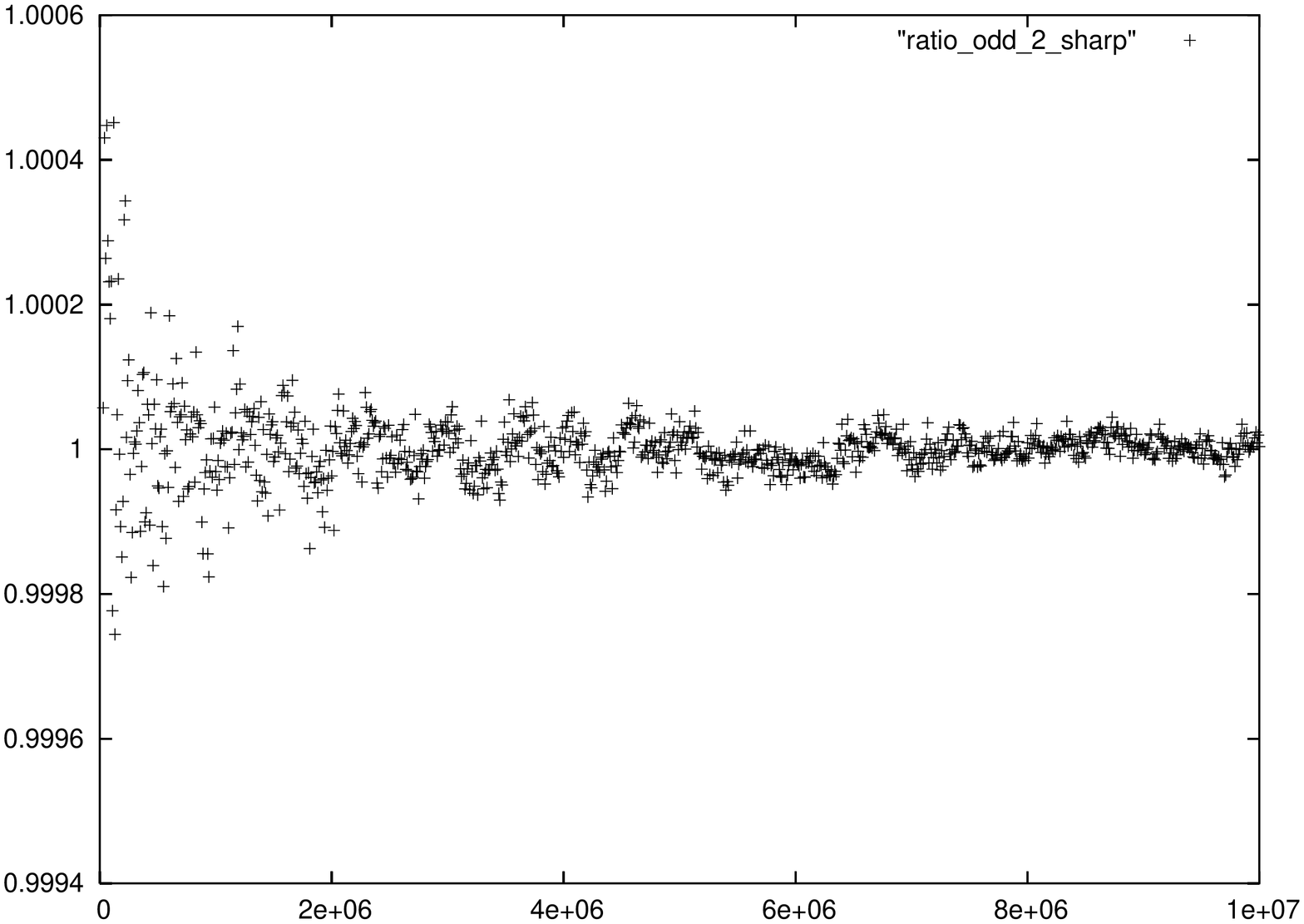,width=3in,angle=-90}
    }
    \centerline{
            \psfig{figure=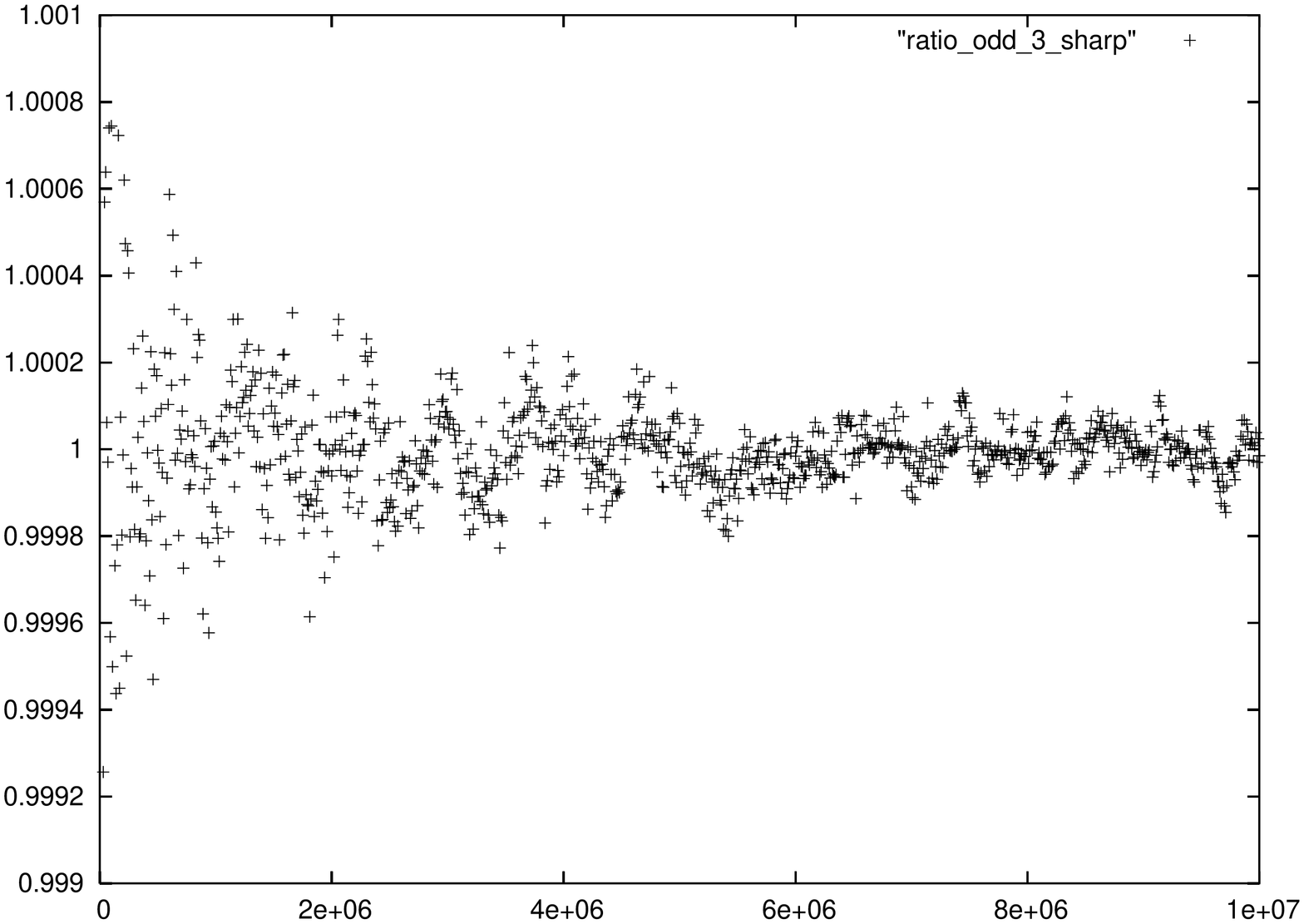,width=3in,angle=-90}
            \psfig{figure=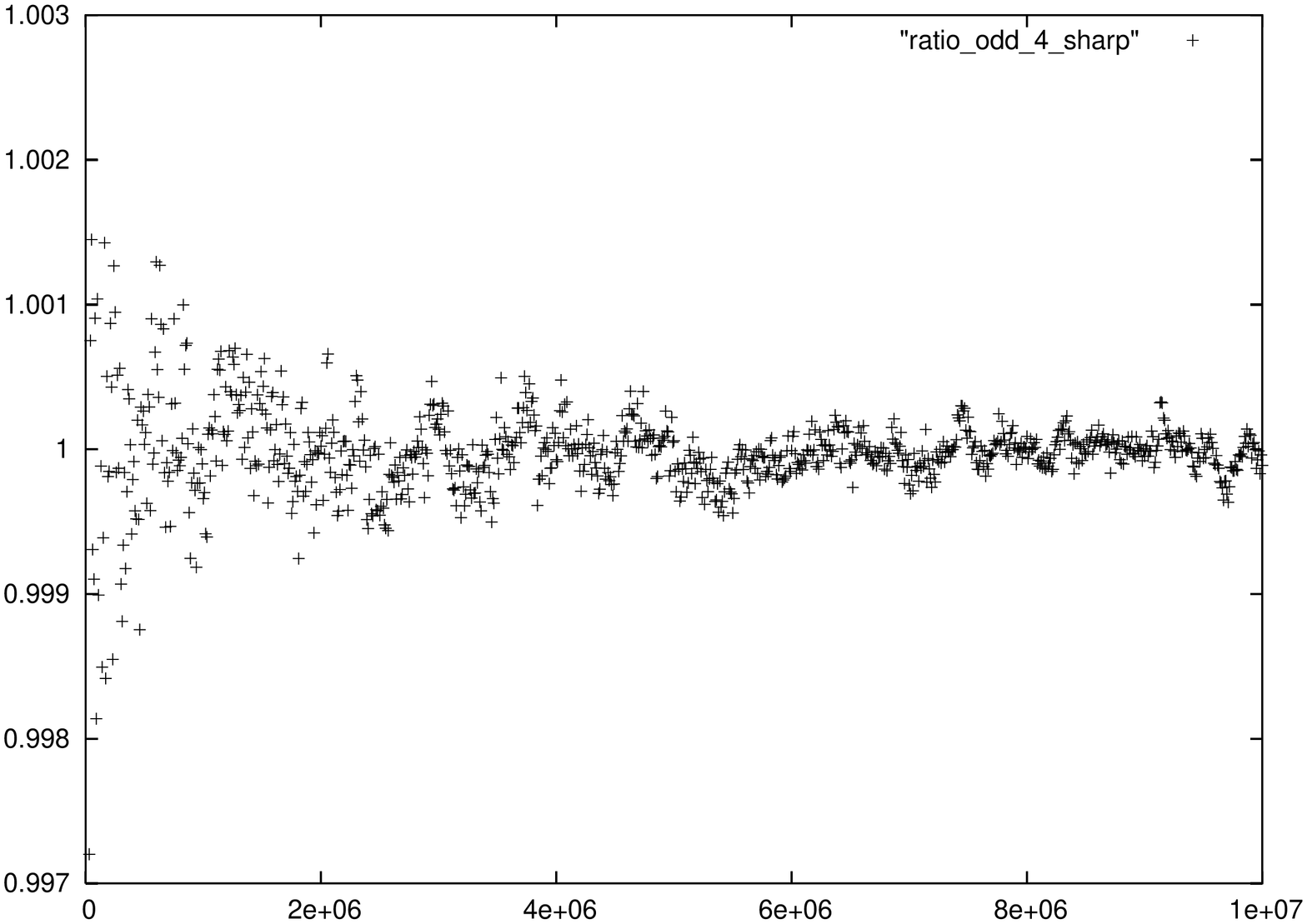,width=3in,angle=-90}
    }
    \centerline{
            \psfig{figure=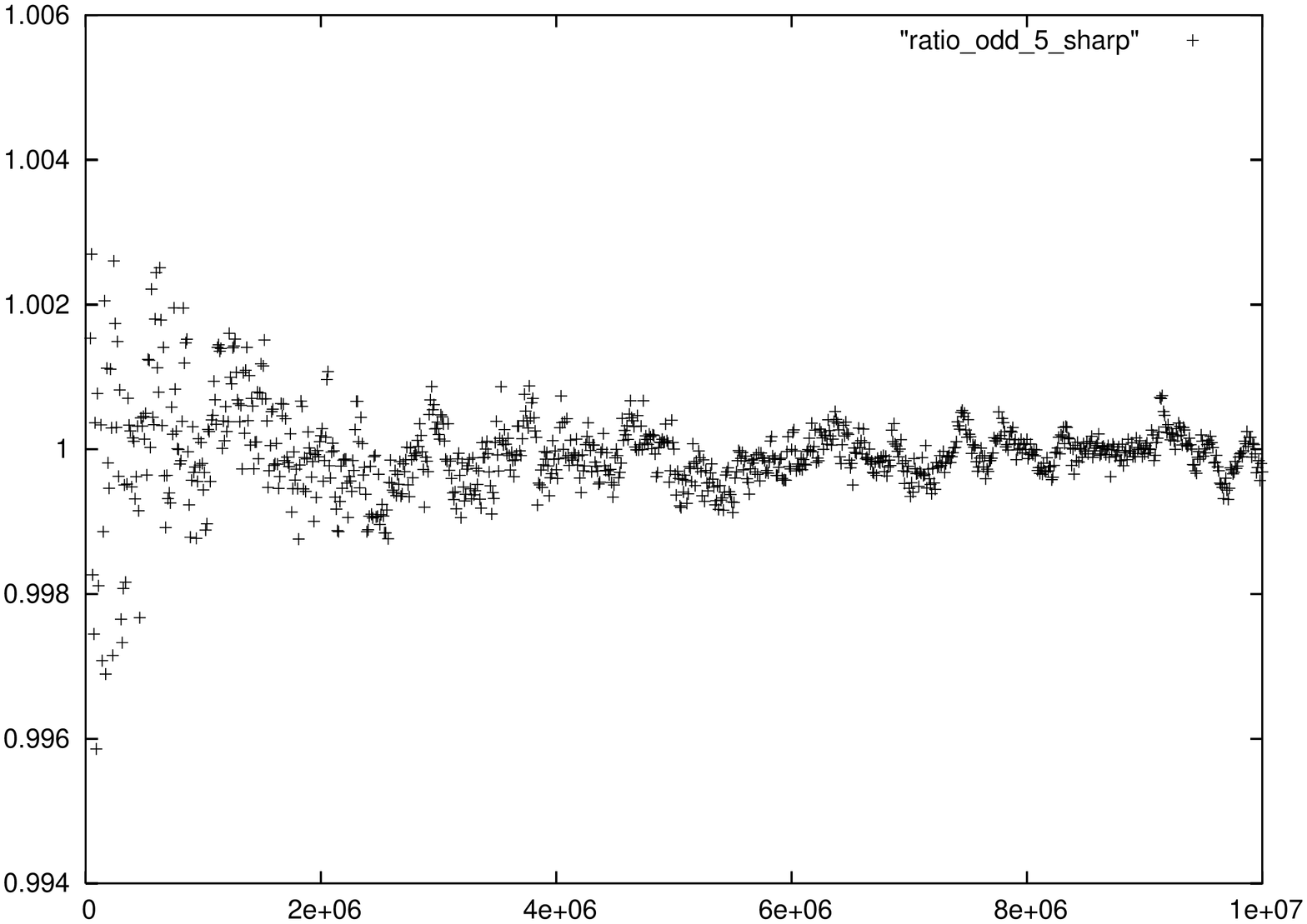,width=3in,angle=-90}
            \psfig{figure=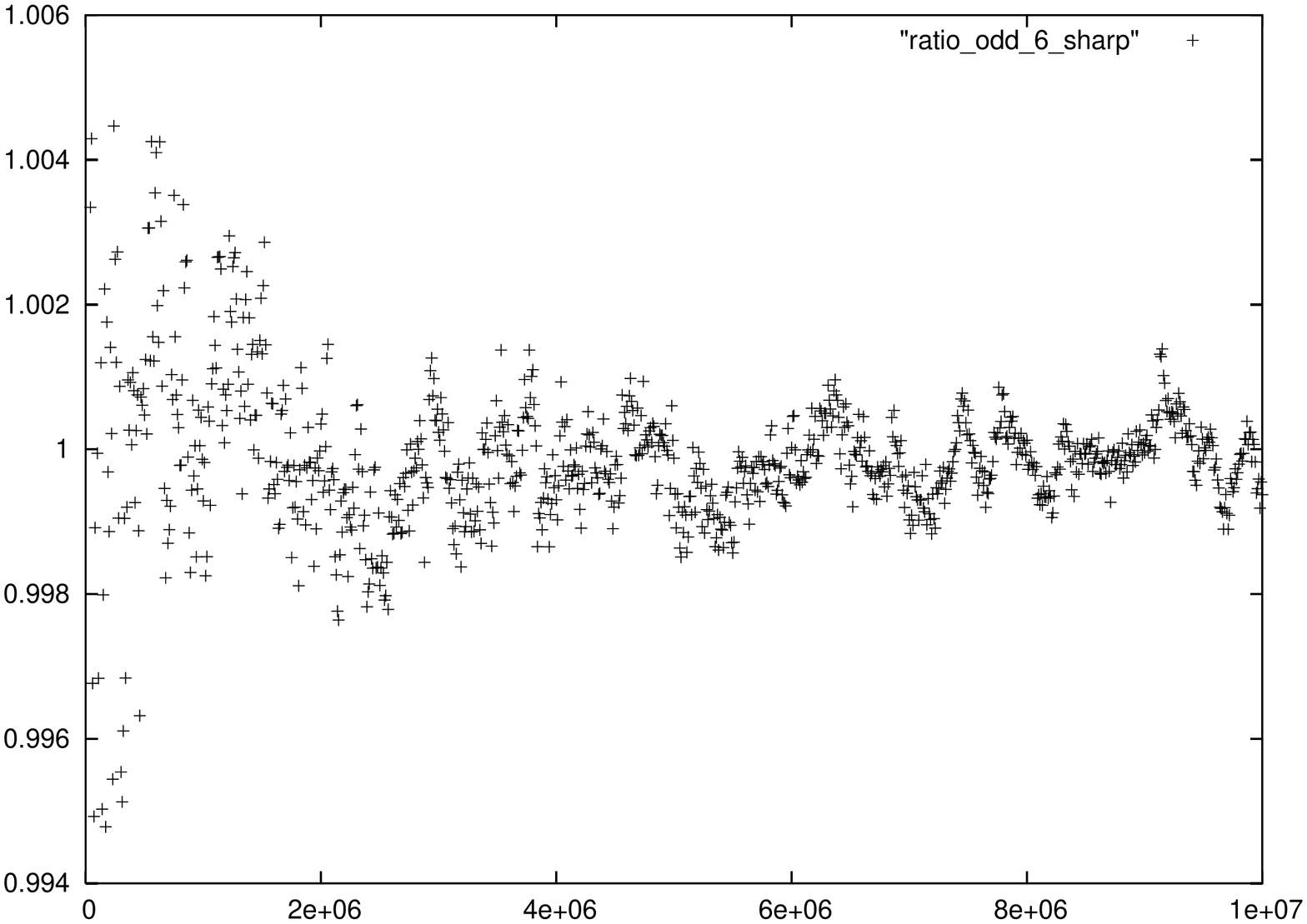,width=3in,angle=-90}
    }
    \centerline{
            \psfig{figure=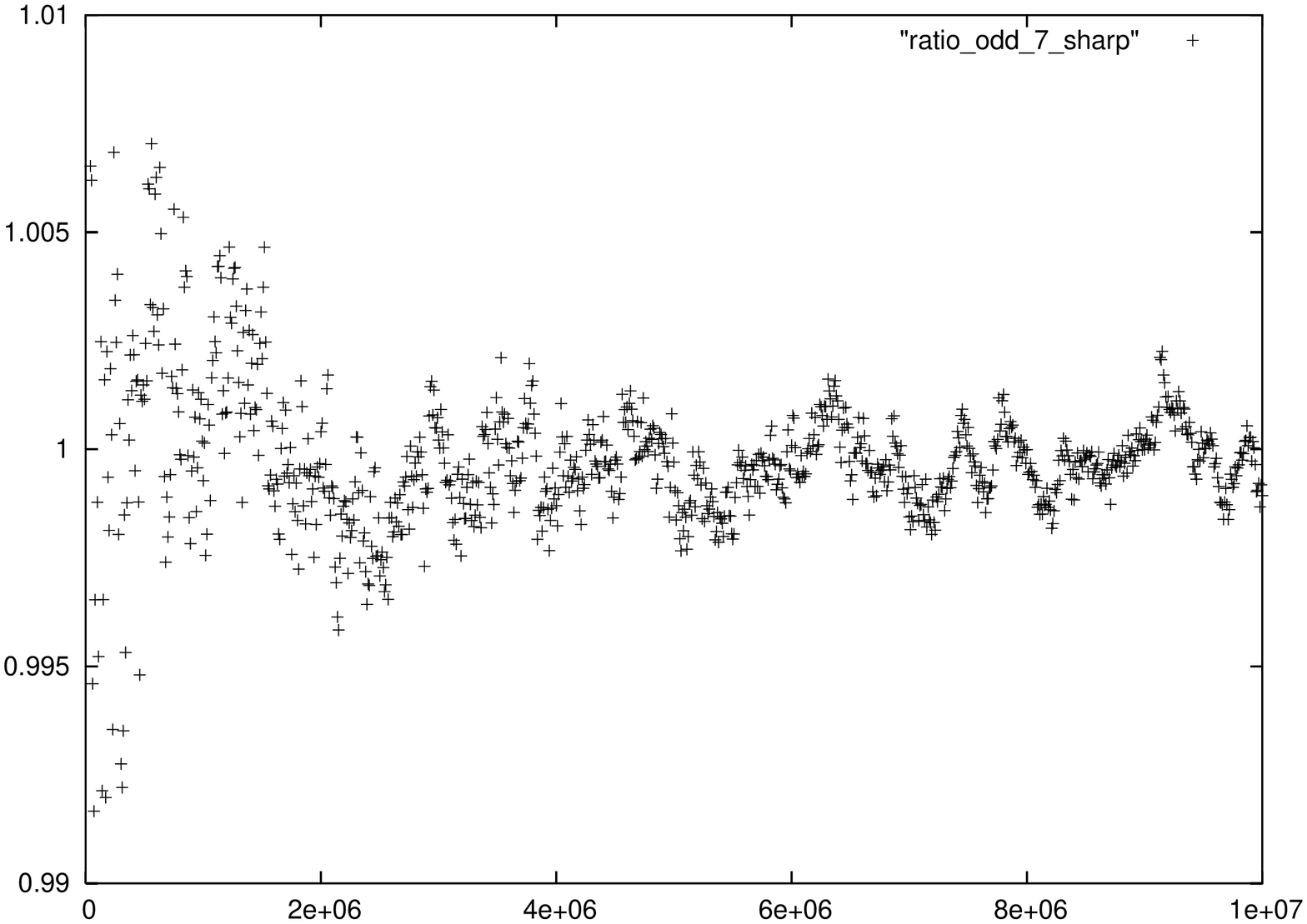,width=3in,angle=-90}
            \psfig{figure=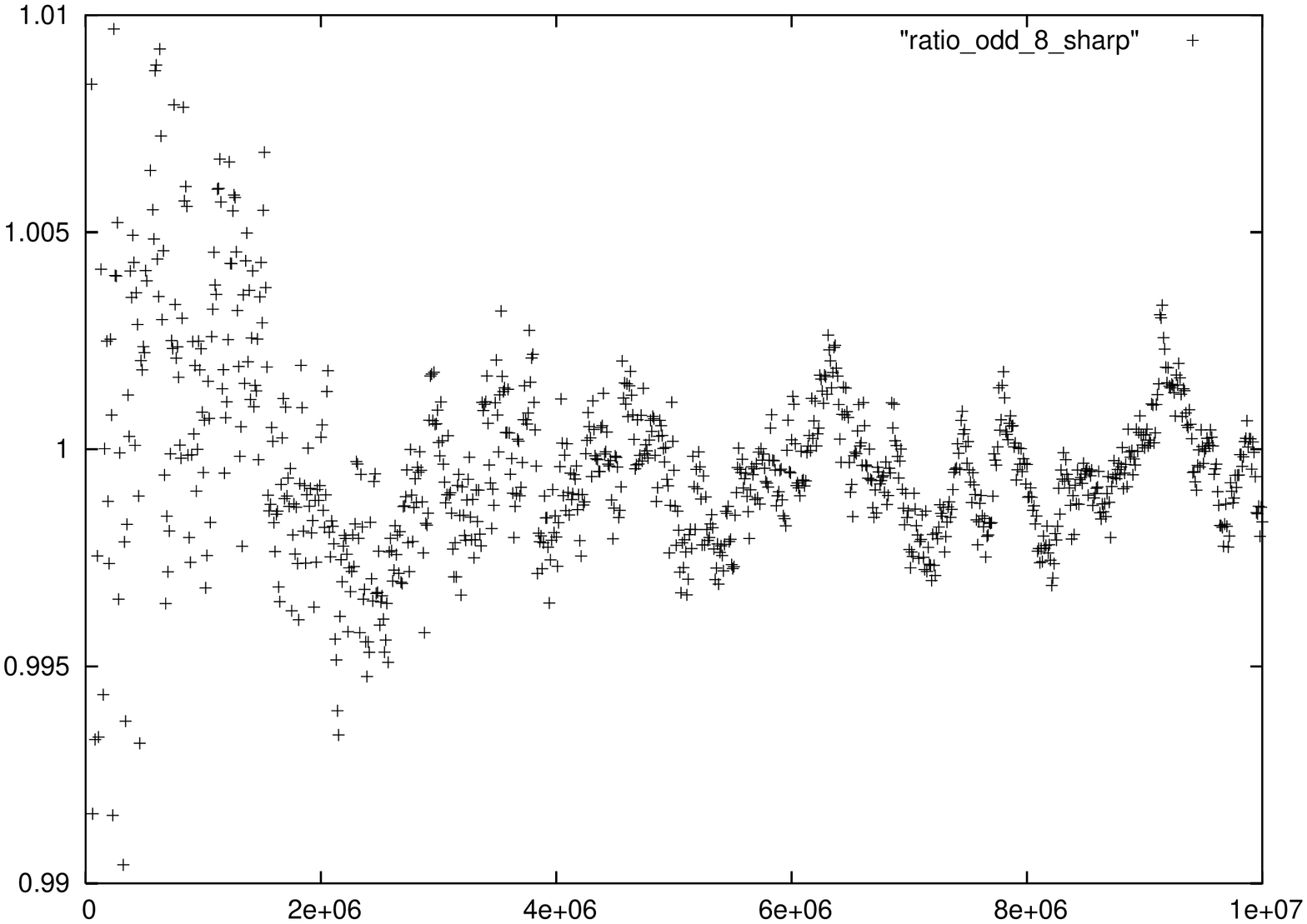,width=3in,angle=-90}
    }
    \caption
    {Horizontal axis in each graph is $X$. These graphs depict the first eight
     moments, sharp cutoff, of $L(1/2,\chi_d)$, $-X\leq d<0$ divided by the conjectured value,
     sampled at $X=10000,20000,\ldots,10^7$.
     One sees the graphs fluctuating above and below one. Notice that the vertical
     scale varies from graph to graph}.
    \label{fig:conj vs reality, odd chi}
\end{figure}

\begin{figure}[h!tb]
    \centerline{
            \psfig{figure=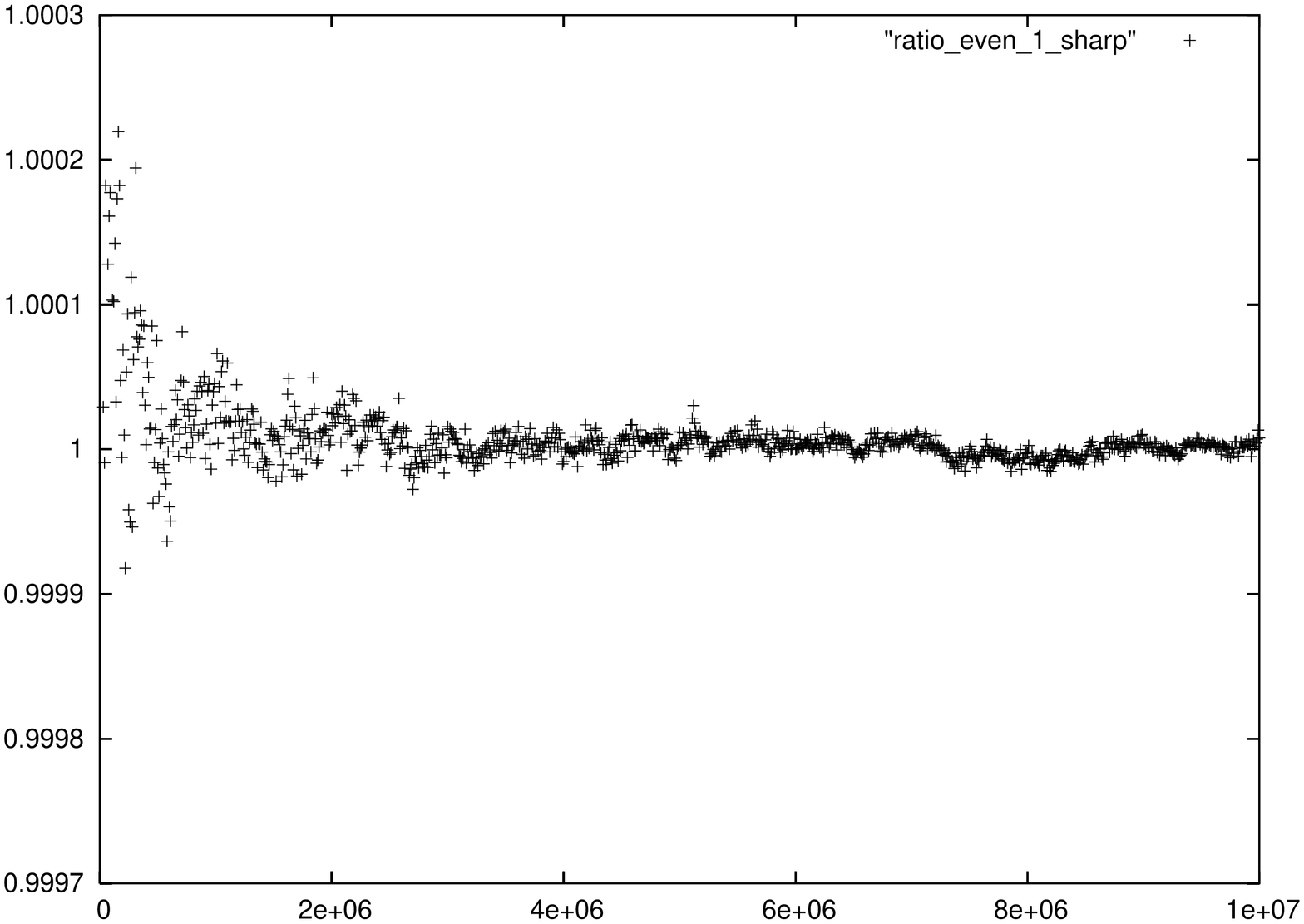,width=3in,angle=-90}
            \psfig{figure=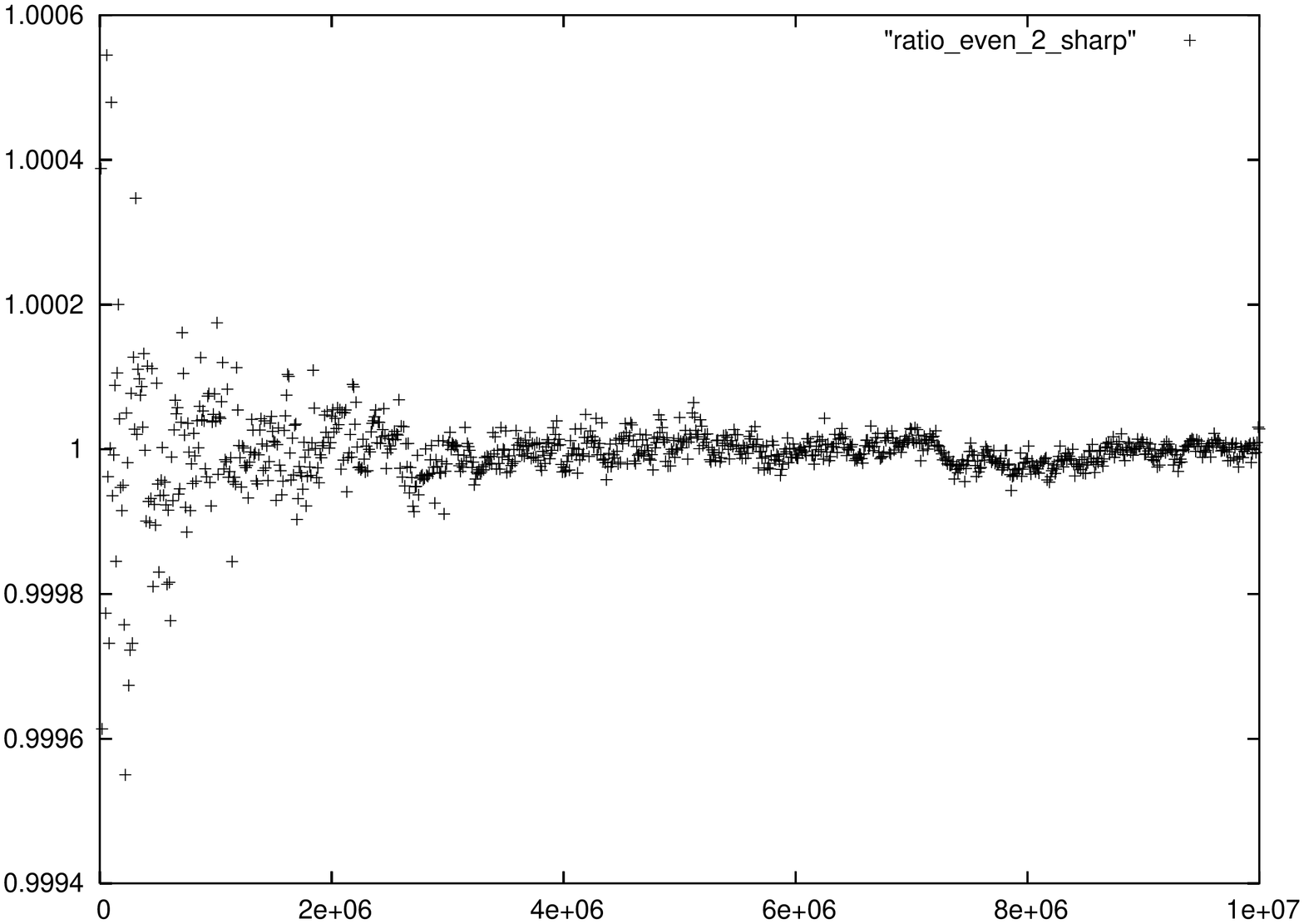,width=3in,angle=-90}
    }
    \centerline{
            \psfig{figure=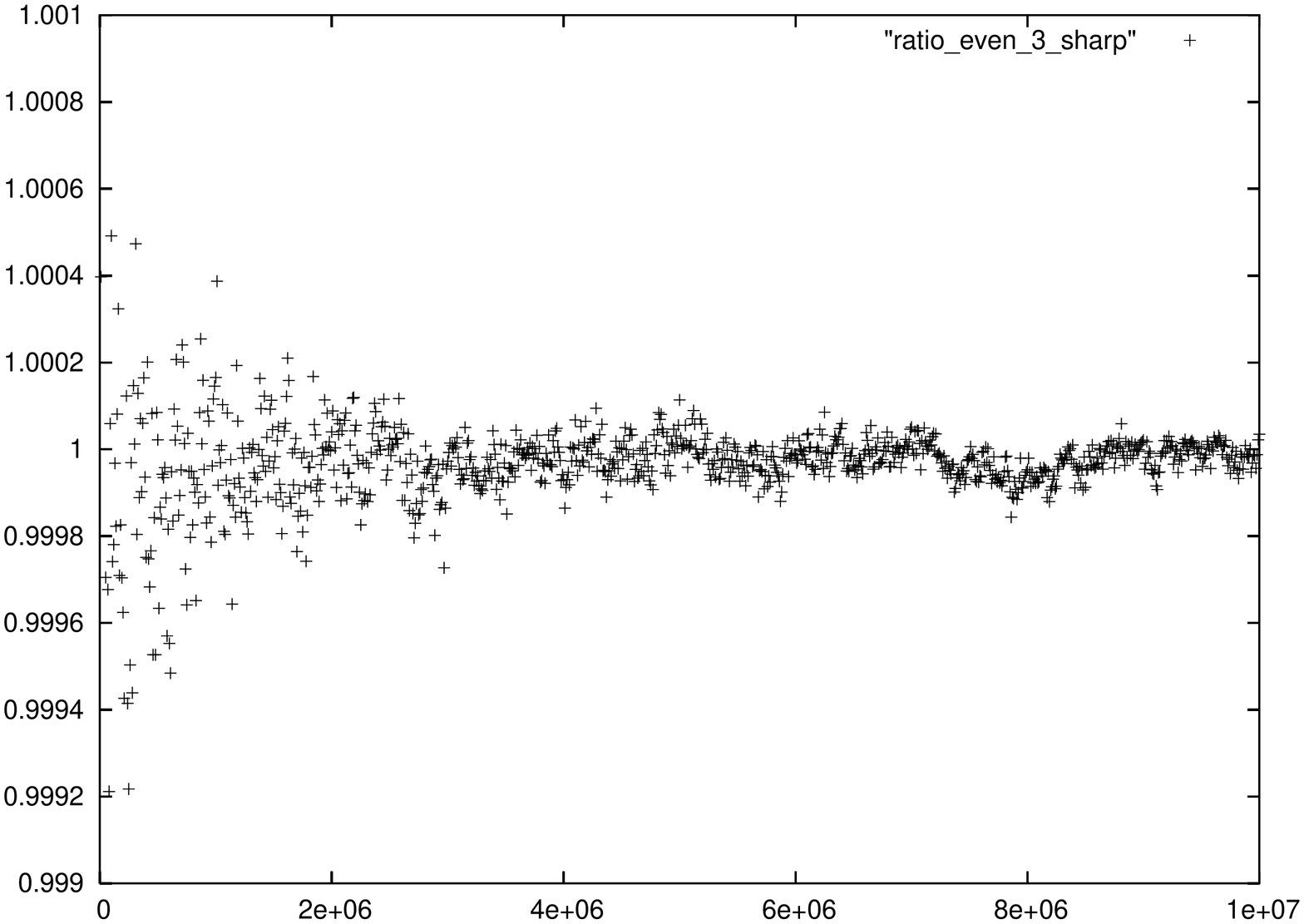,width=3in,angle=-90}
            \psfig{figure=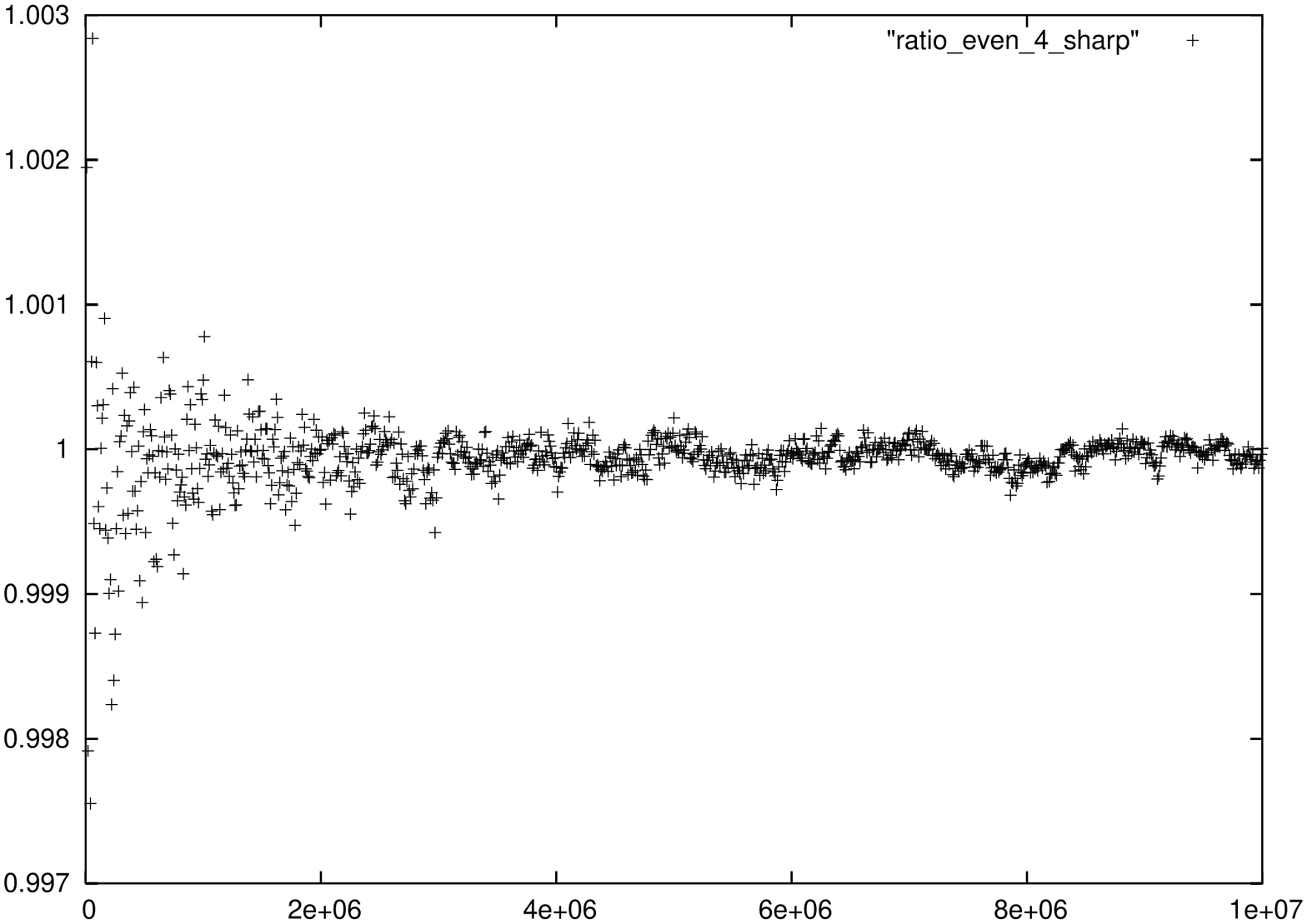,width=3in,angle=-90}
    }
    \centerline{
            \psfig{figure=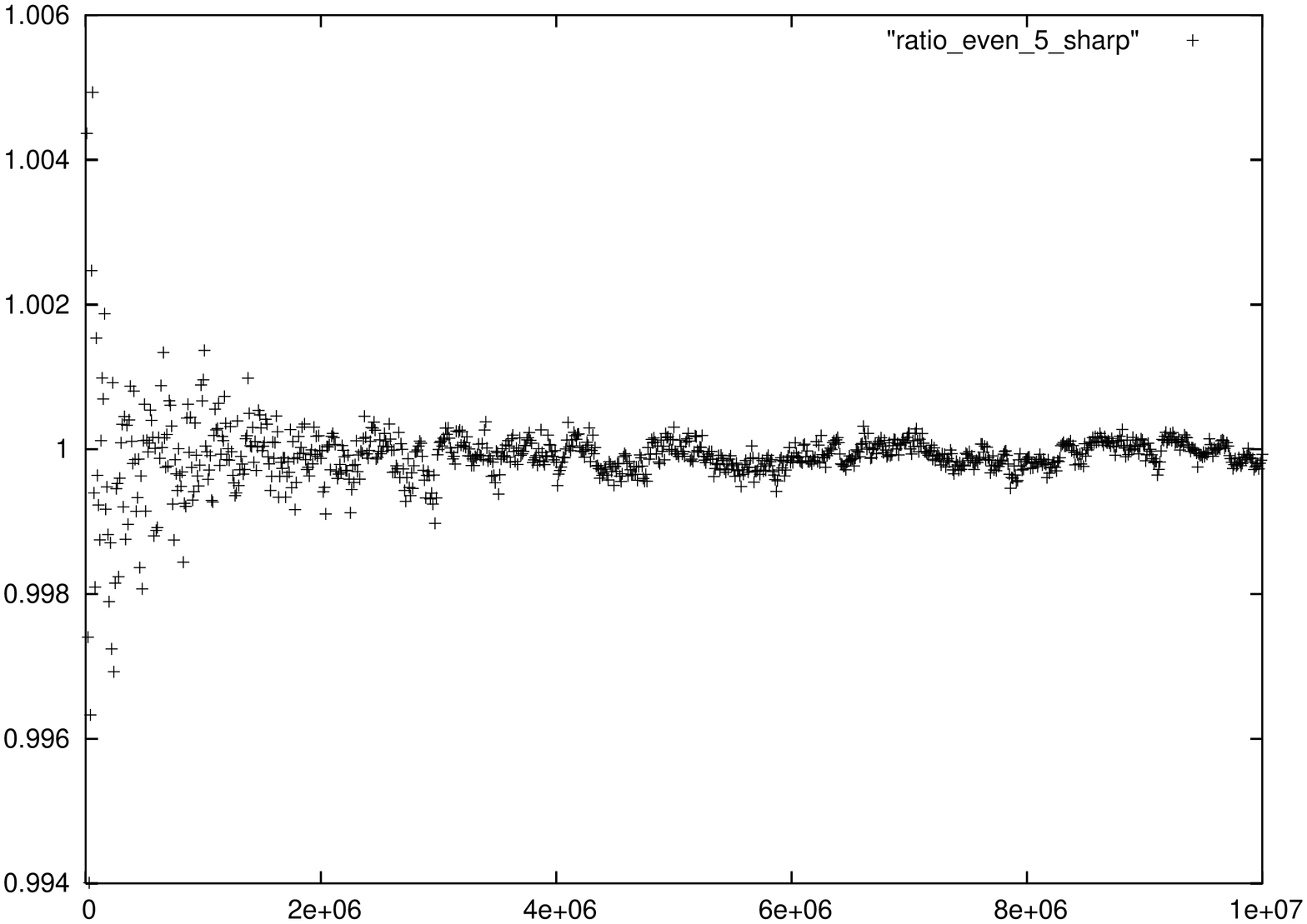,width=3in,angle=-90}
            \psfig{figure=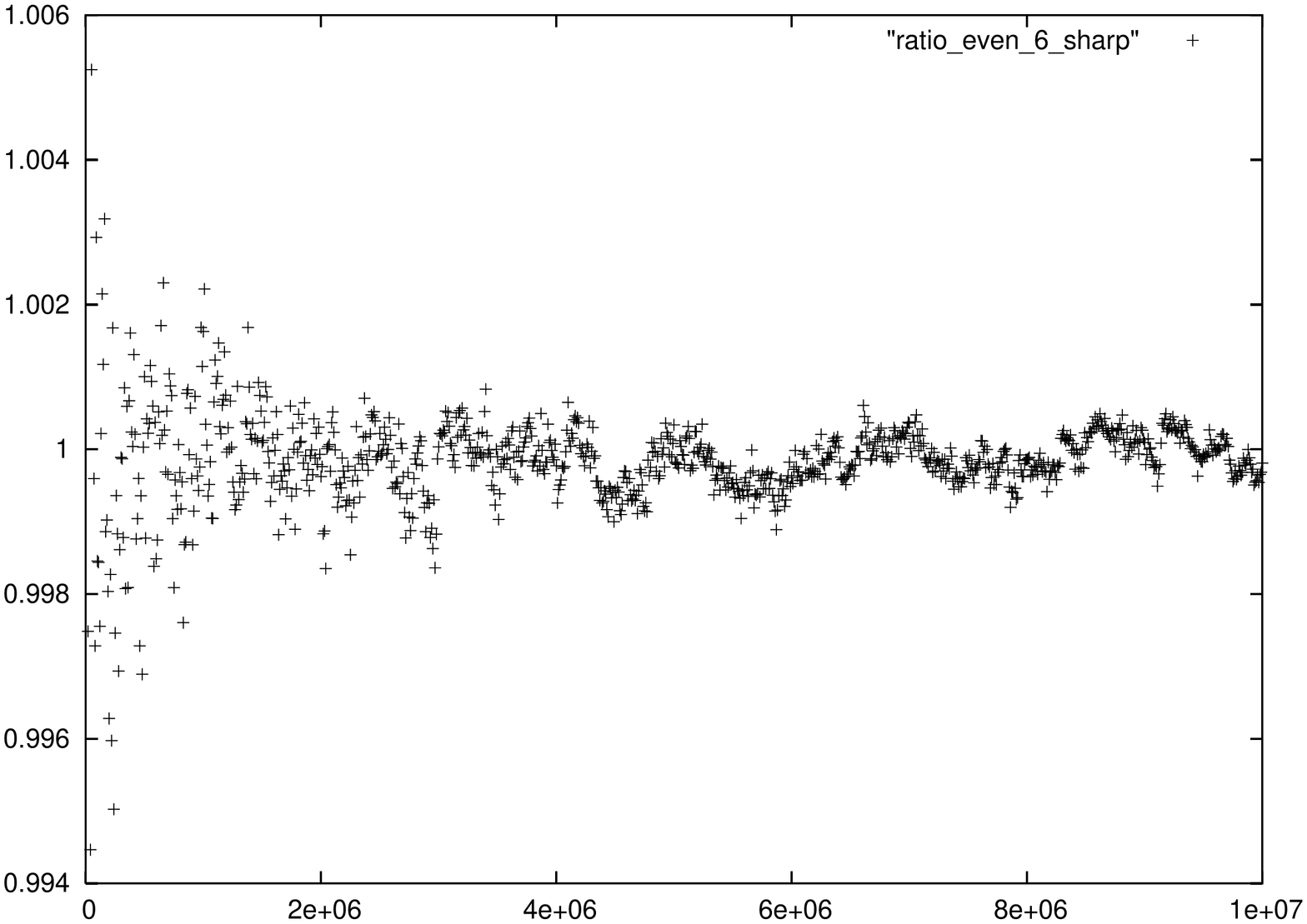,width=3in,angle=-90}
    }
    \centerline{
            \psfig{figure=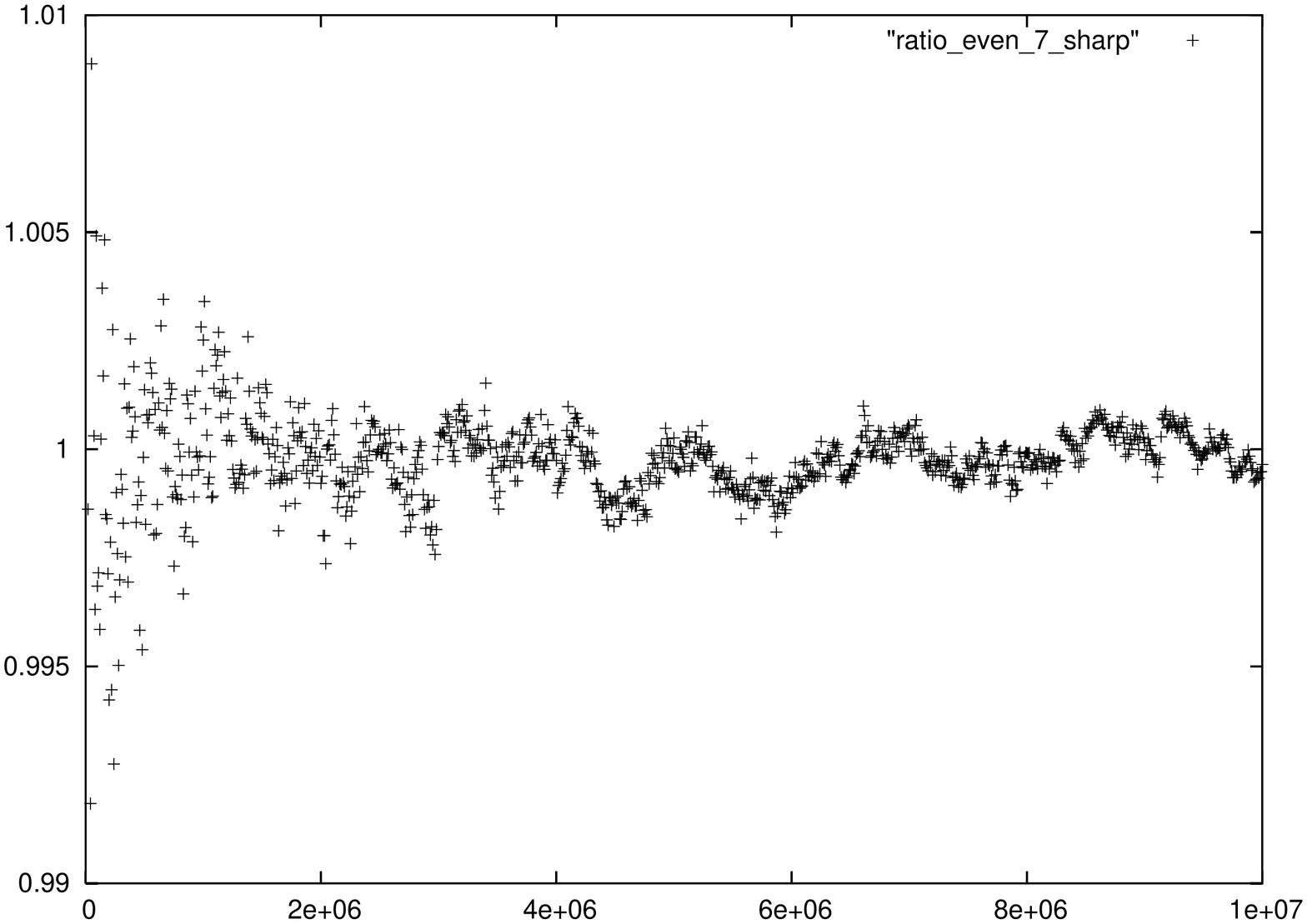,width=3in,angle=-90}
            \psfig{figure=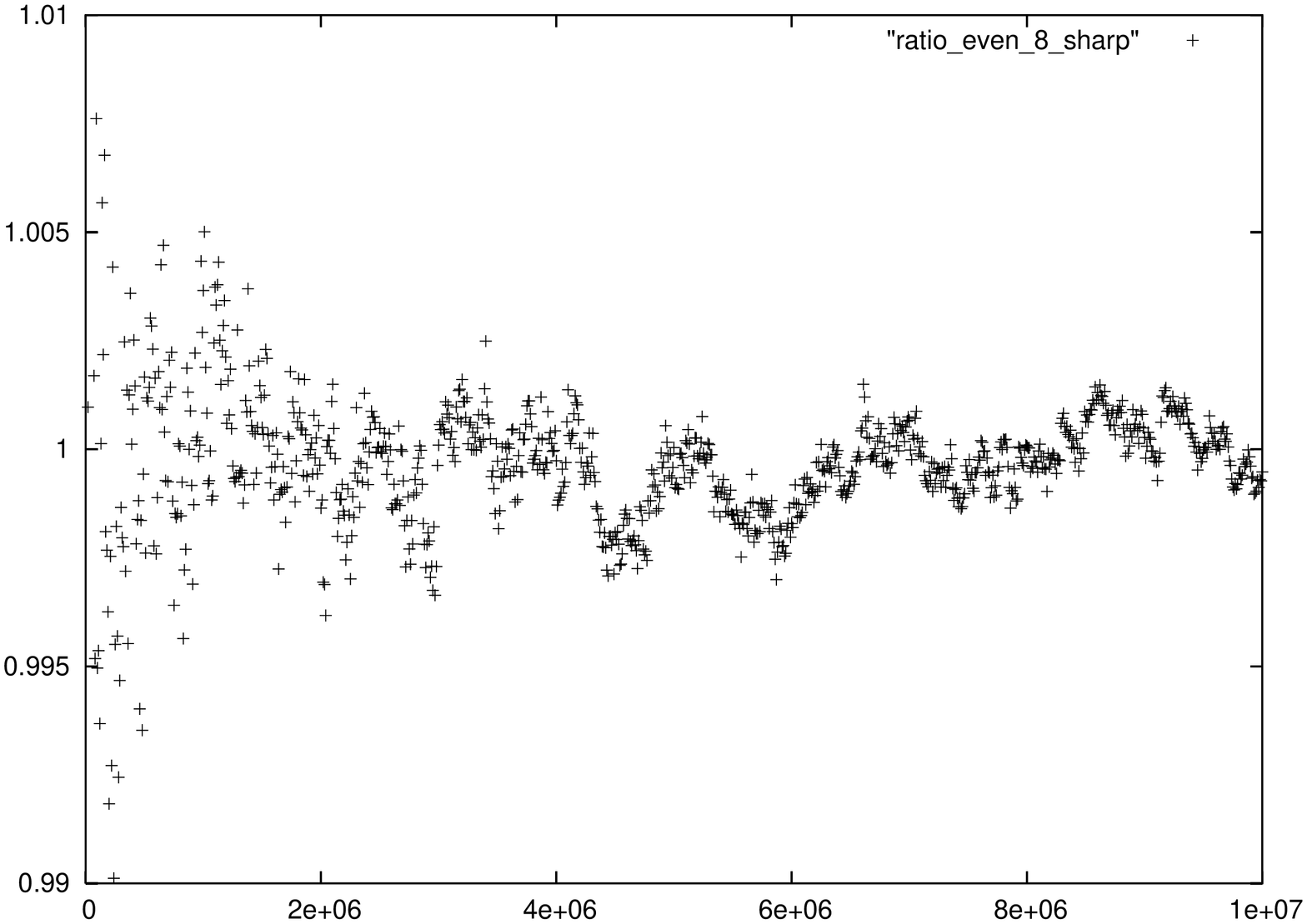,width=3in,angle=-90}
    }
    \caption
    {Same as the previous fingure, but for $0 < d \leq X$.}
    \label{fig:conj vs reality, even chi}
\end{figure}

\begin{table}[h!tb]
\centerline{\small
\begin{tabular}{|c|c|c|c|}
\hline
$k$ & reality (\ref{eqn:lhchidn}) & conjecture (\ref{eqn:Qkn})& ratio \cr \hline
1 &1144563.5        &1144535.5        &1.000024 \cr 2 &9252479.6
&9252229.9        &1.000027 \cr 3 &109917867.0      &109917367.9
&1.0000045 \cr 4 &1622521963.4     &1622508843.4     &1.0000081
\cr 5 &27321430060.     &27320230686.     &1.000043 \cr 6
&501621762060.6   &501542204848.7   &1.000159 \cr 7
&9787833470714.1  &9783848274459.6  &1.000407 \cr 8
&199831160877919 &199664775232854   &1.000833 \cr  \hline
\end{tabular}
}
\caption{ Smoothed moment of $L(\frac12 ,\chi_d)$ versus
Conjecture~\ref{thm:Lhalfchidconjecture}, for fundamental discriminants $0 < d < 1000000$,
and $k=1,\ldots,8$.
}\label{tab:Lhalfchip}\end{table}


\subsection{  Orthogonal: twists of a $GL(2)$ $L$-function
}\label{sec:orthogonalnumerics}

Let
\begin{equation}
    L_{11}(s)=\sum_{n=1}^\infty \frac{a_n}{n^{\frac12 +s}}
\end{equation}
be the $L$-function of conductor 11 of the elliptic curve
\begin{equation}
    y^2+y=x^3-x^2.
\end{equation}
The coefficients $a_n$ are obtained from the cusp form of weight
two and level 11 given by
\begin{equation}
    \sum_{n=1}^\infty a_n q^n = q \prod_{n=1}^\infty (1-q^n)^2
(1-q^{11n})^2.
\end{equation}
Expanding the right side using Euler's pentagonal theorem provides
an efficient means to compute the $a_n$'s.

$L_{11}(s)$ satisfies an even functional equation (i.e.
$\varepsilon=+1$):
\begin{equation}
   \left(\frac{11^{\frac12 }}{2\pi}\right)^s\Gamma(s+\tfrac12 )L_{11}(s) =
   \left(\frac{11^{\frac12 }}{2\pi}\right)^{1-s}\Gamma(\tfrac32-s)L_{11}(1-s),
\end{equation}
and may be written as a product over primes
\begin{equation}
    L_{11}(s)=
    \frac{1}{1-11^{-s-\frac12 }}
    \prod_{p \neq 11}
    \frac{1}{1-a_pp^{-s-\frac12 }+p^{-2s}}.
\end{equation}

Consider now quadratic twists of $L_{11}(s)$,
\begin{equation}
    L_{11}(s,\chi_d)=\sum_{n=1}^\infty \frac{a_n}{n^{\frac12 +s}}
\chi_d(n).
\end{equation}
with $(d,11)=1$.
$L_{11}(s,\chi_d)$ satisfies the functional equation
\begin{equation}
    L_{11}(s,\chi_d) = \chi_d(-11)
                       \frac{\Gamma(\frac32-s)}{\Gamma(s+\frac12 )}
                       \left(\frac{2\pi}{11^{\frac12 }}\right)^{2s-1}
                       |d|^{2(\frac12 -s)}
                       L_{11}(1-s,\chi_d).
\end{equation}
We wish to look at moments of $L_{11}(\frac12 ,\chi_d)$ but only for those
$L(s,\chi_d)$ that have an even functional equation, i.e. $\chi_d(-11)=1$.
We further only look at
$d<0$ since in that case a theorem of Kohnen and Zagier \cite{KZ} enables us
to easily gather numerical data for 
$L_{11}(\frac12 ,\chi_d)$ with which to check our
conjecture.

When $d<0$, $\chi_d(-1)=-1$, hence, in order to have even
functional equation, we require $\chi_d(11)=-1$, i.e.
$d=2,6,7,8,10 \mod 11$.
Conjectured formula~(\ref{thm:Lhalfchidsum}) combined with Lemma~\ref{thm:concisesumsymplectic}
gives an estimate for the sum over fundamental discriminants
\begin{equation}
     \sumstar_{ \ontop{-D<d<0}{d =  2,6,7,8,10 \mod 11} }
      L_{11}(\tfrac12 ,\chi_d)^k
      = \sumstar_{\ontop{-D<d<0}{d =  2,6,7,8,10 \mod 11}}
        \Upsilon_k\left(\log{|d|}\right)
        +O(D^{\frac12 + \varepsilon})
\label{eqn:L11sum}
\end{equation}
where, as in Section~\ref{sec:symplectic},  $\Upsilon_k$ is the polynomial of degree $k(k-1)/2$ given by
the $k$-fold residue
\begin{equation}
\Upsilon_k(x)=\frac{(-1)^{k(k-1)/2}2^k}{k!} \frac{1}{(2\pi i)^{k}}
\oint \cdots \oint \frac{R_{11}(z_1,
\dots,z_{k})\Delta(z_1^2,\dots,z_{k}^2)^2} {\displaystyle
\prod_{j=1}^{k} z_j^{2k-1}} e^{x \sum_{j=1}^{k}z_j}\,dz_1\dots
dz_{k} ,
\end{equation}
where
\begin{equation}
R_{11}(z_1,\dots,z_k)=A_k(z_1,\dots,z_k)
\prod_{j=1}^k
\left(
    \frac{\Gamma(1+z_j)}
         {\Gamma(1-z_j)}
    \left(\frac{11}{4\pi^2} \right)^{z_j}
\right)^{\frac 12} \prod_{1\le i < j\le k}\zeta(1+z_i+z_j),
\end{equation}
and $A_k$ is the Euler product which is absolutely convergent for
$\sum_{j=1}^k |z_j|<\frac12 $,
\begin{equation}
A_k(z_1,\dots,z_k) =
    \prod_p R_{11,p}(z_1,\ldots,z_k)
    \prod_{1\le i < j \le k}
    \left(1-\frac{1}{p^{1+z_i+z_j}}\right)
\end{equation}
with, for $p \neq 11$,
\begin{equation}
R_{11,p} =
           \left(1+\frac 1 p\right)^{-1}
           \left(\frac 1 p +\frac{1}{2}
              \left(
                 \prod_{j=1}^k
                 \frac{1}{1-a_pp^{-1-z_j}+p^{-1-2z_j}}
               + \prod_{j=1}^k
                 \frac{1}{1+a_pp^{-1-z_j}+p^{-1-2z_j}}
              \right)
           \right)
\end{equation}
and
\begin{equation}
    R_{11,11} =
    \prod_{j=1}^k  \frac{1}{1+11^{-1-z_j}}.
\end{equation}

Numerically, it is more challenging to compute the polynomials
$\Upsilon_k$. First, using
\begin{equation}
    \prod_{1\le i < j\le k}\zeta(1+z_i+z_j)
\end{equation}
to estimate the sum over primes of (\ref{eqn:RNp1overp}) makes a poor
approximation and one would do better to use the Rankin-Selberg
convolution $L$-function of $L_{11}(s)$ with itself. However, it
is simpler to work with $\zeta$, and we thus computed the first 4
moment polynomials of $L_{11}(\frac12 ,\chi_d)$ but to low accuracy.
The coefficients of these polynomials are given to 2-5 decimal
place accuracy in Table~\ref{tab:fr}.

In Table~\ref{tab:L11moments} we compare moments computed numerically with moments
estimated by our conjecture. The two agree to within the accuracy we
have for the moment polynomial coefficients. We believe that if one
were to compute the coefficients to higher accuracy, one would see
an even better agreement with the data.

While one can compute $L_{11}(\frac12 ,\chi_d)$ using standard
techniques (see \cite{Coh}), in our case we can exploit a theorem
of Kohnen and Zagier \cite{KZ} which relates $L_{11}(\frac12 ,\chi_d)$,
for fundamental discriminants $d<0$, $d\equiv 2,6,7,8,10 \mod 11$, to
the coefficients $c_{11}(|d|)$ of a weight $3/2$ modular form
\begin{equation}
    L_{11}(\tfrac12 ,\chi_d) = \kappa_{11} c_{11}(|d|)^2/\sqrt{d}
\end{equation}
where $\kappa_{11}$ is a constant.
The weight $3/2$ form in question was determined
by Rodriguez-Villegas (private communication)
\begin{eqnarray}
   \sum_{n=1}^\infty c_{11}(n) q^n&  = & (\theta_1(q)-\theta_2(q))/2 \nonumber \\&  = &\mbox{} -q^3+q^4+q^{11}+q^{12}-q^{15}-2q^{16}-q^{20}\ldots
\end{eqnarray}
where
\begin{equation}
    \theta_1(q) = \sum_{  \ontop{(x,y,z) \in {\mathbb Z}^3}{x\equiv y \bmod 2}}
    q^{x^2+11y^2+11z^2}
    =1+2q^4+2q^{11}+4q^{12}+4q^{15}+2q^{16}+4q^{20}\ldots
\end{equation}
and
\begin{equation}
    \theta_2(q) =  \sum_{  \ontop{(x,y,z) \in {\mathbb Z}^3}{\ontop{x\equiv y \bmod 3}{y\equiv z \bmod 2}}}
    q^{(x^2+11y^2+33z^2)/3}
    =1+2q^3+2q^{12}+6q^{15}+6q^{16}+6q^{20}\ldots.
\end{equation}
This was used to compute the $c_{11}(|d|)$'s for $d<85,000,000$.

Evaluating the left side of (\ref{eqn:L11sum}) in a more traditional manner
for $d=-3$, and comparing with the right side, we determined
\begin{equation}
    \kappa_{11} = 2.917633233876991.
\end{equation}

\begin{table}[h!tb]
\centerline{\small
\begin{tabular}{|c|c|c|c|c|}
\hline
r& $f_r(1)$ &$f_r(2)$ &$f_r(3)$ &$f_r(4)$ \cr \hline
0 & 1.2353 & .3834 & .00804 & .0000058 \cr
1 & & 1.850 & .209 & .000444 \cr
2 & & & 1.57 & .0132 \cr
3 & & & 2.85 & .1919 \cr
4 & & & & 1.381 \cr
5 & & & & 4.41 \cr
6 & & & & 4.3 \cr
\hline
\end{tabular}
}
\caption{
Coefficients of $\Upsilon_k(x) = f_0(k) x^{k(k-1)/2} + f_1(k)
x^{k(k-1)/2-1} + \ldots$, for k=1,2,3,4.
}\label{tab:fr}\end{table}

\begin{table}[h!tb]
\centerline{\small
\begin{tabular}{|c|c|c|c|}
\hline
$k$ & left side (\ref{eqn:L11sum}) & right side (\ref{eqn:L11sum})& ratio \cr \hline
1 &14628043.5       &14628305.       &0.99998 \cr 2 &100242348.8
&100263216.      &0.9998 \cr 3 &1584067116.8     &1587623419.
&0.998 \cr 4 &41674900434.9    &41989559937.    &0.993 \cr \hline
\end{tabular}
}
\caption{ Moments of $L_{11}(\frac12 ,\chi_d)$ versus
their conjectured values, for fundamental discriminants
$-85,000,000 < d < 0$, $d=2,6,7,8,10 \mod 11$, and $k=1,\ldots,4$.
The data agree with our conjectures to the accuracy to which we
have computed the moment polynomials $\Upsilon_k$.
}\label{tab:L11moments}\end{table}




\end{document}